\documentclass{article}	

\usepackage[utf8]{inputenc}				

\usepackage{amsmath,amsfonts}
\usepackage{lmodern} 
\usepackage{microtype}
\usepackage{xargs,hyperref}  
\usepackage[numbers,square,sort&compress]{natbib}
\usepackage{color}
\usepackage{todonotes}
\usepackage{algorithmic}
\usepackage{algorithm2e}
\usepackage[normalem]{ulem}
\usepackage{subcaption}
\usepackage[letterpaper,top=2cm,bottom=2cm,left=3cm,right=3cm,marginparwidth=1.75cm]{geometry}

\newtheorem{remark}{Remark}

\def \x {{\bf x}}

\def\RR{\mathbb{R}}

\def \pmatrix{ \left( \begin{array} }
\def \endpmatrix{ \end{array} \right) }
\def \u {{\bf u}}
\def \v {{\bf v}}

\def \s {{\bf s}}
\def \f {{\bf f}}
\def \g {{\bf g}}
\def \F {{\bf F}}
\def \G {{\bf G}}

\def \R {{\mathbb R}}
\DeclareMathOperator{\argmax}{argmax}
\DeclareMathOperator{\rank}{rank}

\title{Adaptive POD-DEIM correction for Turing pattern approximation in reaction-diffusion PDE systems}

\author{Alessandro Alla$^{1}$, Angela Monti$^{2}$, Ivonne Sgura$^{2}$\\
\small $^{1}$ Università Ca' Foscari Venezia, Dipartimento di Scienze Molecolari e Nanosistemi, Venezia, Italy,\\
\small e-mail: alessandro.alla@unive.it\\
\small $^{2}$ Università del Salento, Dipartimento di Matematica e Fisica ‘‘E. De Giorgi’’, Lecce, Italy,\\ 
\small e-mail: \{angela.monti, ivonne.sgura\}@unisalento.it}

\begin{document}
\maketitle
	\begin{abstract}
We investigate a suitable application of Model Order Reduction (MOR) techniques for the numerical approximation of Turing patterns, that are stationary solutions of reaction-diffusion PDE (RD-PDE) systems. We show that solutions of surrogate models built by classical Proper Orthogonal Decomposition (POD) exhibit an unstable error behaviour over the dimension of the reduced space. To overcome this drawback, first of all, we propose a POD-DEIM technique with a correction term that includes missing information in the reduced models. To improve the computational efficiency, we propose an adaptive version of this algorithm in time that accounts for the peculiar dynamics of the RD-PDE in presence of Turing instability. We show the effectiveness of the proposed methods in terms of accuracy and computational cost for a selection of RD systems, i.e. FitzHugh-Nagumo, Schnackenberg and the morphochemical DIB models, with increasing degree of nonlinearity and more structured patterns. 
\end{abstract}

{\noindent \textbf{Keywords:} Reaction-Diffusion PDEs, Turing patterns, Model Order Reduction, Proper Orthogonal Decomposition, Adaptivity, Discrete Empirical Interpolation Method. }\\
{\noindent \textbf{MSC:} 65M06, 35K57, 65F99, 65M22 }

\section{Introduction} 

In different fields of application, mathematical models can be expressed in the form of Reaction-Diffusion Partial Differential Equations (RD-PDEs). These evolutionary PDE systems are defined in one or more space dimensions, on stationary or evolving domains and surfaces. The nonlinear kinetics of the coupled reaction terms account for physical, biological and other kind of phenomena. In the wide literature on these topics, we cite only a selection of works like \cite{Murray03book,Maini01book,Madzva03,Barreira11} for applications in bio-mathematics, \cite{Malchow08book,Sherratt12,Gilad04} in ecology, \cite{Gerisch06,Garzon11,Lefevre10brain} in biomedicine, and other applications like \cite{Chaplain01,Sherratt01} for tumor growth, \cite{Painter19} for chemotaxis and \cite{DIB12,DIB13,DIB15,DIB17sphere} for metal electrodeposition.

In this paper, we are interested in RD-PDE systems where the interplay between diffusion and reaction is responsible of the so-called \emph{Turing pattern formation}.
For simplicity, we focus on RD-PDE systems of two equations that are defined on a stationary 2D spatial domain given by
\begin{align}
\label{RDPDE}
\begin{aligned}
u_t &= d_u \Delta u + f(u,v) , \quad (x,y) \in \Omega \subset 
\mathbb{R}^2 , \quad t \in (0,T], \\
v_t &= d_v \Delta v + g(u,v), \\
(\mathbf{n} \nabla u)_{| \partial \Omega} & = (\mathbf{n} \nabla v)_{| \partial 
\Omega}  = 0, \\
u(x,y,0)& = u_0(x,y) , \quad v(x,y,0) = v_0(x,y)
\end{aligned}
\end{align}
with $\Omega = [0,L_x] \times [0,L_y]$, nonlinear reaction terms $f$ and $g$ and zero Neumann boundary conditions. The diffusion coefficients are $d_u, d_v>0$ and there exists at least a spatially homogeneous equilibrium
$P_e=(u_e,v_e)$, such that $f(u_e,v_e)=0$ and $g(u_e,v_e)=0$.
It is well known (see e.g. \cite{Murray03book}) that the coupling between diffusion and nonlinear reaction terms can lead to the so-called \emph{diffusion-driven} or \emph{Turing} instability, when $P_e$ is \emph{stable} in absence of diffusion $(d_u=d_v=0)$ and becomes
\emph{unstable} when the diffusion is present, giving rise to inhomogeneous spatial patterns as stationary solutions of \eqref{RDPDE}, called \emph{Turing patterns}. To study pattern formation, the initial condition of \eqref{RDPDE} is typically a small (random) spatially distributed perturbation of $P_e$ that destabilizes such that the dynamics can attain at the steady state different kind of Turing patterns, like labyrinths, spots, stripes, etc., in correspondence of different choices of the parameters involved in the model kinetics.

It is worth emphasizing that the above classical theory of spatial pattern formation is an asymptotic theory, concerned with the long term behaviour of perturbations, but more recently, in
\cite{Neubert97,Neubert02} the authors proved that also the transient dynamics is important for pattern formation. They introduced two quantitative indicators called
\emph{resilience} and \emph{reactivity}, to describe the
asymptotic and short-term time regimes,  respectively, and proved that 
the concept of \emph{reactivity} is a \emph{a necessary condition} for Turing instabilities. In \cite{DSS20}, it has been proposed a quantitative way to identify these two time regimes 
${\cal I}_1 = [0, \tau]$ and ${\cal I}_2 =(\tau, T]$ such that the value of $\tau$ can be obtained a-posteriori looking at the qualitative behaviour of the numerical solution of \eqref{RDPDE} (see Section \ref{sec:4}). In ${\cal I}_1$ the reactivity
holds: the solution departs oscillating from the spatially
homogeneous pattern due to the superimposed (small random)
perturbations and becomes unstable, in ${\cal I}_2$ the solution starts to stabilize towards the steady Turing pattern.

Only for some choices of the kinetics it is possible, by means of weakly nonlinear analysis, to identify analytically the spatial modes of the Turing patterns for special values of the parameters (see e.g. \cite{Gambino19FHN,Sammartino15DIB}). Usually, for a given parameter choice, a numerical approximation of \eqref{RDPDE} is required to simulate Turing patterns until the steady state and by reproducing the reactivity features during the initial transient dynamics.
%
Therefore, the numerical approximation of Turing pattern solutions is
challenging for the following reasons: (i) long time integration is
needed to attain the stationary pattern; (ii) a large domain $\Omega$ is needed to carefully \emph{see} its spatial structures; 
(iii) an accurate spatial discretization is required to identify the {pattern class} (labyrinths, spots, etc..); (iv) the time solver should capture also the reactivity regime. 


The main goal of this paper is the application of Model Order Reduction (MOR) techniques to reduce the computational costs of a standard numerical method for \eqref{RDPDE} while preserving the main features of the Turing dynamics as discussed above. Hence, we apply the Proper Orthogonal Decomposition (POD, see e.g. \cite{Sir87}) with an hyper reduction of the nonlinear terms by the Discrete Empirical Interpolation Method (DEIM, \cite{CS10,DG15}). 
We mention that model reduction for coupled PDEs has been investigated in \cite{RS07, RS08} for linear time invariant coupled systems and in \cite{BF15} a survey on the topic has been presented. The study of model reduction for reaction-diffusion PDEs with pattern formation has been first investigated in \cite{KUK15}, where the authors have applied POD and DEIM for the FitzHugh-Nagumo model. More recently, in \cite{KMUY21} it has been proposed a partitioned POD technique to deal with asymptotic patterns in a nonlinear cross-diffusion system in population ecology. 
A first study of MOR for Turing type dynamics based on POD and Dynamic Mode Decomposition (DMD, \cite{DMDbook}) applied to the Schnackenberg model (\cite{Madzva03}) and to the DIB kinetics (see e.g. \cite{DIB13}) for battery modeling, has been proposed recently in \cite{BMS21}.

Here, we show that a straightforward application of POD and DEIM in usual forms presents non monotone, often oscillating, or unstable approximations for increasing sizes of the reduced space. The main motivation of this paper is to devise a suitable strategy to overcome this bad behaviour. 


Stabilization of model reduction is an active research topic especially for fluid dynamics models. We refer to e.g. \cite{BMQR15} where the use of supremizers solutions stabilizes the reduced problem, to e.g. \cite{AF12} for linear models, to e.g. \cite{BBSK17, WWZXI17} for the use of machine learning techniques and to e.g. \cite{BMQR15} for the use of a minimal rotation of the projection subspace.

In this paper, we follow the stabilization approach proposed in \cite{XMRI18} where a ``correction term'' is added to the reduced model to provide the missing information, possible cause of the instability. This correction is built upon a surrogate model which dimension is close to the rank of the snapshot matrix, whose colums are the numerical solution of the full model at a selection of time steps. We will show that the new reduced problem results to be very accurate, even if it can be very expensive in the offline stage. To speed up the method, we employ DEIM to this corrected version.
The corrected POD-DEIM here proposed turns out to be one of the novelty of the current work. This leads to a stabilization of the DEIM technique that further improves the results in \cite{DG15}. In \cite{XMRI18, MLWI20}, the authors have approximated the correction terms building a quadratic model based on a least square optimization problem. Moreover they have focused on a single equation, whereas we deal with coupled problems.

Furthermore, to lower the computational cost required by the construction of the correction term in the offline stage and to better capture the main features of the dynamics in the reduced spaces, we have introduced an adaptive strategy for the corrected POD-DEIM that operates in the time intervals $\mathcal{I}_1$ (reactivity zone) and $\mathcal{I}_2$ (stationary zone) discussed above. This partition allows to compute SVD for smaller snapshot matrices and also the correction term is built upon a smaller reduced space. We mention that an adaptive approach for DEIM has been introduced in e.g. \cite{PBWB14} for steady parametrized PDEs.

To summarize the novelties in this work are: (1) the corrected POD-DEIM algorithm and (2) adaptive MOR driven by the time dynamics in the reactive and asymptotic regimes. This leads to the following results: (i) POD stabilization for coupled PDEs with Turing pattern solutions; (ii) stabilization of the POD-DEIM algorithm; (iii) computational efficiency for both online and offline stage. In particular, in the section dedicated to numerical examples, we will show that the corrected POD provides a stable and accurate method and the corrected POD-DEIM provides a fast but less accurate algorithm. The effectiveness of our method will be shown for three choices of kinetics, i.e. FitzHugh-Nagumo, Schnackenberg and DIB models, corresponding to reaction terms of increasing order of nonlinearity. In all these cases, we provide also a comparison of computational execution times between the different techniques proposed and the full model approximation.


%
The paper is organized as follows. Section \ref{sec:2}  presents a general form of the model we are interested in and presents its time discretization by an IMEX approach. In Section \ref{sec:3}, we present the POD method and one example that motivates our work. Then, Section \ref{sec:4} explains the novelty of this paper: POD and POD-DEIM with correction for the coupled nonlinear terms and the adaptive method for both techniques. A complete algorithm is presented at the end of the section. Finally, in Section \ref{sec:5}, we present our numerical results.

\section{The full model and its numerical approximation}\label{sec:2}

In this section, we first introduce the kinetics that we will consider along the current work. Hence, we briefly recall how to approximate the system \eqref{RDPDE} by the Implicit-Explicit (IMEX) Euler method and, to highlight the computational challenges, we present a numerical simulation of Turing pattern obtained by solving the full model in a typical case.

\subsection{Model kinetics}\label{sec:21}
In this study we consider three different choices of kinetics, listed below in nondimensional form, with different coupling properties, starting from the simplest linear case. In all these models the Turing instability holds for peculiar parameter choices, as we will describe in details in Section \ref{sec:5}.\\\

\noindent
{\bf FitzHugh-Nagumo (FHN) model.}  The FitzHugh-Nagumo model \cite{Gambino19FHN,KMUY21} describes the flow of an electric current through the surface membrane of a nerve fiber and it is given by
\begin{equation}
\label{fhn_kin}
\begin{aligned}
 f(u,v) & = \gamma \left(-u(u^2-1)-v\right) \\
 g(u,v) & = \gamma \left(\beta (u - \alpha v)\right),
\end{aligned}
\end{equation} 
where $u$ represents the electric potential, $v$ a recovery variable and $\alpha, \beta, \gamma >0$, $\gamma$ parameter scale. The FHN system has been proposed also in population dynamics for modeling of predator-prey interaction \cite{Padilla14}.
In \eqref{fhn_kin}, the two RD equations are coupled in linear way: the second equation has a linear kinetics and the first one is cubic only in the first variable.\\

\noindent
{\bf Schnackenberg model.} As second example, we consider the reaction-diffusion system with \emph{activator-depleted} kinetics, known also as Schnackenberg model (\cite{Madzva03}), where the simplest nonlinear coupling is present, that is quadratic in the first unknown:
\begin{equation}
\label{Schnak_kin}
\begin{aligned}
 f(u,v) &= \gamma \left(a-u+u^2 v\right), \\
 g(u,v) &= \gamma\left(b-u^2 v\right).
\end{aligned}
\end{equation}
In this model,  $u$ and $v$ represent two chemical concentrations in autocatalytic reactions and $a, b, \gamma >0$. In the biological interpretation the term $u^2v$ represents nonlinear activation of $u$ and nonlinear consumption of $v$. This model is well known in the literature as a prototype of RD-PDE system with Turing patterns of cosine-like spot type, see for example \cite{Madzva03}.\\

\noindent
{\bf DIB model.} As third example, we consider the morpho-chemical model for metal growth in electrodeposition, also known as DIB model \cite{SLB19,DIB13,DIB15,DIB17sphere}, to model phenomena arising in recharge processes in batteries with metal electrodes. The nonlinear kinetics are given by
\begin{equation}
\label{DIB_kin}
\begin{aligned}
 f(u,v) &= \rho \big{(}A_1(1-v)u - A_2u^3-B(v-\alpha) \big{)}, \\
 g(u,v) &= \rho \big{(} C(1+k_2 u)(1-v)[1-\gamma(1-v)]-Dv(1+k_3u)(1+\gamma v) \big{)},
\end{aligned}
\end{equation}
where $u$ describes the morphology of the electrodeposit, $v$ its chemical composition; here the scaling parameter is $\rho >0$. The electrochemical meaning of the parameters can be found, for example,  in \cite{DIB17sphere} and references therein. \\

For all the above models, more details about the homogeneous equilibria, which instability leads to Turing patterns, will be provided in the next sections.

\subsection{Numerical approximation of the RD-PDE system} 
\label{sec:22}
Here we present the numerical approximation of \eqref{RDPDE} in space and time that will be used in the sequel of the paper. Towards our aim, first of all, we define the so called \emph{full model}, that is the ODE system arising from the space semi-discretization that will be projected on the reduced space.

For the space discretization of \eqref{RDPDE} on the rectangular domain $\Omega$ we consider classical finite differences. Given $n_x$, $n_y$ interior mesh points along the x and y directions, respectively, the Method of Lines (MOL) with step sizes $h_x = \frac{L_x}{n_x+1}$ , $h_y = \frac{L_y}{n_y+1}$  yields to the following ODE system
\begin{equation}
\label{dyn_system}
\begin{cases}
\dot{\u} = d_u A \u+\f(\u,\v), \quad t \in (0,T],\\
\dot{\v} = d_v A \v + \g(\u,\v), \\
\u(0) = \u_0, \quad \v(0) = \v_0,
\end{cases}
\end{equation}
where the unknowns are organized in the usual vector form as $\u =\u(t) =(\u_1(t), \dots, \u_{n_x}(t))^T \in \RR^n,  \u_i(t)=(u_{i1}, \dots, u_{i,n_y})^T$ such that $u_{ij}(t) \approx u(x_i,y_j,t)$ on the given spatial meshgrid of $n = n_x n_y$ interior points. Similarly for the unknown $\v$. The discrete operator $A \in \RR^{n \times n}$ accounts for the approximation of the Laplace operator $\Delta = \partial_{xx} + \partial_{yy}$, as follows. Let $T_1 \in \RR^{n_x \times n_x}$ and $T_2 \in \RR^{n_y \times n_y}$ be matrices for the classical second order finite difference discretization of the second order derivatives along $x$ and $y$, respectively, that also include a contribution due to the approximation of the zero Neumann boundary conditions (see e.g. \cite{DSS20}). Therefore, the discrete Laplace operator on a rectangular domain can be written in Kronecker form as
$A = \frac{1}{h_x^2}(I_{n_y} \otimes T_1) + \frac{1}{h_y^2} (T_2 \otimes I_{n_x}) \in \RR^{n \times n}$, with $I_{n_x}, I_{n_y}$ being the identity matrices of dimension $n_x \times n_x$ and $n_y \times n_y$, respectively.\\
In the literature on pattern formation, several time integrators have been used to solve \eqref{RDPDE}, see e.g. \cite{Madzva03,Madzva06,DSS20}, for example the class of implicit-explicit (IMEX) schemes have been considered to approximate each equation in the ODE system \eqref{dyn_system} as a sequence of large sparse linear systems of $n$ equations. Here, to show the features of the MOR techniques in exam, we consider the simplest IMEX scheme, that is the IMEX Euler method. Given the time meshgrid $t_k =k h_t$ of time step $h_t = T/n_t$, IMEX Euler in vector form applied to \eqref{dyn_system} yields:
\begin{equation}
\label{IMEXvec}
\begin{cases}
(I_n-h_t d_u A) \u_{k+1}= \u_k+h_t \f(\u_k,\v_k), \quad k=0, \dots, n_t-1,\\
(I_n-h_t d_v A) \v_{k+1}= \v_k+h_t \g(\u_k,\v_k) \\
\end{cases}
\end{equation}
where the diffusion part of \eqref{dyn_system} is treated implicitly, while
the reaction (in general nonlinear) parts are treated explicitly
\cite{Ruuth95,Ascher95}, $\u_0, \v_0 $ are given by the initial conditions in \eqref{dyn_system}.

To recover the pattern structure we need a meshgrid with a sufficiently large number of meshpoints, e.g. $n_x, n_y \ge 50$ (see \cite{DSS20}) and final integration time $T$ sufficiently large to attain the steady state solution. Moreover, due to the explicit component of the method, it is well known (see e.g. \cite{DSS20}) that a stability bound for the choice of $h_t$ can be present. Therefore, the IMEX approach in many cases can be too expensive because the coefficient matrices in \eqref{IMEXvec} can have dimension $n \geq 2500$ and a very large number $n_t$ of discrete problems must be solved in time.

For this reason, recently in \citep{DSS20} the authors introduced matrix-oriented methods to build the time solver as a sequence of Sylvester matrix equations of dimension $n_x \times n_y$ that can be solved in the spectral space yielding a significant saving of computational execution time.
Briefly, the matrix-oriented approach can be resumed as follows. If $vec(Z)=\u$, $vec(W)=\v$, the differential matrix system equivalent to \eqref{dyn_system} is given by:
\begin{equation}
\label{matrix_system}
\begin{cases}
\dot{Z} = d_u(T_1 Z + Z T_2) + F(Z,W), \quad t \in (0,T]\\
\dot{W} = d_v(T_1 W + W T_2) + G(Z,W),\\
Z(0) = Z_0,  \quad W(0) = W_0
\end{cases}
\end{equation}
where $T_1$ and $T_2$ are described before, $F(Z,W), G(Z,W)$ are the kinetics evaluated componentwise in each spatial (interior) grid point $(x_i,y_j)$, $i=1, \dots, n_x$, $j=1,\dots,n_y$ and $(Z_0)_{i,j} = u_0(x_i,y_j)$, $(W_0)_{i,j}  = v_0(x_i,y_j)$. Then the IMEX Euler method in matrix form applied to \eqref{matrix_system} corresponds to
\begin{equation}
\label{Sylvestereq}
\begin{cases}
M_1^u Z_{k+1} + Z_{k+1}M_2^u = C_k ,\quad k=0, \dots, n_t-1\\
M_1^v W_{k+1} + W_{k+1}M_2^v = D_k , \quad
\end{cases}
\end{equation}
where $C_k = Z_k + h_t F(Z_k,W_k), \ D_k = W_k + h_t G(Z_k,W_k)$, while the coefficient matrices do not change during time evolution and are given by
$$M_1^u = I_{n_x} - h_t d_u T_1, \ M_1^v = I_{n_y} - h_t d_v T_1 \in \RR^{n_x \times n_x}  \quad  
M_2^u = -h_t d_u T_2 , \ M_2^v = -h_t d_v T_2 \ \in \RR^{n_y \times n_y}.$$ 

The solutions of the Sylvester equations in \eqref{Sylvestereq} are the matrices $Z_k$, $W_k \in \RR^{n_x \times n_y}$ which entries approximate the solutions of \eqref{RDPDE}, i.e. $(Z_k)_{ij} \approx u(x_i,y_j,t_k)$,  $(W_k)_{ij} \approx v(x_i,y_j,t_k)$ in each point $(x_i,y_j)$ at the time $t_k$. In this paper, we calculate the solutions of \eqref{Sylvestereq} (and then of the full model \eqref{dyn_system}) by the \emph{rEuler} method described in \cite{DSS20} that solves the Sylvester equations in the spectral space in a very efficient, fast and accurate way.

To highlight the computational challenges required by the approximation of Turing patterns, we present here a typical simulation for the DIB morpho-chemical model, whose kinetics are given in \eqref{DIB_kin}. For all the parameter choices, with $D = \frac{C(1-\alpha)(1-\gamma+\gamma \alpha)}{\alpha(1+\gamma \alpha)}$, there is the homogeneous equilibrium $(u_e,v_e) = (0,\alpha)$ that can undergo Turing instability \cite{DIB13}. For this example, in \eqref{RDPDE}-\eqref{DIB_kin}, we consider the parameter values
$$A_1 = 10, \ A_2 = 1, \ \alpha = 0.5, \ B = 66, \ C = 3, \ \gamma = 0.2, \ d_u = 1, \ d_v = 20, \ k_2 = 2.5, \ k_3 = 1.5, \ \rho = \frac{25}{4}.$$
The initial conditions are spatially random perturbation of the homogeneous equilibrium, i.e.
$u_0(x,y) = u_e + 10^{-5} {\tt rand}(x,y), \quad v_0(x,y) = v_e + 10^{-5} {\tt rand}(x,y),$
where {\tt{rand}} indicates the default Matlab function to generate random values with uniform distribution.
\begin{figure}[t]
\centering
\begin{subfigure}{0.9\textwidth}
\centering
\captionsetup{justification=centering}
\includegraphics[scale=0.45]{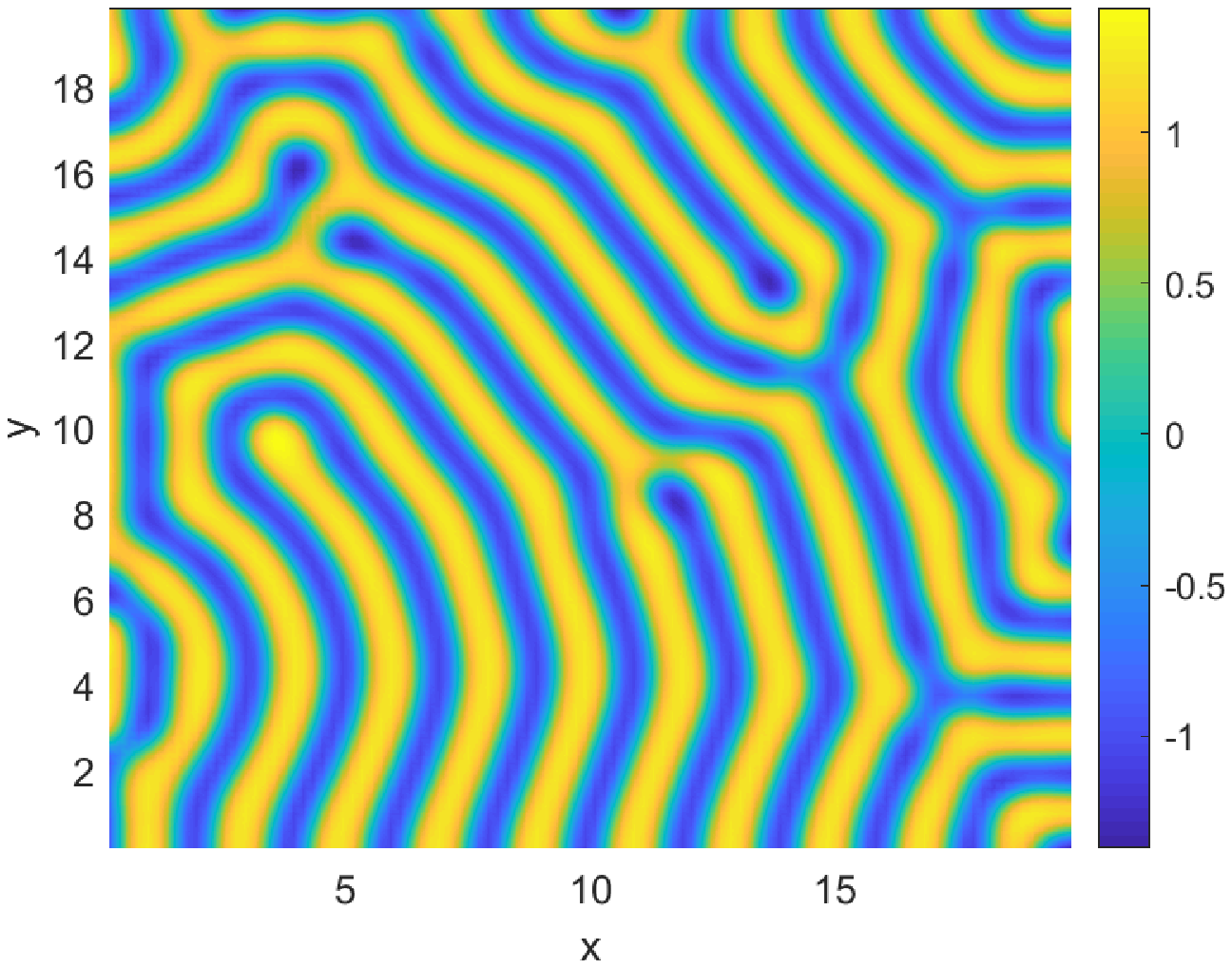}
\caption{}
\end{subfigure}
\begin{subfigure}{0.45\textwidth}
\centering
\captionsetup{justification=centering}
\includegraphics[scale=0.45]{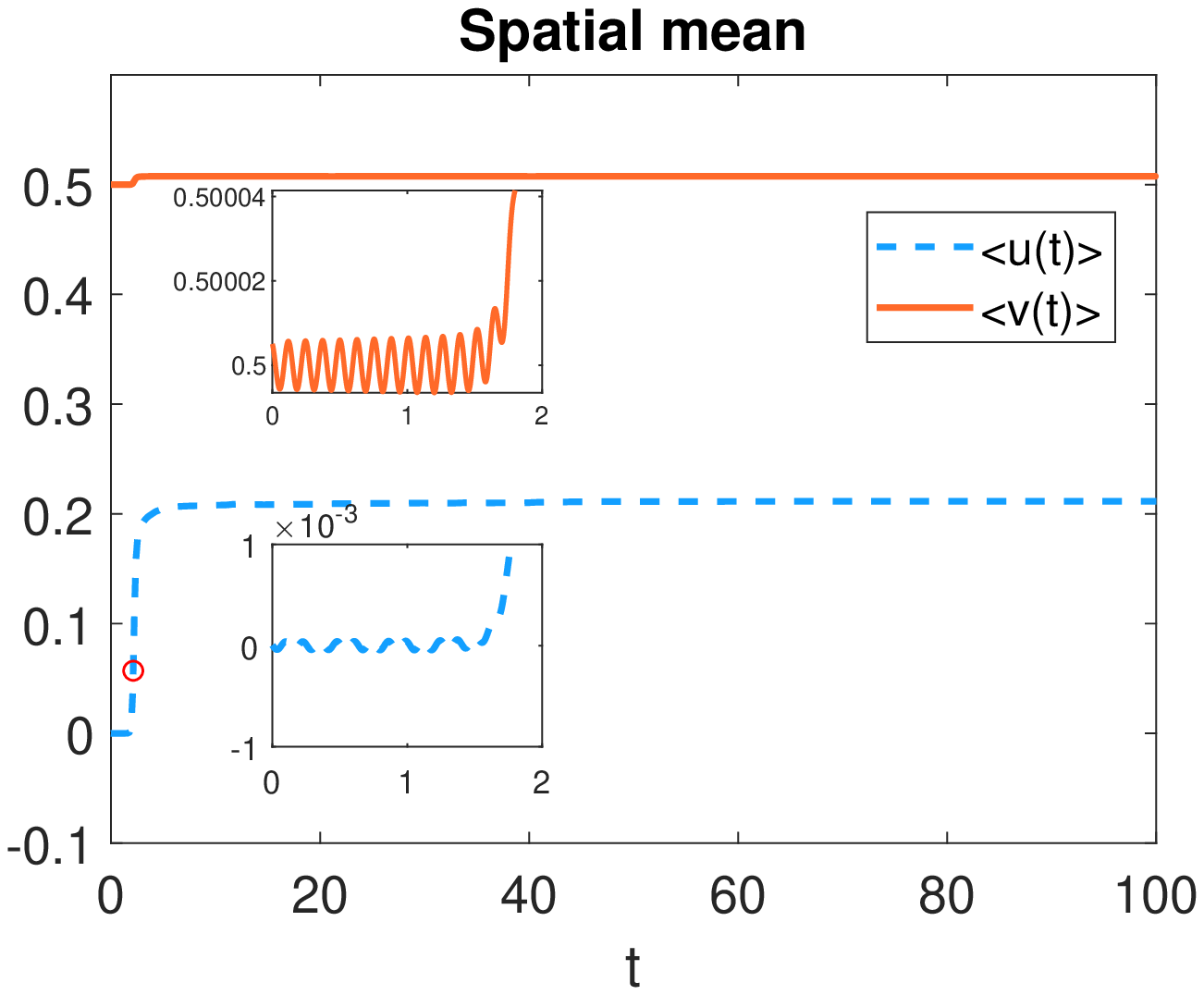}
\caption{}
\end{subfigure}
\begin{subfigure}{0.45\textwidth}
\centering
\captionsetup{justification=centering}
\includegraphics[scale=0.45]{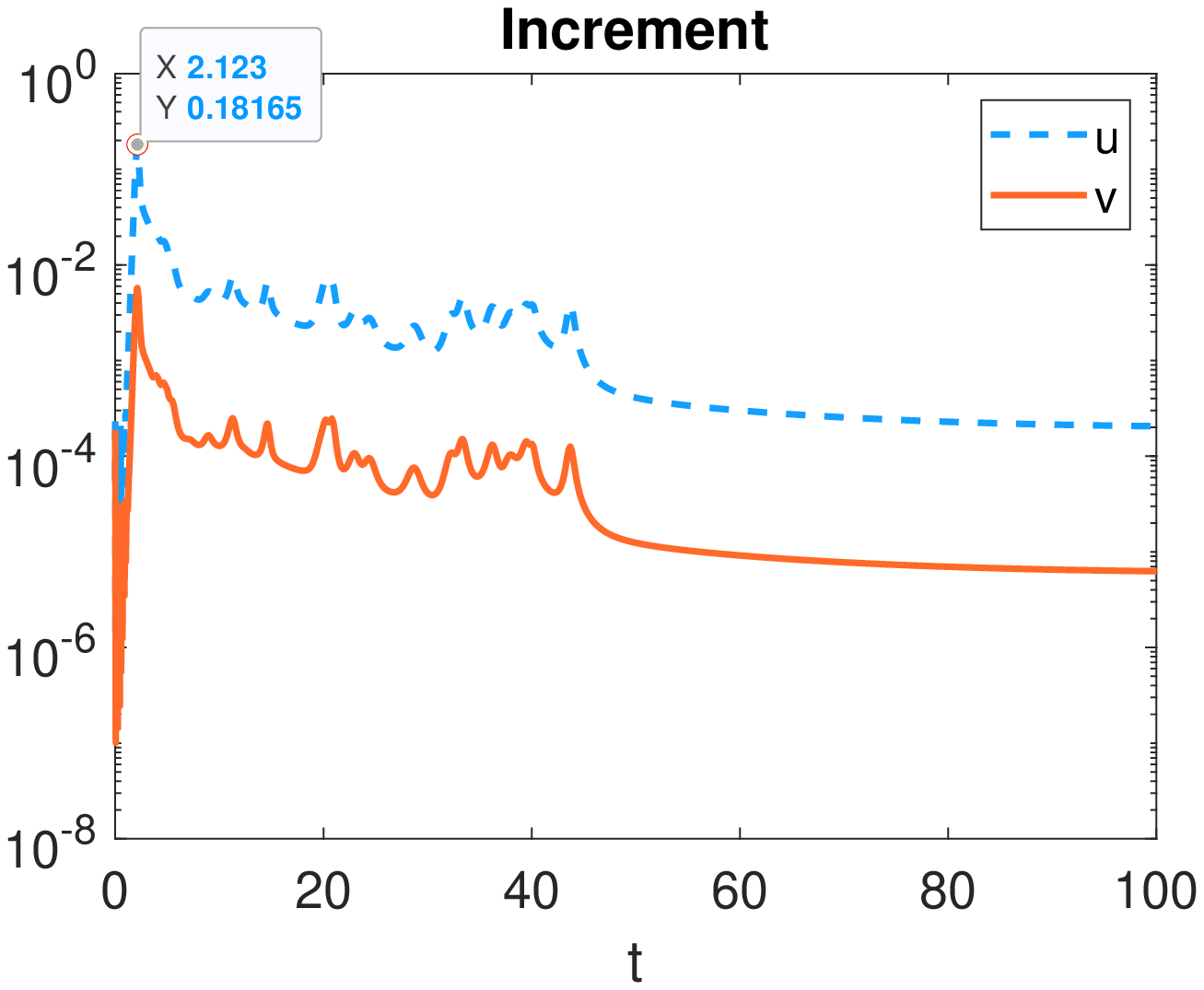}
\caption{}
\end{subfigure}
\captionsetup{justification=justified}
\caption{DIB model: (a) full model solution $u$ at the final time $T = 100$; (b) time dynamics of the spatial mean and (c) of the increment $\delta_k$ for the full model solutions $u$ and $v$. The value $\tau$, where the increment has its maximum, is indicated with a red 'o' symbol and also with a data tip in (c). The zoom insets in (b) show that in the reactivity regime the solutions depart from the initial condition with an oscillating behaviour.}
\label{fig_dib}
\end{figure}

We discretize the spatial domain $\Omega =[0,20] \times [0, 20]$ with $n_x = n_y = 100$ spatial meshpoints, such that $n = n_x n_y = 10000$ and, for stability reasons, we consider the time step $h_t = 10^{-3}$ until the final time $T = 100$, such that $n_t = 10^5$ discrete problems must be solved. 
We solve the full model \eqref{dyn_system} by using the IMEX-Euler method both in vector \eqref{IMEXvec} and matrix form \eqref{Sylvestereq} which lead to the same final labyrinth pattern shown in the Figure \ref{fig_dib}(a). The CPU time needed to solve the full model is about $326$ seconds for the matrix form and $583$ seconds for the vector one.

Moreover, we compute two indicators that will be useful also in the rest of the paper, to track and check the time dynamics until the steady state. They are given by the spatial mean $\langle u(t) \rangle$ and by the time increment of $u$, defined as follows
\begin{equation}
\label{mean}
\langle u(t) \rangle := \frac{1}{|\Omega|} \int_\Omega u(x,y,t) dx dy \approx 
\text{mean}(\u_k) = \text{mean}(Z_k), \quad  \quad k =0, \dots, n_t. 
\end{equation}
\begin{equation}
\label{increment}
\delta_k = \| \u_{k+1} - \u_k \|_F, \quad k = 0, \dots, n_t-1.
\end{equation}
In Figure \ref{fig_dib}(b)-(c), we report the dynamics of the spatial mean and of the increment for both the unknows $u$ and $v$. These behaviours are very similar for $u$ and $v$; for this reason, in the following we will consider only $u$ as reference solution.
Since a Turing pattern is an asymptotic solution of \eqref{RDPDE}, the spatial mean must attain a constant value and the increment must go to zero (see also e.g. \cite{DIB17sphere}).

In conclusion of this section, we focus on the meaning of the two time regimes discussed in the Introduction, also in view of the adaptive algorithm proposed in Section \ref{sec:4}. If $\tau$ is the time value where the maximum of the increment $\delta_k$ is attained (see Figure \ref{fig_dib}(c)), along the dynamics we can identify:
\begin{itemize}
\item the \emph{reactivity zone} $\mathcal{I}_1=[0,\tau]$: the solution departs from the initial condition showing an initial phase of oscillations and becomes unstable;
\item the \emph{stabilizing zone} $\mathcal{I}_2=(\tau, T]$: the solution starts to stabilize towards the asymptotic pattern.  
\end{itemize}

\section{Model Order Reduction and POD instability}\label{sec:3}
In this section, we briefly recall a method which reduces the dimension of the problem \eqref{dyn_system} by means of orthogonal projections. In particular, we discuss the Proper Orthogonal Decomposition (POD). The interested reader may refer to e.g. \cite{GHV20,BGW15} for a more detailed description for both continuous and discrete problems.

A numerical method to approximate a PDE already reduces the dimension of the problem: we switch from an infinite dimensional problem to a finite dimensional problem of the form \eqref{dyn_system}. However, the dimension $n$ of the semi-discretized problem is usually very large, as in the case of Turing pattern approximation presented in the previous section. The main goal of MOR is to approximate accurately the solution of \eqref{dyn_system}, that will be called \emph{full model}, by reducing its dimension. In this work, $r \ll n$ will be the dimension of the reduced problem.

For a given $r \geq 1$, a general MOR technique starts by considering two fixed matrices $\Psi_\u,\Psi_\v\in\R^{n\times r}$ such that their columns $\{(\psi_\u)_i\}_{i=1}^r$ are orthonormal vectors, that is $\Psi_\u^T\Psi_\u = I\in\R^{r\times r}$, and they form a basis for a $r-$dimensional subspace $V_{\u}^r=\mbox{span}\{(\psi_{\u})_1,\ldots, (\psi_{\u})_r\}\subset\R^n$. The matrix $\Psi_\v$ has the same properties of $\Psi_\u$ and $V_\v^r$ will be the corresponding space. Whenever we want to stress the rank $r$ of the bases we will use the notation  $\Psi_\u^r, \Psi_\v^r$. 
An appropriate choice of $V_\u^r$ and $V_\v^r$ would require that the solutions $\{\u(t),\v(t)\}$ of \eqref{dyn_system} can be approximated by a linear combination of $\{\Psi_\u,\Psi_\v\}$, that is:
 \begin{equation}\label{ansatz}
\u(t) \approx \Psi_\u \tilde{\u}_r(t),\qquad \v(t) \approx \Psi_\v \tilde{\v}_r(t),
\end{equation}
where $\{\tilde{\u}_r(t), \tilde{\v}_r(t)\}$ are unknown functions from $[0,T]$ to $\mathbb{R}^r.$ To simplify the notations, in the sequel we will use $\{\u, \v\}$ instead of $\{\u(t), \v(t)\}$.
The choice $\Psi_\u = \Psi_\v$ is also possible and corresponds to a unique subspace onto which the full model \eqref{dyn_system} can be projected. In this work, we will focus on the use of two different subspaces as discussed below.

If we plug the assumption \eqref{ansatz} into our reference problem \eqref{dyn_system} and employ the orthogonality of the bases, we obtain the following reduced system:
\begin{equation}
\label{pod_galerkin}
\begin{cases}
\dot{\tilde{\u}}_r = d_u A_r \tilde{\u}_r + \f_r(\tilde{\u}_r, \tilde{\v}_r),\qquad \tilde{\u}_r(0) = \Psi_\u^T \u_0, \\
\dot{\tilde{\v}}_r = d_v B_r \tilde{\v}_r + \g_r(\tilde{\u}_r, \tilde{\v}_r), \qquad \tilde{\v}_r(0) = \Psi_\v^T \v_0,
\end{cases}
\end{equation}
where $A_r =  \Psi_\u^T A \Psi_\u$, $B_r = \Psi_\v^T A \Psi_\v \in \RR^{r \times r}$ 
and \begin{equation} \label{fr_gr} \f_r(\tilde{\u}_r, \tilde{\v}_r) = \Psi_\u^T \f(\Psi_\u \tilde{\u}_r, \Psi_\v \tilde{\v}_r)\in \RR^{r}, \quad \g_r(\tilde{\u}_r, \tilde{\v}_r) = \Psi_\v^T \g(\Psi_\u \tilde{\u}_r, \Psi_\v \tilde{\v}_r)\in \RR^{r}.\end{equation}

We note that the Kronecker structure of the matrix $A$ in the full  model is not preserved after projection, moreover $A$ is large and sparse, while $A_r$ and $B_r$ will be small but dense. In the next simulations, the reduced system \eqref{pod_galerkin} will be solved always by the IMEX Euler method, but in vector form like in \eqref{IMEXvec}, because the matrix form and the Sylvester formulation in \eqref{Sylvestereq} are not still possible. We mention that recently a two-sided POD approach based on the matrix form has been proposed in \cite{KS20,K21}.\\

{\bf Proper Orthogonal Decomposition.} To obtain appropriate matrices $\Psi_\u,\Psi_\v$ in \eqref{ansatz},
we recall POD introduced by in \cite{Sir87}. We first need to collect data from \eqref{dyn_system}, say the analytical (if known) or approximate solution $\{\u_k,\v_k\}, k=0, \dots, n_t$ for some time instances $\{t_0,\ldots,t_{n_t}\}$. These data are usually called {\em snapshots}.  We build two snapshot matrices for the variables $\u$ and $\v$ as:
\begin{equation}\label{snap:mat}
S_\u = \begin{bmatrix} | & | & \dots & | \\
\u_0 & \u_1 & \dots & \u_{n_t} \\
| & | & \dots & |  
\end{bmatrix}, \quad  S_\v = \begin{bmatrix} | & | & \dots & | \\
\v_0 & \v_1 & \dots & \v_{n_t} \\
| & | & \dots & |  
\end{bmatrix} \in \RR^{n \times (n_t+1)}
\end{equation}
where $\u_k \approx \u(t_k)$, $\v_k \approx \v(t_k)$. 
It turns out that the left singular vectors of the (truncated) Singular Value Decomposition (SVD) \cite{GVL96} of $ S_\u$ are the {\em POD basis} of rank $r$ we are looking for. In fact, if $S_\u \approx \Psi_\u \Sigma_\u  V_\u^T$, where $\Psi_\u\in\R^{n\times r},  V_\u\in\R^{(n_t+1) \times r}$ and $\Sigma_\u \in\R^{r\times r}$ the diagonal matrix with the singular values $\sigma_i, i=1, \dots, r$, the POD basis is formed by the $r$ columns $\{(\psi_{\u})_1,\ldots, (\psi_{\u})_r\}$ of the matrix $\Psi_\u$. Equivalently, $\Psi_\v$ is obtained from the snapshot matrix $S_\v$. In this work, we prefer to use two different subspaces for $\u$ and $\v$ to compute the SVD of matrices of dimension $n\times (n_t+1)$  rather than one subspace where the SVD would be for the matrix $S=[S_\u, S_\v]\in\R^{n\times 2(n_t+1)}$. Thus, when dealing with two subspaces we obtain smaller snapshot matrices (in our simulations we will consider $n_t>10^3$). It is worth recalling that the error due to the SVD based projection is related to the neglected singular values $\sigma_i, i=r+1, \dots, {n_t+1}$ as explained in \cite{GHV20}.

As usually done in MOR there are some quantities computed {\em offline}; here we compute once and store the following objects: (i) the snapshot matrices \eqref{snap:mat}, (ii) the POD bases $\Psi_\u,\Psi_\v$ and (iii) the projected quantities $A_r, B_r$ in \eqref{pod_galerkin}. It is allowed (see e.g. \cite{BGW15}) an expensive offline stage to generate a fast (or real time) approximation of \eqref{pod_galerkin} which usually represents the {\em online stage}.\\

{\bf Discrete Empirical Interpolation Method.} The nonlinear functions $\f_r, \g_r$ in \eqref{fr_gr} require to evaluate the full, high-dimensional model reactions in the terms $\Psi_\u {\tilde \u}_r, \Psi_\v {\tilde \v}_r \in \RR^n$, thus the reduced model still depends on the full dimension $n$. To circumvent this inconvenience, the {\em Empirical Interpolation Method} (EIM, \cite{BMNP04}) and its discrete counterpart, the {\em Discrete Empirical Interpolation Method} (DEIM, \cite{CS10}), were introduced. The idea is to interpolate the nonlinear functions using only $\ell$ points. Typically, the dimension $\ell$ is much smaller than the dimension of the original problem $n$. For this goal we need to compute the interpolation basis and an operator which selects the interpolation points.
To set the algorithm for coupled problems we calculate the empirical bases $\Phi_\f=\{(\phi_\f)_1,\ldots,(\phi_\f)_\ell\}$ and $\Phi_\g=\{(\phi_\g)_1,\ldots,(\phi_\g)_\ell\}$ for the functions $\f$ and $\g$ respectively. Thus, we need to build the snapshot matrices for the nonlinear kinetics terms as: 
\begin{equation}\label{non_snap}
S_\f= \begin{bmatrix} | & | & \dots & | \\
\f(\u_0,\v_0) & \f(\u_1,\v_1) & \dots & \f(\u_{n_t},\v_{n_t}) \\
| & | & \dots & |  
\end{bmatrix},\,
S_\g = \begin{bmatrix} | & | & \dots & | \\
\g(\u_0,\v_0) & \g(\u_1,\v_1) & \dots &  \g(\u_{n_t},\v_{n_t}) \\
| & | & \dots & |  
\end{bmatrix}
\end{equation}
with $S_\f,S_\g \in \RR^{n \times (n_t+1)}$ and then we compute $\Phi_\f$ and $\Phi_\g \in \RR^{n \times \ell}$ as the POD bases of $S_\f$ and $S_\g$, respectively. 

Let us define two matrices $P_\f,P_\g \in\mathbb{R}^{n \times \ell}$ by taking $\ell$ columns of a $n \times n$ {\em permutation matrix}. Following the approach suggested in \cite{DG15}, we compute these matrices by a QR decomposition with pivoting of $\Phi_\f^T$ and $\Phi_\g^T$. Then the DEIM approximation for the nonlinear terms is given by
\begin{align*}
\begin{aligned}
\f(\Psi_\u \tilde{\u}_r, \Psi_\v \tilde{\v}_r)&\approx \Phi_\f(P_\f^T\Phi_\f)^{-1}\f(P_\f^T\Psi_\u \tilde{\u}_r,  P_\f^T\Psi_\v \tilde{\v}_r),\\
\g(\Psi_\u \tilde{\u}_r, \Psi_\v \tilde{\v}_r)&\approx \Phi_\g(P_\g^T\Phi_\g)^{-1}\g(P_\g^T\Psi_\u \tilde{\u}_r, P_\g^T\Psi_\v \tilde{\v}_r).
\end{aligned}
\end{align*}
\noindent
The matrices here involved
\begin{equation}\label{pre-comp}
P_\f^T \Psi_\u, \ P_\g^T \Psi_\u, \ P_\f^T \Psi_\v, \  P_\g^T \Psi_\v\in\mathbb{R}^{\ell\times r},\qquad
(P_\f^T \Phi_\f)^{-1}, \ (P_\g^T \Phi_\g)^{-1}\in\mathbb{R}^{\ell\times\ell} 
\end{equation}
can be pre-computed independently of the full dimension $n$. We recall that our kinetics $\f, \g$ are nonlinear functions evaluated at $\{\u(t), \v(t)\}$ component-wise. 

For given $r$ and $\ell$, the reduced ODE system for the DEIM approximation can be written as follows
\begin{equation}
\label{deim_2sottospazi}
\begin{cases}
\dot{\tilde{\u}}_r = d_u A_r \tilde{\u}_r + \Psi_\u^T \Phi_{\f}^D \f(P_\f^T\Psi_\u {\tilde \u}_r, P_\f^T\Psi_\v {\tilde \v}_r), \qquad \tilde{\u}_r(0) = \Psi_\u^T \u_0, \\
\dot{\tilde{\v}}_r = d_v B_r \tilde{\v}_r + \Psi_\v^T \Phi_{\g}^D \g(P_\g^T\Psi_\u {\tilde \u}_r, P_\g^T\Psi_\v {\tilde \v}_r), \qquad \tilde{\v}_r(0) = \Psi_\v^T \v_0,
\end{cases}
\end{equation}
where 
\begin{equation}\label{pre-comp2}
\Phi_\f^D:= \Phi_\f (P_\f^T \Phi_\f)^{-1},\qquad \Phi_\g^D := \Phi_\g (P_\g^T \Phi_\g)^{-1}\in \RR^{n \times \ell}
\end{equation}
and the quantities $\Psi_\u^T \Phi_\f^D, \ \Psi_\v^T \Phi_\g^D\in \RR^{r \times \ell}$ are also precomputed in the offline stage. This approach implies that the functions $\f,\g$ are evaluated only on $ \ell \ll n$ selected points. 

\subsection{POD instability: a numerical example}\label{sec:31}
\begin{figure}[tbp]
\centering
\begin{subfigure}{0.49 \textwidth}
\captionsetup{justification=centering}
\centering
\includegraphics[scale = 0.5]{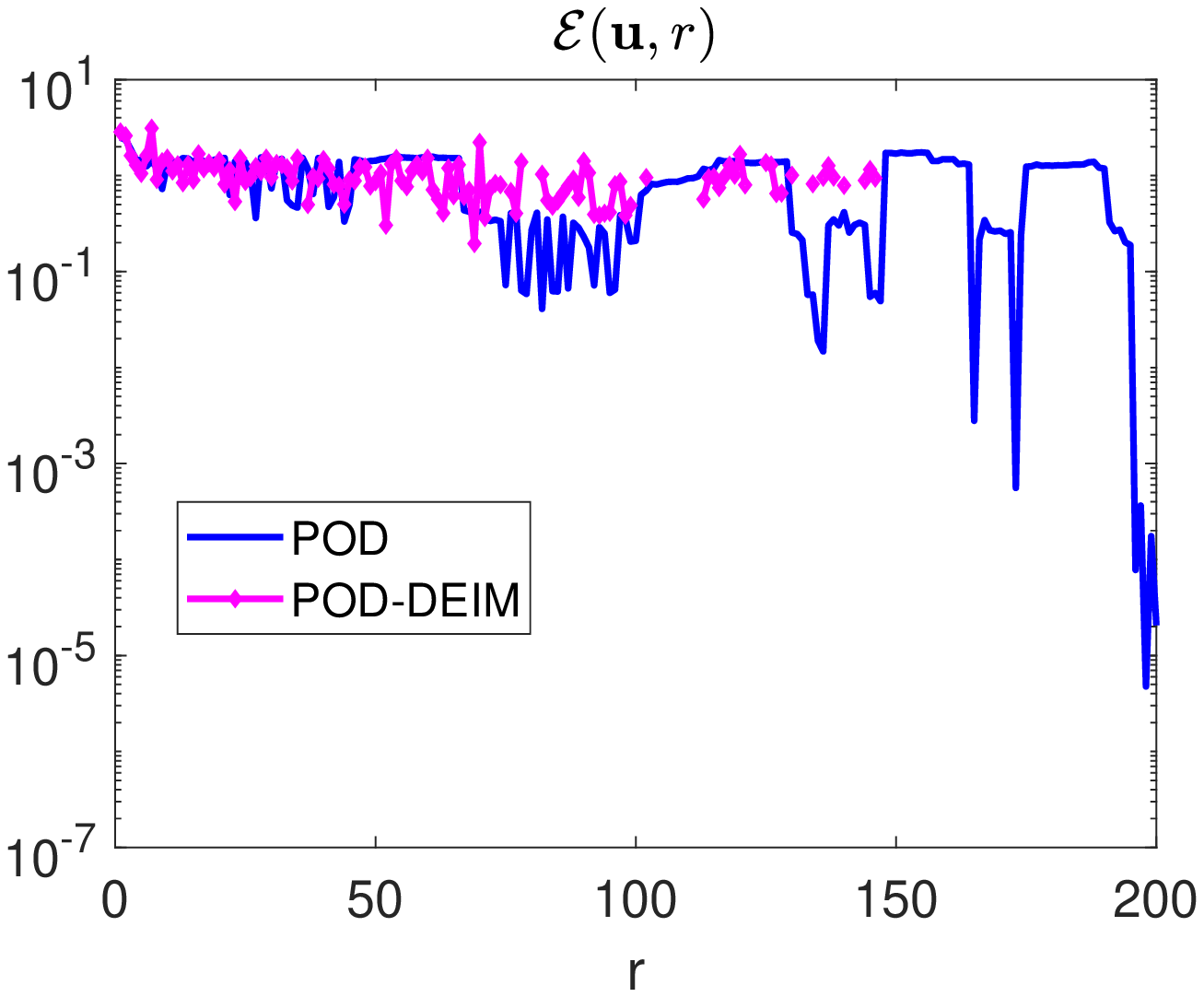}
\caption{}
\end{subfigure}
\begin{subfigure}{0.49 \textwidth}
\captionsetup{justification=centering}
\centering
\includegraphics[scale=0.5]{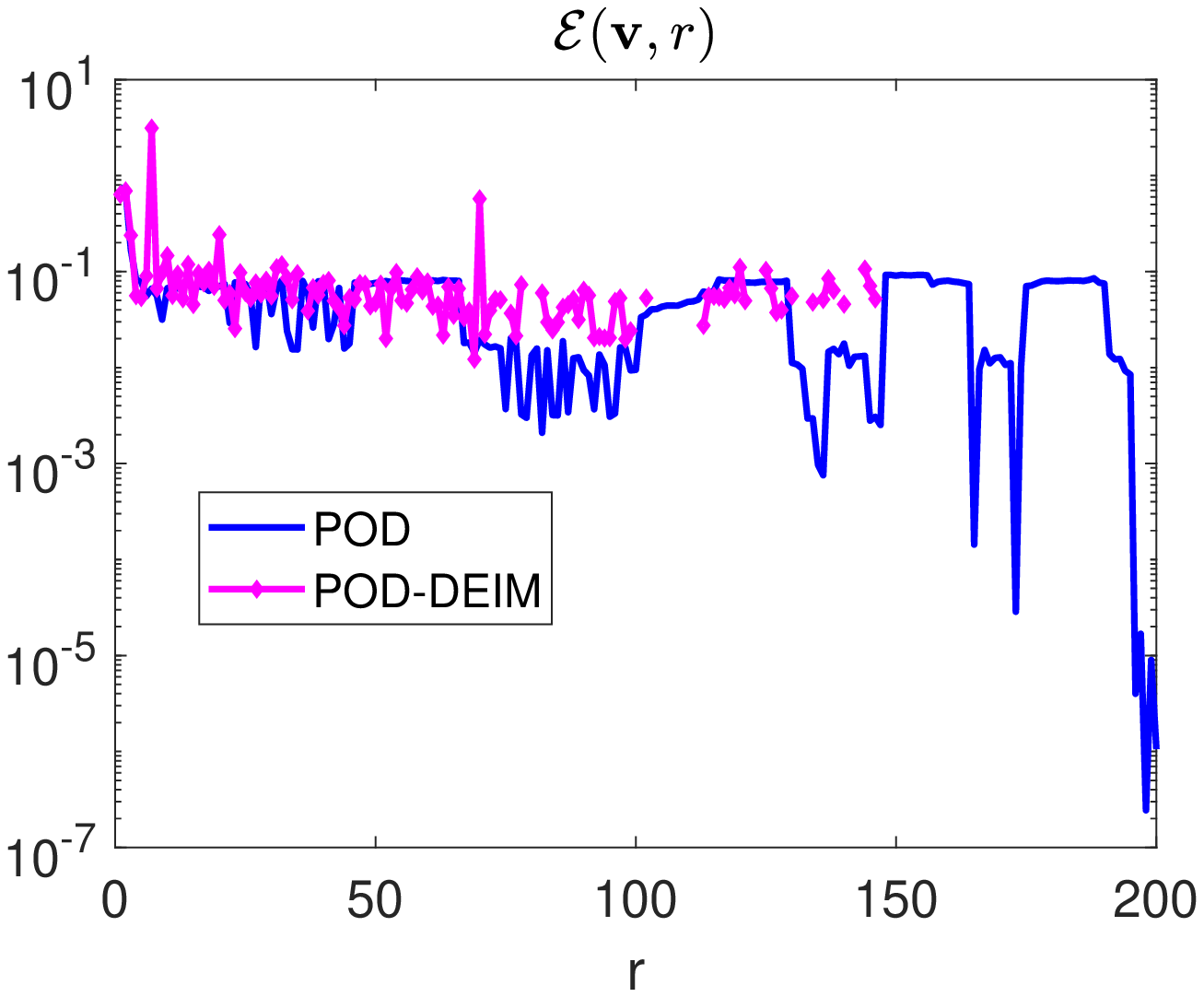} 
\caption{}
\end{subfigure}
\captionsetup{justification=justified}
\caption{DIB model. POD and POD-DEIM relative errors in \eqref{err:frob} for the variables $u$ (a) and $v$ (b) in the Frobenius norm with respect to the final pattern (reported in Fig \ref{fig_dib}(a)). Both techniques exhibit an erratic and sometimes unstable behaviour for increasing values of the dimension $r$ of the reduced space.}
\label{fig:motex2}
\end{figure}

\begin{figure}[tbp]
\centering
\begin{subfigure}{0.49 \textwidth}
\captionsetup{justification=centering}
\centering
\includegraphics[scale = 0.5]{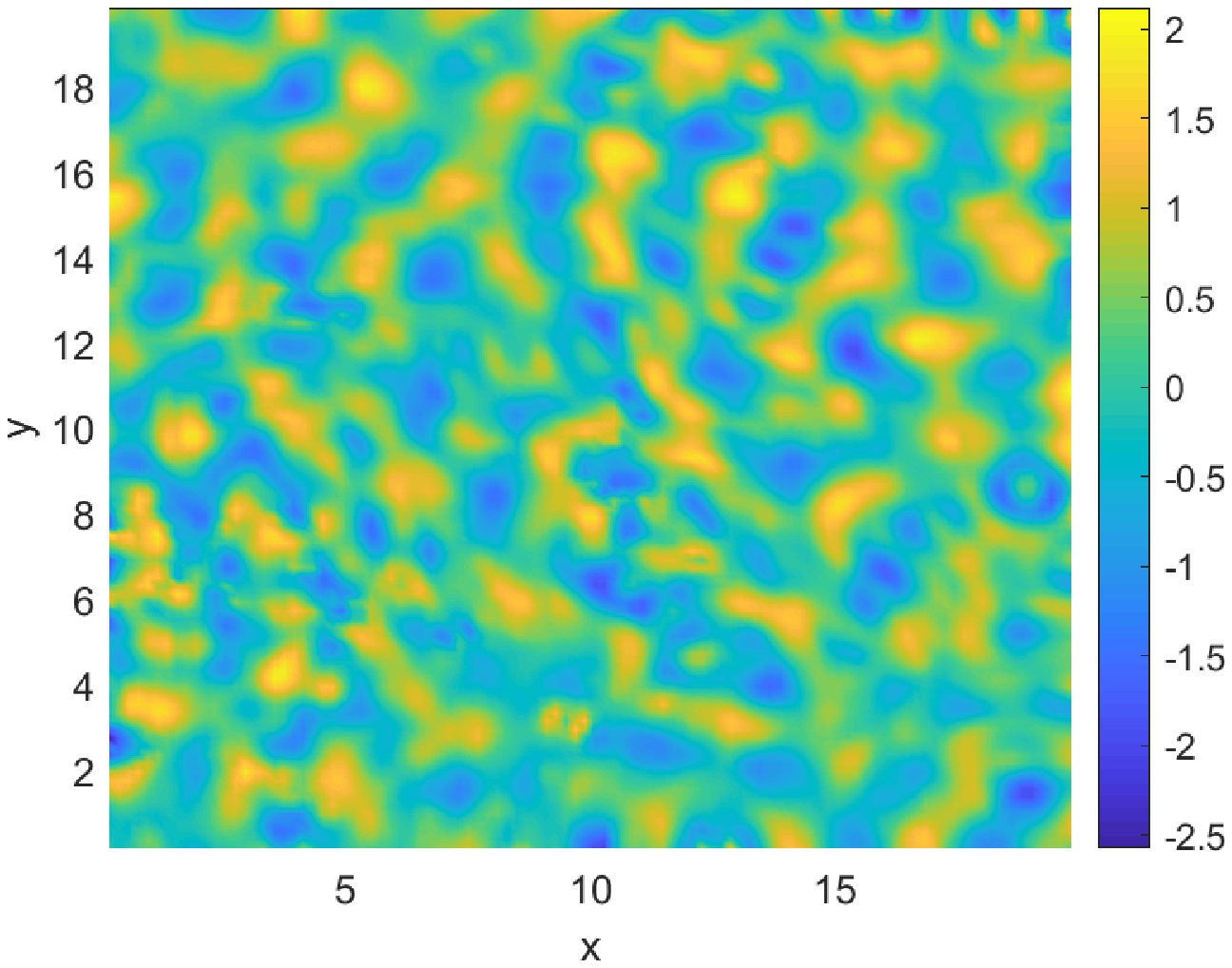}
\caption{}
\end{subfigure}
\begin{subfigure}{0.49 \textwidth}
\centering
\captionsetup{justification=centering}
\includegraphics[scale=0.5]{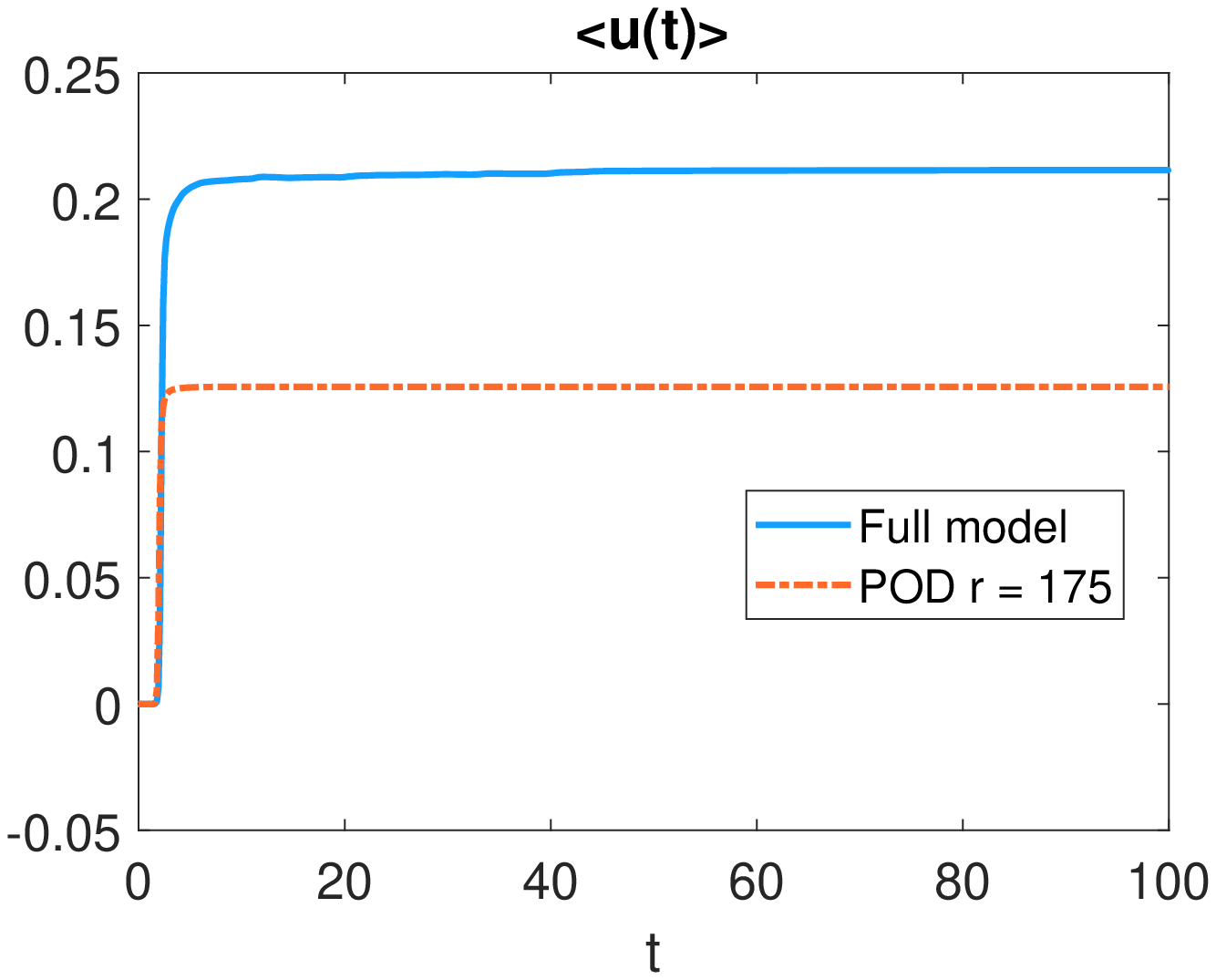}
\caption{}
\end{subfigure}
\captionsetup{justification=justified}
\caption{POD instability example \ref{sec:31} for the DIB model: (a) approximation of the pattern for $u$ at the final time $T = 100$ obtained by POD with $r = 175$  and relative error $\mathcal{E}(\u,175)=1.24$; (b) comparison of the spatial mean with respect to that of the full model (see Fig \ref{fig_dib}(b)).}
\label{fig:motex}
\end{figure}

In this section, we present a straightforward application of POD and POD-DEIM for the DIB model with kinetics in \eqref{DIB_kin} and parameters described in Section \ref{sec:22}.
For given values of $r$ and $\ell$, we calculate the reduced numerical solutions of the POD \eqref{pod_galerkin} and POD-DEIM \eqref{deim_2sottospazi} models by the IMEX Euler method described in Section \ref{sec:2} by using $h_t = 10^{-3}$. The snapshot matrices are computed by saving $\u_k$, $\v_k$ every four time steps for memory reasons.
We denote by $\u$, $\v$ the obtained solutions at the final time $T$ of integration reconstructed by \eqref{ansatz}, we calculate the relative errors $\mathcal{E}(\u,r), \mathcal{E}(\v,r)$ in the Frobenius norm with respect to the {\em reference solutions} $\u^*, \v^*$ of the full model, defined by:
\begin{equation}
\label{err:frob}
\mathcal{E}(\u,r):= \dfrac{\|\u - \u^*\|_F}{\|\u^*\|_F},\qquad \mathcal{E}(\v,r):= \dfrac{\|\v- \v^*\|_F}{\|\v^*\|_F}.
\end{equation}
In Figure \ref{fig:motex2}(a)-(b) we show the behaviour of $\mathcal{E}(\u,r)$ and $\mathcal{E}(\v,r)$ for $r=1, \dots, 200$. The reference Turing pattern $\u^*$ is reported in Fig. \ref{fig_dib}(a).
The number of DEIM points $\ell = 363$ is chosen as the maximum between the ranks of the matrices $S_\f$ and $S_\g$ defined in \eqref{non_snap}. 
The error behaviour for both POD and POD-DEIM is very erratic. Moreover, in Figure \ref{fig:motex2}, the missing information for some values of $r$ in the case of DEIM means that the corresponding reduced solutions are unstable. One usually expects a decreasing decay of the errors when the dimension of the reduced space increases. Here, we clearly show that this property does not hold also for POD, and, even worse, there are big jumps of several orders of magnitude. To give an idea of the meaning of this instability with respect to the Turing dynamics, in Figure \ref{fig:motex}(a), we show the pattern at the final time $T = 100$ obtained by POD with $r = 175$. In Figure \ref{fig:motex}(b), we show also the corresponding spatial mean $\langle u(t) \rangle$ with respect to that of the full model solution: the POD approximation completely fails because the dynamics tends towards another spatial structure, as the relative error $\mathcal{E}(\u,175)$ = $1.24$ confirms.
The goal of this example is to justify that a stabilization approach is needed for the POD and POD-DEIM methods in order to guarantee a monotone error decay and to assure that a Turing pattern is attained also in the reduced spaces. This is, in fact, the main motivation of this paper.

\section{Stabilization and adaptivity for the POD-DEIM approach}\label{sec:4}

In Section \ref{sec:31}, we have shown that both POD and POD-DEIM exhibit an irregular trend of the error for several values of $r$ to approximate the final pattern of \eqref{RDPDE} with kinetics given by \eqref{DIB_kin} (see Figures \ref{fig:motex2} and \ref{fig:motex}). We have actually found a similar behaviour for the FHN \eqref{fhn_kin} and for the Schnackenberg models \eqref{Schnak_kin} (see Section \ref{sec:5}). For this reason, our purpose here is threefold: (i) first of all to propose an algorithm to stabilize the POD behaviour of the surrogate model, (ii) to extend this algorithm to the DEIM approach and (iii) to improve the efficiency by applying these algorithms in adaptive way. 

\subsection{Correction and stabilization}
For the first goal listed above, following the approach in \cite{MLWI20}, we propose to construct a new $R-$dimensional surrogate model, $R > r$, and add a suitable correction term to the $r$-dimensional reduced system \eqref{pod_galerkin} obtained by the POD-Galerkin projection.


We consider $R$ such that $r < R \ll n$ and $\{\tilde{\u}_R, \tilde{\v}_R \}$ solutions of the reduced model \eqref{pod_galerkin} of dimension $R$.  
Let $\F(\u,\v)= d_u A \u + \f(\u,\v)$ and $\G(\u,\v)= d_v A \v + \g(\u,\v)$ be the right-hand sides (RHS) of the full model \eqref{dyn_system}. Its solution $\u$ can be approximated as in \eqref{ansatz} by choosing both $r$ and $R$. This approximation also holds true for {\v} and for the derivatives $\dot{\u}, \dot{\v}$. If we project the full model \eqref{dyn_system} onto a subspace of dimension $r$, we obtain $(\Psi_\u^r)^T \dot{\u} = (\Psi_\u^r)^T \F(\u,\v)$ and the same holds for $\v$. If we plug in the different approximations for $\u$ and $\v$, we obtain two $r$-dimensional systems. In particular, for the unknown $\u$ it holds
\begin{align}\label{corr_2}
\begin{aligned}
\centering
(\Psi_\u^r)^T (\Psi_\u^r \dot{\tilde{\u}}_r) &= (\Psi_\u^r)^T \F(\Psi_\u^r \tilde{\u}_r, \Psi_\v^r \tilde{\v}_r) + \mathcal{R}(r)\\
(\Psi_\u^r)^T (\Psi_\u^R \dot{\tilde{\u}}_R) &= (\Psi_\u^r)^T \F(\Psi_\u^R \tilde{\u}_R, \Psi_\v^R \tilde{\v}_R) + \mathcal{R}(R),
\end{aligned}
\end{align}
where $\mathcal{R}(\cdot)$ is the residual that vanishes whenever $r$ or $R$ go to infinity and it is a consequence of the error generated by \eqref{ansatz}. Since the POD bases are orthonormal, that is $(\Psi_\u^r)^T (\Psi_\u^r \tilde{\u}_r) = \tilde{\u}_r$, it follows that
\begin{equation}
\label{eq:just}
(\Psi_\u^r)^T (\Psi_\u^R \dot{\tilde{\u}}_R) = (\Psi_\u^r)^T \Psi_\u^R \dot{\tilde{\u}}_R = [I_r \ {\bf{0}} ] \dot{\tilde{\u}}_R = \dot{\tilde{\u}}_r,
\end{equation}
because the first $r$ components of $\tilde{\u}_R$ correspond to $\tilde{\u}_r$.
Therefore, by subtracting the second equation in \eqref{corr_2} from the first one, we obtain
$$ (\Psi_\u^r)^T \F(\Psi_\u^R \tilde{\u}_R, \Psi_\v^R \tilde{\v}_R)  - (\Psi_\u^r)^T \F(\Psi_\u^r \tilde{\u}_r, \Psi_\v^r \tilde{\v}_r) := \widetilde{ \mathcal{R}}_\u$$
where $\widetilde{ \mathcal{R}}_\u= \mathcal{R}(r)-\mathcal{R}(R)$. Our idea is to use this residual $\widetilde{\mathcal{R}}_\u$ to correct the original reduced model \eqref{pod_galerkin} by adding the correction terms defined as follows 
\begin{equation}
\label{correction terms}
\begin{aligned}
\f_c(\tilde{\u}_r, \tilde{\v}_r, \tilde{\u}_R, \tilde{\v}_R) &:= (\Psi_\u^r)^T \big{(} \F(\Psi_\u^R \tilde{\u}_R, \Psi_\v^R \tilde{\v}_R) - \F(\Psi_\u^r \tilde{\u}_r, \Psi_\v^r \tilde{\v}_r) \big{)}, \\ 
\g_c(\tilde{\u}_r, \tilde{\v}_r, \tilde{\u}_R, \tilde{\v}_R) &:= (\Psi_\v^r)^T \big{(} \G(\Psi_\u^R \tilde{\u}_R, \Psi_\v^R \tilde{\v}_R) - \G(\Psi_\u^r \tilde{\u}_r, \Psi_\v^r \tilde{\v}_r) \big{)}.
\end{aligned}
\end{equation}
Hence, we consider the corrected reduced system as the $r$-dimensional model using information from the $R-$dimensional solution as follows
\begin{equation}
\label{pod_correction}
\begin{cases}
\dot{\tilde{\u}}_r = d_u A_r \tilde{\u}_r + \f_r(\tilde{\u}_r, \tilde{\v}_r) + \f_c(\tilde{\u}_r, \tilde{\v}_r, \tilde{\u}_R, \tilde{\v}_R), \qquad \tilde{{\u}}_r(0) = (\Psi_\u^r)^T \u_0, \\
\dot{\tilde{\v}}_r = d_v B_r \tilde{\v}_r + \g_r(\tilde{\u}_r, \tilde{\v}_r) + \g_c(\tilde{\u}_r ,\tilde{\v}_r, \tilde{\u}_R, \tilde{\v}_R), \qquad \tilde{\v}_r(0) = (\Psi_\v^r)^T \v_0, \\
\end{cases}
\end{equation}
where $\f_r, \g_r$ are defined in \eqref{fr_gr} whereas $\f_c,\g_c$ in \eqref{correction terms}.

Here, thanks to \eqref{eq:just}, it is easy to show that the surrogate model \eqref{pod_correction}-\eqref{correction terms} is still a $r-$ dimensional model, but built upon an $R-$ dimensional one and given by:
\begin{equation}
\label{podc_simpl}
\begin{cases}
\dot{\tilde \u}_r = d_u(\Psi_\u^r)^T A \Psi_\u^R\tilde{ \u}_R +  (\Psi_\u^r)^T \f(\Psi_\u^R\tilde{ \u}_R, \Psi_\v^R\tilde{ \v}_R), \qquad {\tilde \u}_r(0) = (\Psi_\u^r)^T \u_0, \\
\dot{\tilde \v}_r = d_v (\Psi_\v^r)^T A  \Psi_\v^R\tilde{ \v}_R  + (\Psi_\v^r)^T \g(\Psi_\u^R\tilde{ \u}_R, \Psi_\v^R\tilde{ \v}_R), \qquad {\tilde \v}_r(0) = (\Psi_\v^r)^T \v_0, \\
\end{cases}
\end{equation}
where the quantities $(\Psi_\u^r)^T A \Psi_\u^R, (\Psi_\v^r)^T A  \Psi_\v^R \in\R^{r\times R}$ can be precomputed in the offline stage. Formally, the new ODE system \eqref{podc_simpl} is obtained by cancelling out the terms in \eqref{pod_correction}-\eqref{correction terms}.
The correction approach has been introduced for POD in \cite{XMRI18} where the method requires the solution of a least square optimization problem at each time instance to approximate $\f_c$ and $\g_c$. The authors in \cite{XMRI18,MLWI20} have not focused on coupled problems, but on a single equation.

It is worth noting that the \emph{new information} provided in \eqref{podc_simpl} by the correction is indeed effective if a (high) value for $R$, typically related to the rank of the snapshot matrices, is chosen. In fact, we must provide sufficiently accurate additional information. On the other hand, although this correction can stabilize the original POD technique and can yield a very accurate reduced solution, as we will show in Section \ref{sec:5}, it could make the offline stage more expensive due to the computation of $\{\tilde{\u}_R,\tilde{ \v}_R\}$.

\begin{remark}
In principle, for the correction, one could use the solution of the full model that is $R = n$ but that would make hard to store and manage all the original data. Indeed, to construct the snapshot matrices and then the POD bases $\Psi_\u^k, \Psi_\v^k$, for $k=r, R$ only a subsampling in time of the numerical solution is stored for memory reasons. To conclude it is important to choose $R\ll n$. In the rest of the paper we will refer to the surrogate model for the corrected POD technique by the acronym PODc. 
\end{remark}


In line with the objectives of a DEIM approach, as explained in Section \ref{sec:3}, to make the PODc computationally more efficient, we propose a POD-DEIM approach also for the correction model \eqref{podc_simpl}. The corresponding reduced system can be written as follows:
\begin{equation}
\label{podeim_correction}
\begin{cases}
\dot{\tilde{\u}}_r = d_u(\Psi_\u^r)^T A \Psi_\u^R\tilde{ \u}_R  + (\Psi_\u^r)^T \Phi_\f^D \f(P_\f^T \Psi_\u^R\tilde{ \u}_R, P_\f^T \Psi_\v^R\tilde{ \v}_R), \qquad {\tilde \u}_r(0) = (\Psi_\u^r)^T \u_0, \\
\dot{\tilde{\v}}_r = d_v(\Psi_\v^r)^T A \Psi_\v^R\tilde{ \v}_R  + (\Psi_\v^r)^T \Phi_\g^D\g(P_\g^T\Psi_\u^R\tilde{ \u}_R, P_\g^T\Psi_\v^R\tilde{ \v}_R), \qquad {\tilde \v}_r(0) = (\Psi_\v^r)^T \v_0,
\end{cases}
\end{equation}
where $\Phi_\f^D, \Phi_\g^D$ are defined in \eqref{pre-comp2}. 
As far as we know, the POD-DEIM correction is a novelty and it might also be interpreted as a new way to stabilize the original POD-DEIM. 
In the sequel, for brevity, we will use the acronym POD-DEIMc to indicate the surrogate system \eqref{podeim_correction} and the solution of the corrected POD-DEIM.
In Section \ref{sec:5}, we will discuss the accuracy and stabilization properties of both PODc and POD-DEIMc when applied to all three RD-PDE models introduced in Section \ref{sec:21}.
Moreover, if the computation of $\{\tilde{\u}_R, \tilde{\v}_R\}$ is considered offline as in \cite{XMRI18}, the presented numerical results will show that both PODc and POD-DEIMc are faster than the full model also for large choices of $r$. 
%
\subsection{Adaptivity} 
To lower the computational cost when integrating \eqref{podc_simpl} and \eqref{podeim_correction}, and in particular the offline costs to compute $\{\tilde{\u}_R,\tilde{ \v}_R\}$, we propose an adaptive strategy where the size of the bases for POD and PODc are updated according to the qualitative behaviour of the time dynamics. As already discussed at the end of Section \ref{sec:2}, in presence of Turing instability the time dynamics exhibits essentially two regimes \cite{DSS20}, the so-called reactivity for small times in the transient regime and the asymptotic steady state for long times. To capture the entire dynamics in $[0,T],$ a large size of the POD basis even for the correction could be required, that is very large values of $R$. For this reason, our idea is to split the integration of \eqref{pod_correction} on the two time intervals $\mathcal{I}_1:=[0, \tau]$ and $\mathcal{I}_2:=[\tau, T]$, described in Section \ref{sec:2}, where $\tau$ is calculated by the increment \eqref{increment}, as follows:
 \begin{equation}\label{tau:crit}
 \tau = \argmax_{k\in\{0,\dots, n_t-1\} }\delta_k.
 \end{equation}
Note that the increment can be easily computed by using the original snapshot data. In \cite{DSS20}, it has also been shown that $\tau$ can be related to the inflection point of the curve $\langle u(t) \rangle$ (see Figure 1). 
 
To apply the adaptive strategy, first of all, the original snapshot matrices are splitted in new ones with lower number of columns reflecting the snapshots present in $\mathcal{I}_1$  and $\mathcal{I}_2$, respectively.
We can expect that the dimensions $R_1$ and $R_2$ needed for the corrections in $\mathcal{I}_1$ and $\mathcal{I}_2$, respectively, will be much lower than $R$ required on the whole interval $[0,T]$ and, in addition, the two corresponding surrogate models will be solved on smaller time intervals.
Furthermore, for the above considerations, to find all the required projection spaces in the offline stage, the SVD of smaller matrices must be computed, thus reducing the computational costs. Randomized approaches as done in \cite{AK19} could be used to speed up the computation of the SVD algorithm, but in this work we prefer to keep our algorithm completely deterministic.



In Section \ref{sec:5}, we will show that the adaptive approach can provide a speed up of PODc and POD-DEIMc, for both offline and online stages. In fact, to sum up, in the offline stage we benefit from cheaper computations of the SVD and of the approximations of $\{\tilde{\u}_{R_i}, \tilde{\v}_{R_i}\}, i=1,2$.

Let us use the notation $\Psi_\u^{(i),r_i}, \ i=1,2$ for the adaptive POD basis where $r_i$ is the number of bases in the interval $\mathcal{I}_i$. The initial conditions $\tilde{\u}_{r_2}^{(2)}(\tau)$, and $ \tilde{\v}_{r_2}^{(2)}(\tau)$ in $\mathcal{I}_2$, are given by 
\begin{equation}\label{ic:i2}
\tilde{\u}_{r_2}^{(2)}(\tau)=(\Psi_\u^{(2),r_2})^T \Psi_\u^{(1),r_1}\tilde{ \u}_{r_1}^{(1)}(\tau), \quad \tilde{ \v}_{r_2}^{(2)}(\tau)=(\Psi_\v^{(2),r_2})^T \Psi_\v^{(1),r_1}\tilde{ \v}_{r_1}^{(1)}(\tau),
\end{equation}
such that $\tilde{ \u}_{r_1}^{(1)}(\tau), \tilde{ \v}_{r_1}^{(1)}(\tau)$ are first projected back to the original space and then projected into the reduced space of the second region.

\subsection{Our algorithm}
We summarize our adaptive POD-DEIMc procedure in Algorithm \ref{Alg:1}. We comment below our algorithm step by step distinguishing between the offline and the online stages. The adaptive PODc can be obtained neglecting the steps related to the DEIM approach.\\
\begin{algorithm}[H]
\label{Alg:1}
\hrulefill
	\caption{Adaptive POD-DEIM correction}
		\begin{algorithmic}[1] 
		\REQUIRE $\tilde h_t, h_t, A, \f, \g$
		\STATE Computation of the snapshots $\{(\u_0,\v_0),\ldots , (\u_{n_t}, \v_{n_t})\}$ with time step $h_t$ from \eqref{Sylvestereq} and stored every $\tilde h_t$ time steps, 
		\STATE Compute the increment $\delta_k$, $k=0,\ldots, n_t-1,$ and $\tau:= \argmax\limits_{k} \delta_k,$
		\STATE Set $h_t \le \tilde h_t$
		\FOR{i = 1,2}
		\IF{i = 1}
		\STATE $S_\u^{(i)}=[\u_0,\ldots, \u_{n_\tau}],\, S_\v^{(i)}=[\v_0,\ldots, \v_{n_\tau}]$ 
		\STATE $S_\f^{(i)}=[\f(\u_0, \v_0),\ldots, \f(\u_{n_\tau},\v_{n_\tau})],\, S_\g^{(i)}=[\g(\u_0, \v_0),\ldots, \g(\u_{n_\tau}, \v_{n_\tau})],$
		
		\ELSE
		\STATE $S_\u^{(i)}=[\u_{n_\tau}, \ldots, \u_{n_t}],\, S_\v^{(i)}=[\v_{n_\tau}, \ldots, \v_{n_t}]$
				\STATE $S_\f^{(i)}=[\f(\u_{n_\tau}, \v_{n_\tau}), \ldots, \f(\u_{n_t}, \v_{n_t})],\, S_\g^{(i)}=[\g(\u_{n_\tau}, \v_{n_\tau}),\ldots, \g(\u_{n_t}, \v_{n_t})],$
				\ENDIF
	\STATE Fix $R_i \approx \max \{\rank(S_\u^{(i)}), \rank(S_\v^{(i)})\}$, 
	\STATE Compute POD bases $\Psi_\u^{(i),R_i}, \Psi_\v^{(i),R_i}$
	\STATE Fix $\ell_i \approx \max \{\rank(S_\f^{(i)}), \rank(S_\g^{(i)})\}$, 
	\STATE Set $r_i < R_i$,
		\STATE Integrate problem \eqref{pod_galerkin}, with $r = R_i$ to obtain $\{\tilde{\u}_{R_i}^{(i)}, \tilde{\v}_{R_i}^{(i)}\}$
	\STATE Compute DEIM bases $\Phi_\f^{(i)}, \Phi_\g^{(i)}$ 
	\STATE Compute DEIM points $P_\f^{(i)}, P_\g^{(i)}$	
	\IF{i = 1}
	\STATE Set $\{\tilde{\u}_{r_1}^{(1)}(0),\tilde{\v}_{r_1}^{(1)}(0)\}=\{(\Psi_\u^{(1),r_1})^T \u_0, (\Psi_\v^{(1),r_1})^T \v_0\}$,
	\ELSE
	\STATE Set $\{\tilde{\u}_{r_2}^{(2)}(\tau),\tilde{ \v}_{r_2}^{(2)}(\tau)\}$ as \eqref{ic:i2}
	\ENDIF
	\STATE Integrate the model \eqref{podeim_correction} with temporal step size $h_t$
		\ENDFOR
	\end{algorithmic}
	\hrulefill
\end{algorithm}
\hrulefill

\noindent
\subsection*{Offline Stage}
{\bf Inputs.} The inputs of the algorithm are the kinetics $f,g$ for the RD-PDE models in Section \ref{sec:21}, together with the time step size $\tilde{h}_t$ and the matrix $A$ for the discrete Laplace operator in \eqref{dyn_system}.\\
\noindent
{\bf Snapshots and splitting value $\tau$.} We build the snapshot matrices using the matrix method recalled in Section \ref{sec:22}. The computation of the snapshots allows us to obtain the increment $\delta_k$ in \eqref{increment} and the time value $\tau$ to split our problem into subdomains.\\
\noindent
{\bf Set the surrogate models.} Since we split our problem, we possess (also in parallel) the different snapshot matrices $S_\u^{(i)},\ S_\v^{(i)}, \ S_\f^{(i)},\ S_\g^{(i)}$ both for $\mathcal{I}_1$ and $\mathcal{I}_2$. We then compute in each region the POD bases $\Psi_\u^{(i),R_i}, \Psi_\v^{(i),R_i}$ of rank $R_i$ and the DEIM ingredients $\Phi_\f^{(i)},\ \Phi_\g^{(i)},\ P_\f^{(i)},\ P_\g^{(i)}$ as explained in Section \ref{sec:3}. Hence in each region we compute, by applying classical POD, the solution $\{\tilde{\u}_{R_i}^{(i)},\tilde{ \v}_{R_i}^{(i)}\}, i=1,2$ to build the correction terms and the corrected surrogate models. Note that the POD bases of rank $r_i < R_i$ corresponds to the first $r_i$ columns of $\Psi_\u^{(i),R_i}, \Psi_\v^{(i),R_i}$. The choices of $R_i, \ell_i$ are described in step 12 and 14, respectively, of Algorithm \ref{Alg:1}.\\
\noindent
{\bf Projected quantities.} Once the bases are computed, we can store all the projected quantities, $A_{r_i}, B_{r_i}$, $i = 1,2$, the terms in \eqref{pre-comp} and \eqref{pre-comp2}. The initial condition for the reduced problem in $\mathcal{I}_1$ can also be stored at this stage. 
\subsection*{Online Stage}
{\bf Integration of the reduced model.} Integrate the reduced model \eqref{podeim_correction} in $\mathcal{I}_1$. Set the initial conditions in $\mathcal{I}_2$ as \eqref{ic:i2}. Finally, we can integrate the reduced problem in $\mathcal{I}_2.$ The reduced-corrected models \eqref{podeim_correction} will be integrated using the IMEX-Euler scheme in vector form as in \eqref{IMEXvec}.

\section{Numerical results} \label{sec:5}

In this section, we apply the algorithms proposed for PODc and POD-DEIMc in simple and adaptive versions to the RD-PDE models \eqref{RDPDE} with the kinetics described in Section \ref{sec:22}. Then we present three tests for the approximation of Turing patterns where the reactions are coupled with an increasing level of nonlinearity and the expected patterns have even more spatial structures. The numerical
simulations have been performed in MATLAB (ver. 2019a) on a computer DELL, i7 Intel Core processor 2.8 GHz and 16Gb RAM. In each case the ODE systems has been solved by the IMEX Euler method, in matrix form for the full models, in vector form for all typologies of reduced models occuring.
In what follows we will consider the maximum rank of the snapshot matrices \eqref{snap:mat} and \eqref{non_snap}, defined as $$\rho_{sol} = \max \{\rank(S_\u), \rank(S_\v) \},\qquad  \rho_{kin} =\max \{\rank(S_\f), \rank(S_\g) \}.$$
For each test, for POD and POD-DEIM with and without correction, we present the following results:
\begin{itemize}
\item the errors $\mathcal{E}(\u,r)$ defined in \eqref{err:frob},  obtained for all the surrogate models proposed and for $r \leq r_{max}$ where $r_{max} \leq \rho_{sol};$ 
\item the same errors for the adaptive algorithms in each interval $\mathcal{I}_1, \mathcal{I}_2$ for the variable $\u$; 
\item comparisons of the CPU execution times with and without adaptive strategy;
\item a table reporting the CPU execution times needed to achieve a certain level $tol$ of accuracy by PODc and POD-DEIMc both in adaptive and non adaptive versions.
\end{itemize}

\subsection{Test 1: FitzHugh-Nagumo model}\label{run1}
We consider the FHN reaction-diffusion model with kinetics \eqref{fhn_kin} where the expected asymptotic pattern is known thanks to the weakly nonlinear analysis in \cite[Fig.2(b)]{Gambino19FHN}. There is a unique homogeneous equilibrium $(u_e,v_e) = (0,0)$ that will undergo Turing instability. Hence, in \eqref{RDPDE}-\eqref{fhn_kin}, we consider the parameter values from \cite{Gambino19FHN}
$$d_u = 1, \quad d_v = 42.1887, \quad \alpha = 0.1, \quad \beta = 11, \quad \gamma = 65.731$$
and as initial conditions the spatially random perturbation of the homogeneous equilibrium
$u_0(x,y) = u_e + 10^{-3} {\tt rand}(x,y), \quad v_0(x,y) = v_e + 10^{-3} {\tt rand}(x,y).$ We discretize the spatial domain $\Omega=[0,\pi] \times [0, \pi]$ with $n_x = n_y = 100$ meshpoints, such that $n = n_x n_y = 10000$. For stability reasons, the full model is integrated with time step $h_t = 10^{-4}$ until the final time $T = 50$ such that $n_t = 5 \cdot 10^5$ discrete problems must be solved.

\begin{figure}[tbp]
\centering
\begin{subfigure}{0.34 \textwidth}
\captionsetup{justification=centering}
\centering
\includegraphics[scale=0.37]{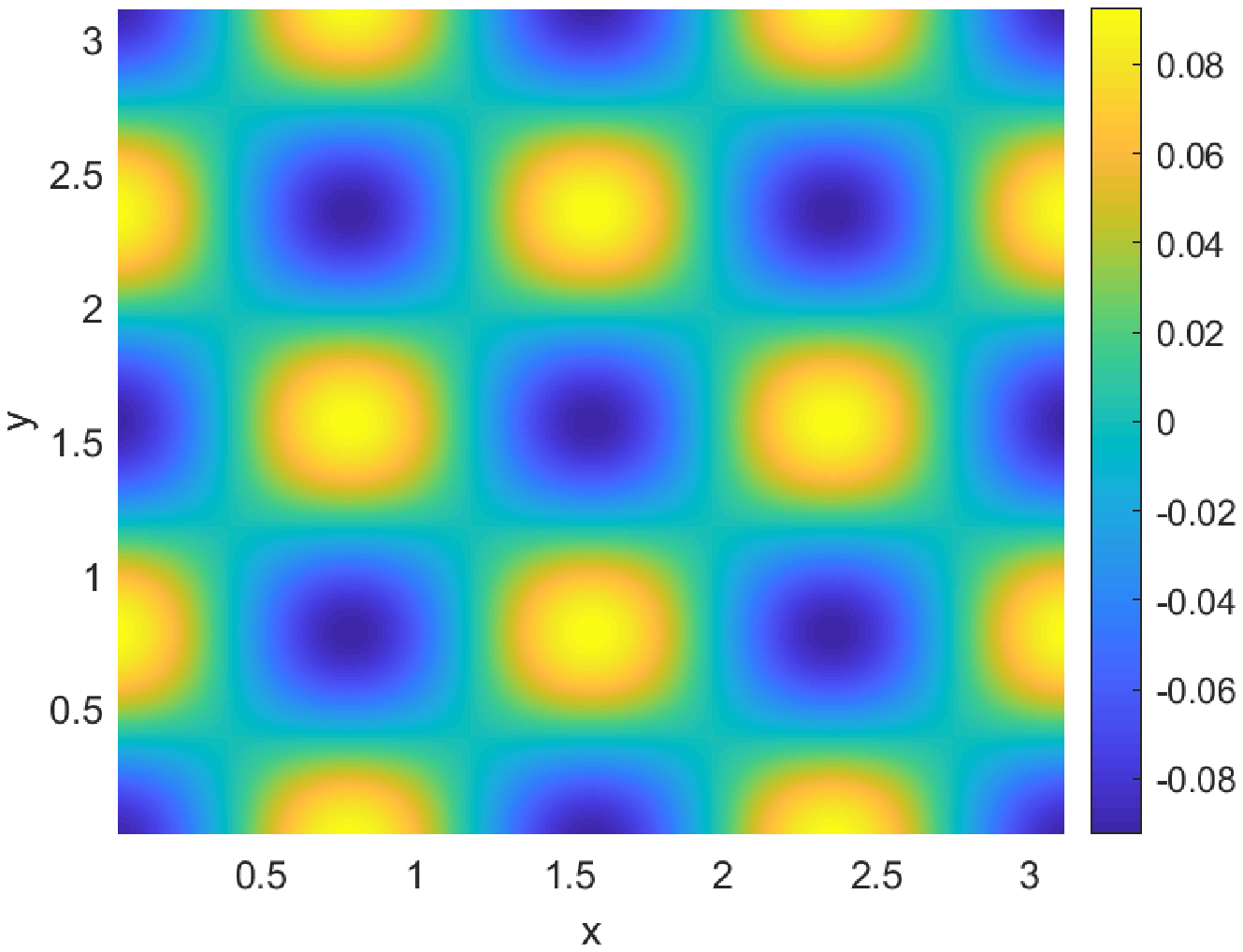}
\caption{}
\end{subfigure}
\begin{subfigure}{0.32 \textwidth}
\captionsetup{justification=centering}
\centering
\includegraphics[scale=0.38]{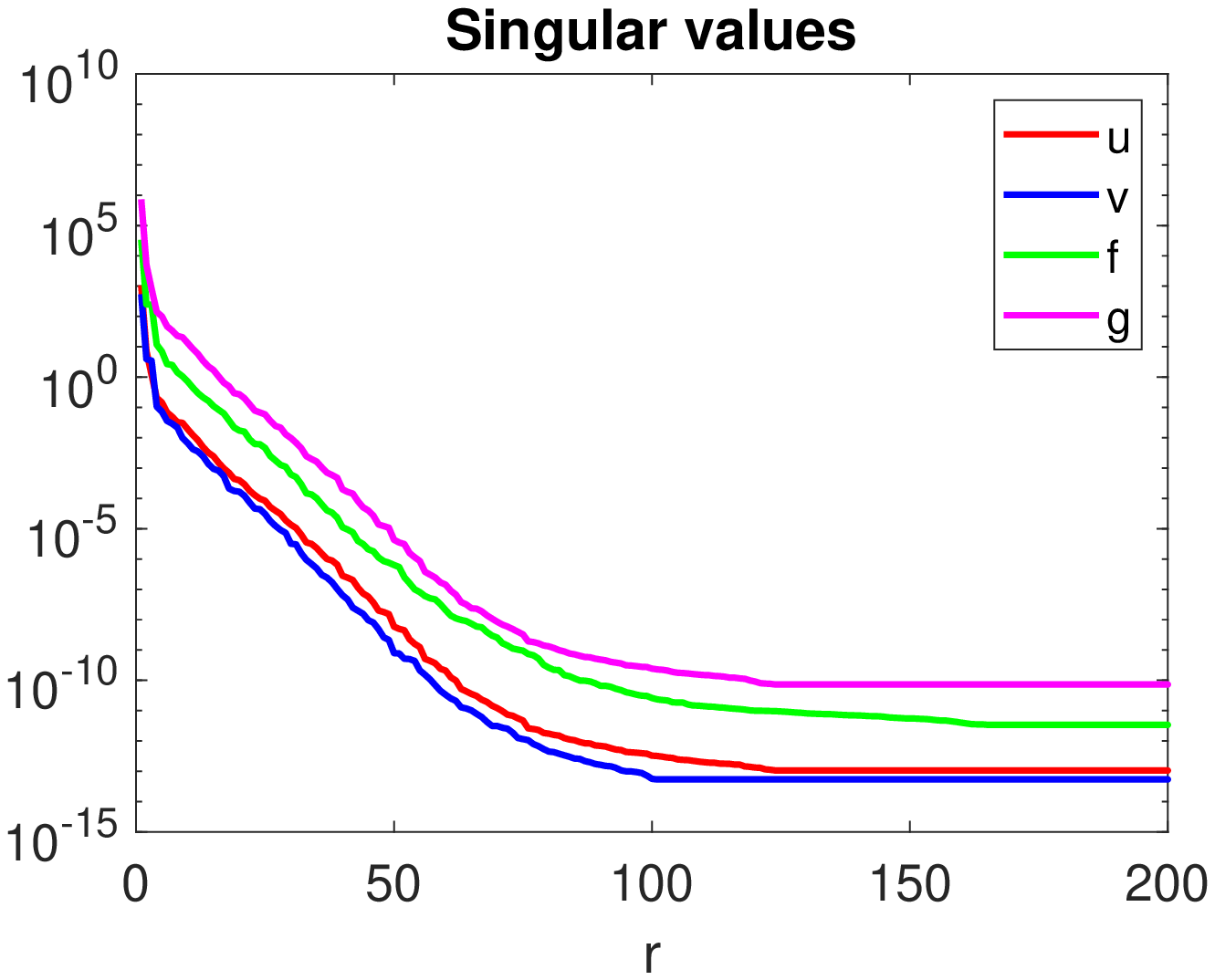}
\caption{}
\end{subfigure}
\begin{subfigure}{0.32 \textwidth}
\captionsetup{justification=centering}
\centering
\includegraphics[scale=0.38]{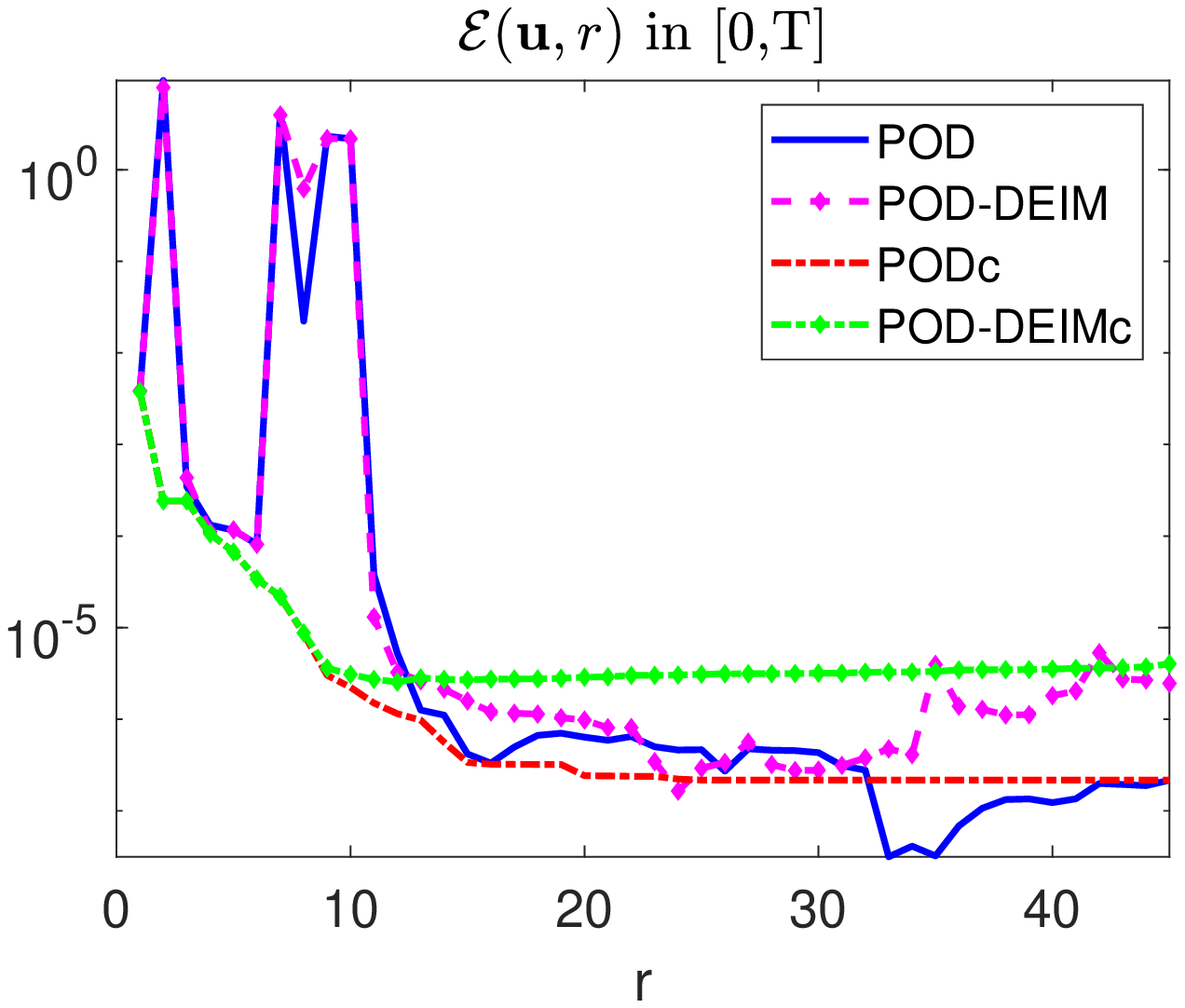}
\caption{}
\end{subfigure}
\captionsetup{justification=justified}
\caption{Test 1: FHN model. (a) Full model solution $u$ at the final time $T = 50$. (b) Singular values decay of the snapshot matrices in \eqref{snap:mat} and \eqref{non_snap}. The maximum rank is $48$. (c) Relative errors $\mathcal{E}(\u,r)$ at the final time $T = 50$ for all MOR techniques. The correction terms for PODc and POD-DEIMc are computed for $R = 45$, while the DEIM interpolation is applied with $\ell = 48$, both values near to the rank of the snapshot matrices. The non-monotone behaviour of POD and POD-DEIM errors is overcome by the corresponding corrected techniques.}
\label{fig_fhn3}
\end{figure}

The starting point of Algorithm \ref{Alg:1} is the construction of the snapshot matrices \eqref{snap:mat}. For memory reasons, we obtain them by saving the full model solutions $\u_k$, $\v_k$ every four time steps, such that $S_\u$, $S_\v \in \RR^{10000 \times 125001}$. In Figure \ref{fig_fhn3}(a)-(b), we show the full model solution $u$ at the final time $T =50$ and the singular values decay for the snapshot matrices for the solutions $u,v$ in \eqref{snap:mat} and for the kinetics $f,g$ in \eqref{non_snap}, respectively. It is easy to see that a fast decay is present, such that the machine precision is almost attained at $r=200 \ll 10^4$. Here and in what follows we calculate the rank by using the Matlab default function {\tt{rank}}\footnote{We note that the built-in Matlab function {\tt rank}(M) is implemented with a default tolerance that depends on the norm and size of the matrix $M$. This can justify why, in some cases, a subsampling of a given snapshot matrix might have larger rank than the original snapshot matrix.}. In this case we find that the maximum rank among all snapshot matrices is 48. 
This behaviour further confirms that the FHN model is the simplest one we are dealing with, due to the linear coupling of the kinetics.

{\bf Stabilization.} We begin by solving each reduced model \eqref{pod_galerkin} and \eqref{deim_2sottospazi} in the time interval $[0,T]$ with time step $h_t = 10^{-4}$. The number of DEIM points $\ell = 48$ is chosen equal to $\rho_{kin}$. 
The corrected systems \eqref{podc_simpl} and \eqref{podeim_correction} are integrated by choosing $R = 45 \approx \rho_{sol}$.
 
For increasing values of $r$ until $r_{max} = R = 45$, we calculate the relative errors \eqref{err:frob} for POD and POD-DEIM in the classical and corrected versions here introduced. The results are shown in Figure \ref{fig_fhn3}(c) for the variable $u$. The non-monotone behaviour of POD and POD-DEIM is overcome by the corresponding corrected techniques. 
In particular, the PODc error decreases and for $r \geq 20$ is almost constant around $10^{-7}$, while the error of the POD-DEIMc decreases until $r=10$ and after is slightly increasing but less than $10^{-5}$. This behaviour of the DEIM technique is well known, see  e.g. \cite{CS10}. It is worth saying that, here and for the other RD-PDE models, the errors $\mathcal{E}(\u,r)$  and $\mathcal{E}(\v,r)$ have almost the same trends for all techniques, for this reason we report only the results concerning the variable $u$.


\begin{figure}[tbp]
\centering
\begin{subfigure}{0.32 \textwidth}
\captionsetup{justification=centering}
\centering
\includegraphics[scale=0.37]{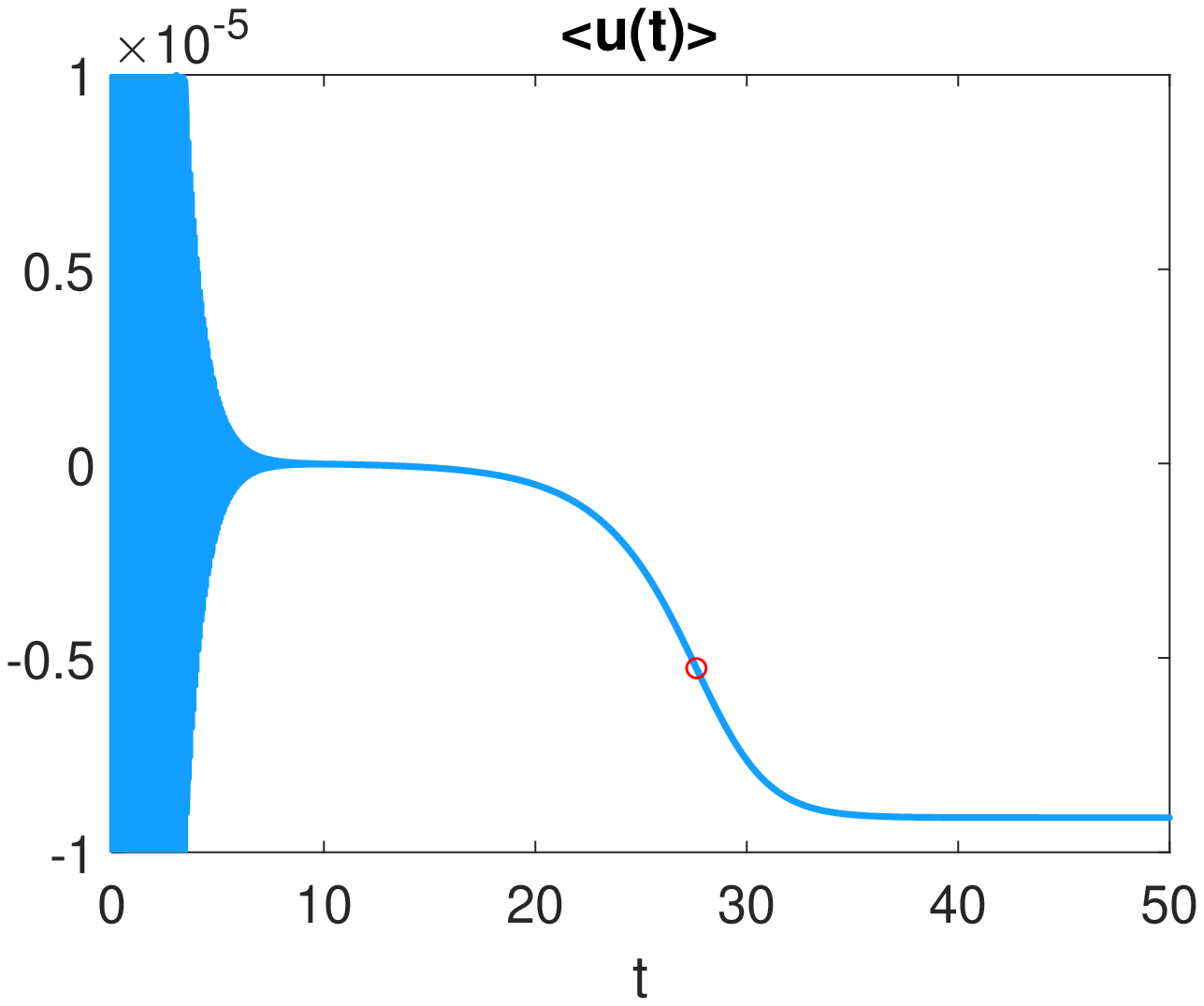}
\caption{}
\end{subfigure}
\begin{subfigure}{0.32 \textwidth}
\captionsetup{justification=centering}
\centering
\includegraphics[scale=0.38]{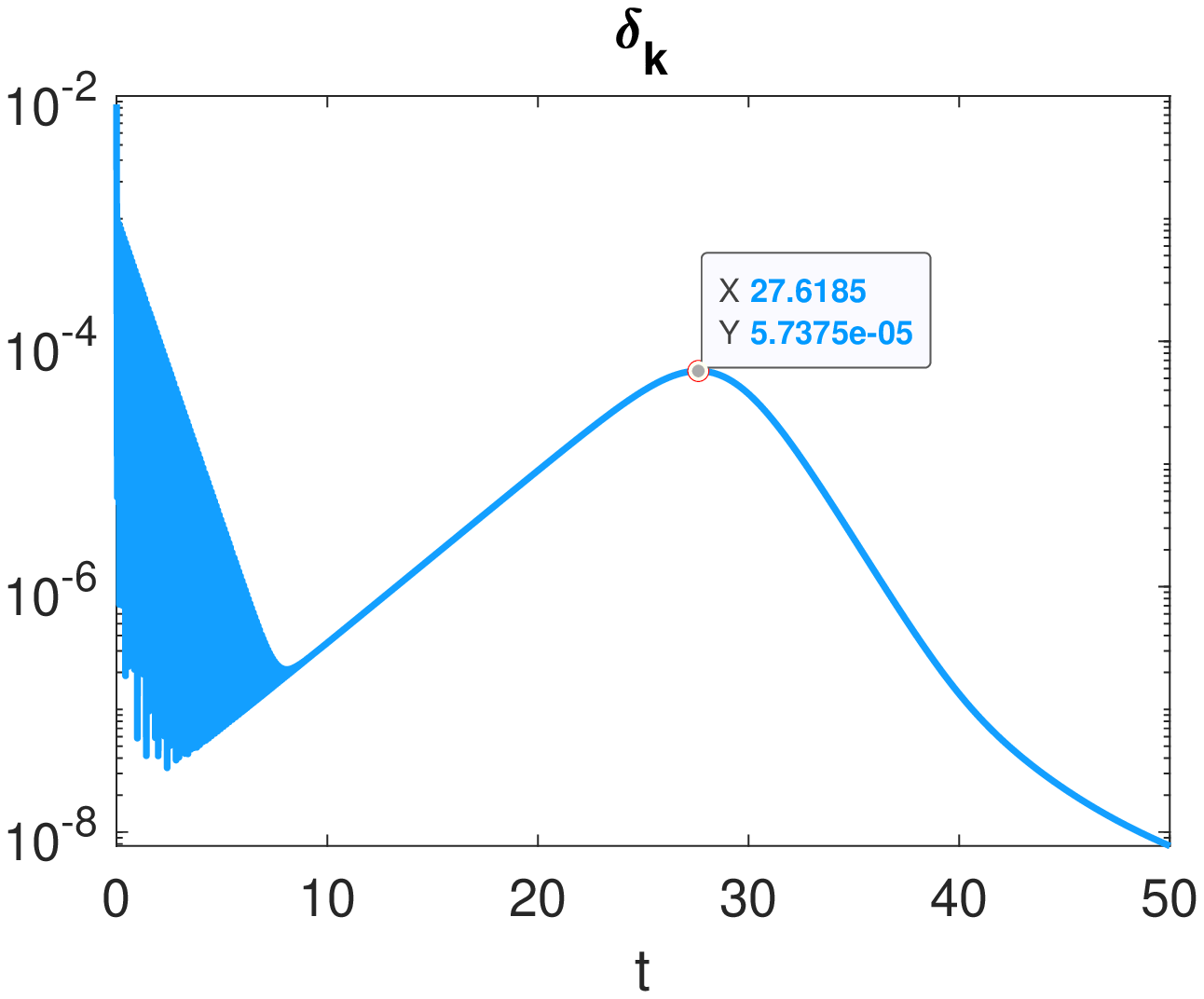}
\caption{}
\end{subfigure}
\begin{subfigure}{0.32 \textwidth}
\captionsetup{justification=centering}
\centering
\includegraphics[scale=0.37]{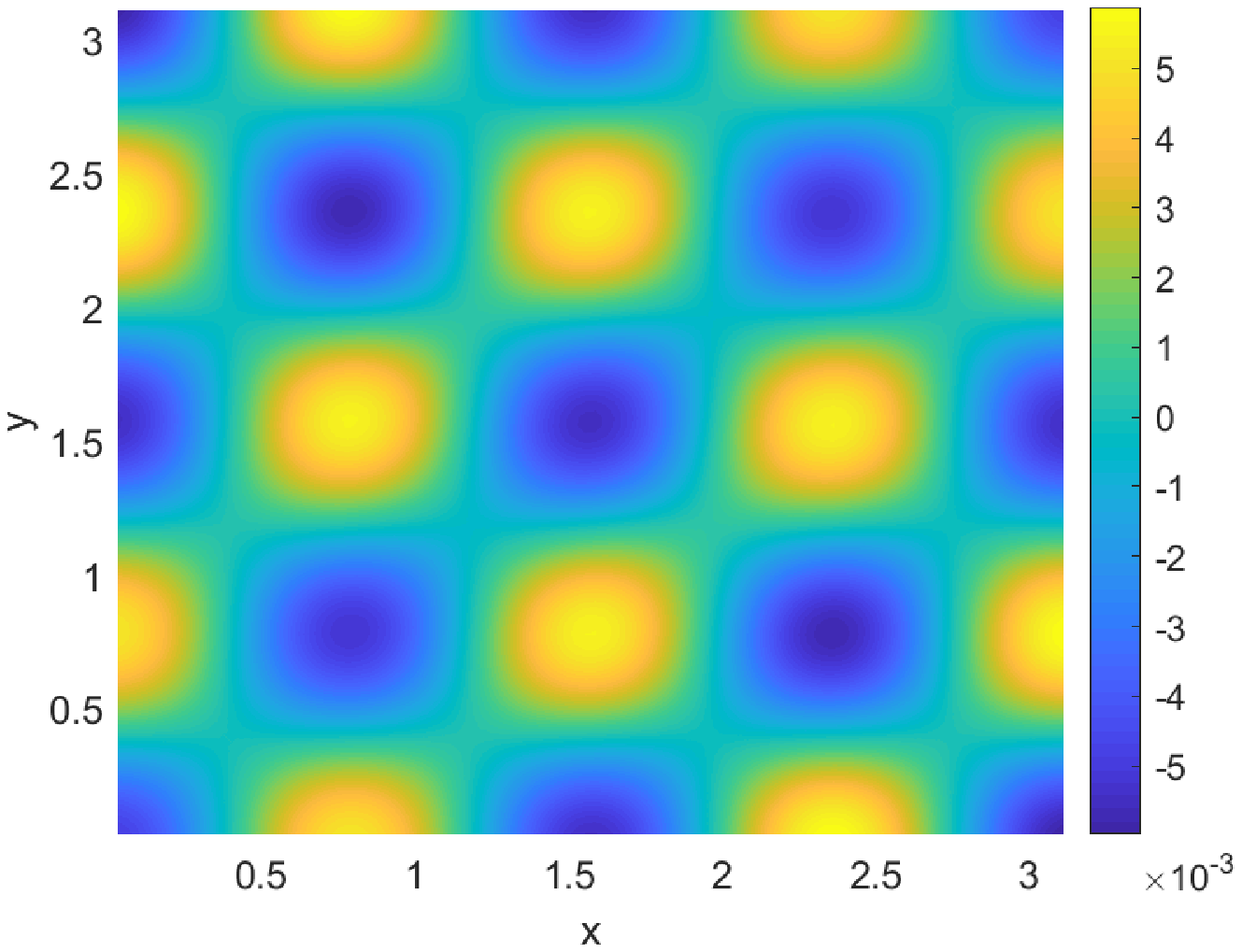}
\caption{}
\end{subfigure}
\captionsetup{justification=justified}
\caption{Test 1: FHN model. Time dynamics of the spatial mean $\langle u(t)\rangle$ (a) and of the increment $\delta_k = ||\u_{k+1}-\u_k ||_F$ (b) for the full model solution $u$. (c) Full model solution $u$ at $\bar{t} = 20$ : the pattern structure is already formed, but its amplitude is not the right one (see colorbar in Fig. \ref{fig_fhn3}(a)).}
\label{fig_fhn2}
\end{figure}

{\bf Adaptivity.} In Figure \ref{fig_fhn2}(a)-(b), we report the dynamics of the spatial mean $\langle u(t) \rangle$ in \eqref{mean} and of the increment $\delta_k$ in \eqref{increment} for the unknown $u$.  
In particular, this example shows that it is important to check information from both indicators. (For example, in \cite{KMUY21} only the mean is analysed.)


In fact, if we look in Figure \ref{fig_fhn2}(a) at the spatial mean $\langle u(t) \rangle$ for $9 \leq t \leq 20$, the stationary pattern seems to be reached. Nevertheless, it is evident that on the same interval an increasing behaviour of the increment is present; therefore the solution $u$ is not stabilized for example at $\bar{t} \approx 20$.
In Figure \ref{fig_fhn2}(c), we show the solution $u$ at $\bar{t} = 20$. It is evident that the pattern structure is already formed, but its amplitude is not the right one (see the colorbar and compare it with Figure \ref{fig_fhn3}(a)).  Looking at the increment in Figure \ref{fig_fhn2}(b), it is evident that for $t < 20$ the Turing dynamics is still in the \emph{reactivity zone} and that it attains a maximum value at $\tau \approx 27.6185$ (see red `o' symbol). Only for $t > \tau$ the increment starts to decrease towards zero and the asymptotic regime begins; in the meantime the spatial mean $\langle u(t) \rangle $ moves towards another constant value that truly indicates that the solution is stabilized.

Therefore we can apply our adaptive PODc and POD-DEIMc on the integration time intervals $\mathcal{I}_1 = [0,\tau]$ and $\mathcal{I}_2 = [\tau, T]$. In Table \ref{tab_adapt_fhn} we summarize the values $\ell$ for DEIM and $R$ for the correction chosen to solve each reduced subsystem \eqref{podc_simpl} and \eqref{podeim_correction} in $I_1$ and $I_2$. The initial conditions in $\mathcal{I}_2$ are those defined in \eqref{ic:i2} where $\u_{r_1}^{(1)}(\tau)$ and $\v_{r_1}^{(1)}(\tau)$ are the solutions of the POD-DEIMc system in $\mathcal{I}_1$ with $R_1 = 45\approx \rho_{sol}$ and $r_1 = 10\ll R_1$. In Table \ref{tab_adapt_fhn}, we show also that the adaptive strategy is able to reduce the computational cost in the offline stage. 

\begin{table}[tbp]
\centering
\begin{tabular}{ c | c | c | c }
Time interval & $\ell$ (DEIM) & $R$ (correction) & CPU time for $\{\tilde{\u}_R, \tilde{\v}_R\}$ \\
\hline
$[0,T]$ & $48$ & $45$ & $2112.5$s \\
$\mathcal{I}_1$ & $50$ & $45$ & $1234.5$s\\
$\mathcal{I}_2$ & $13$ & $12$ & $\phantom{x}188.1$s\\ 
\end{tabular}
\captionsetup{justification=justified}
\caption{Test 1: FHN model. MOR parameter values used to solve each reduced system \eqref{podc_simpl} and \eqref{podeim_correction} in the whole interval $[0,T]$ or by the adaptive strategy. In the last column is reported the cost of the offline stage for the computation of $\{\tilde{\u}_R, \tilde{\v}_R\}$. It is evident the advantage of the adaptive strategy. }
\label{tab_adapt_fhn}
\end{table}

\begin{figure}[!t]
\centering
\begin{subfigure}{0.49 \textwidth}
\captionsetup{justification=centering}
\centering
\includegraphics[scale=0.45]{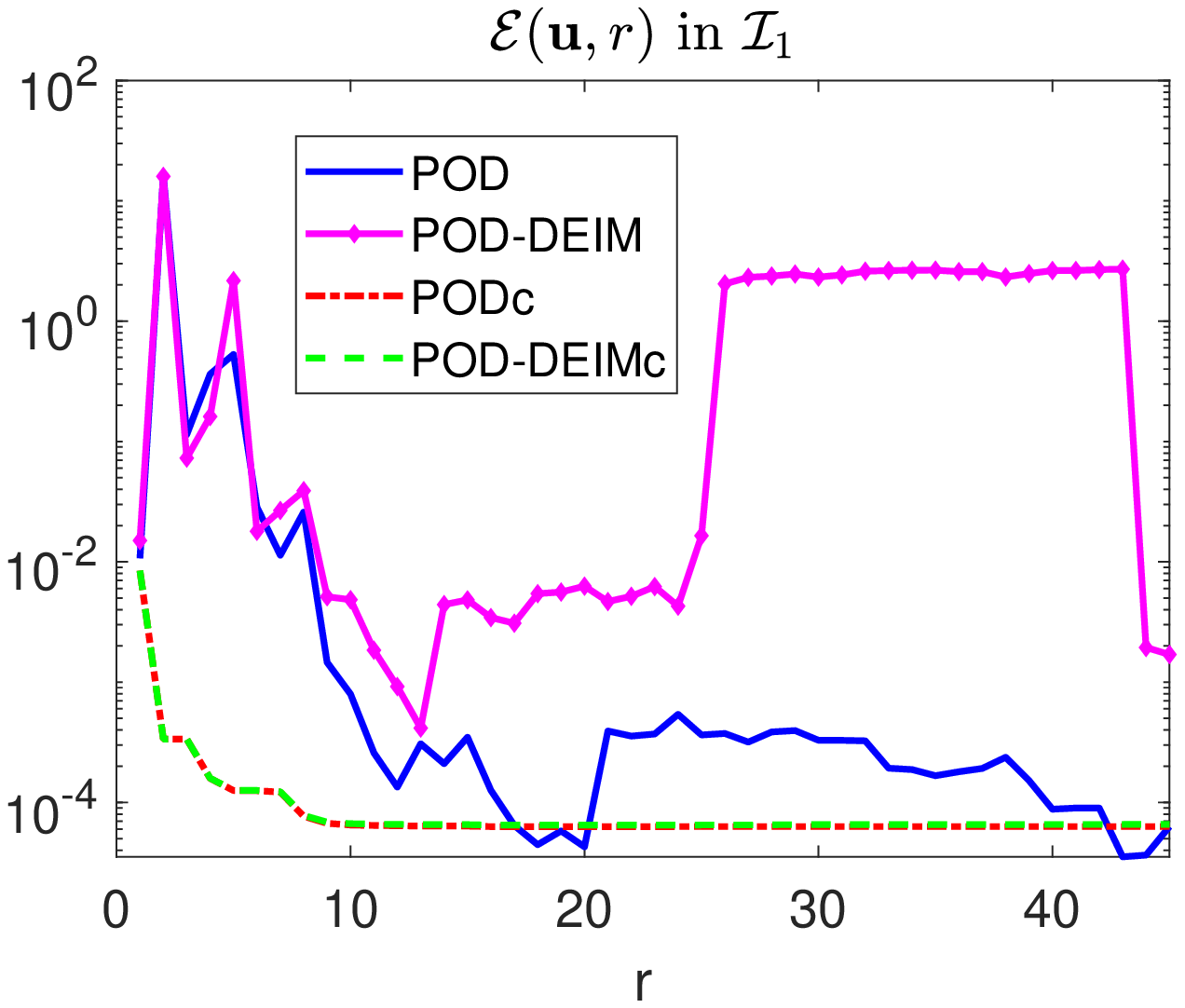}
\caption{}
\end{subfigure}
\begin{subfigure}{0.49 \textwidth}
\captionsetup{justification=centering}
\centering
\includegraphics[scale=0.45]{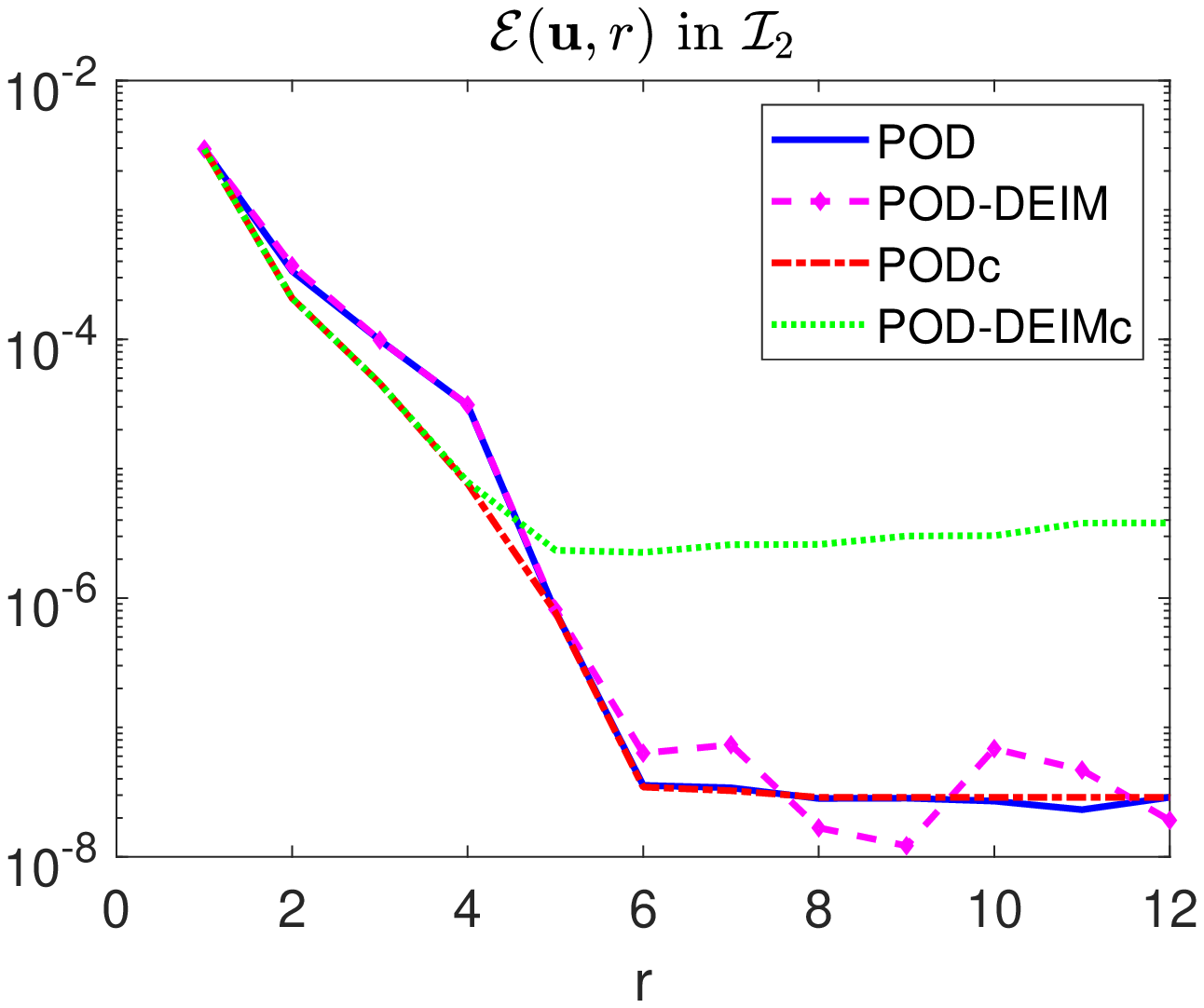}
\caption{}
\end{subfigure}
\captionsetup{justification=justified}
\caption{Test 1: FHN model adaptivity. (a) Zone $\mathcal{I}_1$: relative error $\mathcal{E}(\u,r)$ at $\tau = 27.6185$ for the unknown $u$. (b) Zone $\mathcal{I}_2$: relative error $\mathcal{E}(\u,r)$ at the final time $T = 50$. In Table \ref{tab_adapt_fhn} are reported more details for the application of POD-DEIM and its correction. In $\mathcal{I}_1$ POD-DEIM instability is clearly overcome by its corrected counterpart, that has a similar behaviour of the PODc. In $\mathcal{I}_2$ the PODc reaches the same accuracy as for the POD without correction, whereas the POD-DEIMc is also stable but less accurate for larger $r$.}
\label{fig_fhn5}
\end{figure}
In Figure \ref{fig_fhn5} we show the error behaviours. In $\mathcal{I}_1$ (see Figure \ref{fig_fhn5}(a)) the instability of POD-DEIM is more evident than in the whole interval $[0,T]$, as shown in Fig \ref{fig_fhn3}(c). Indeed in $\mathcal{I}_1$ the time dynamics has its major variability. PODc and POD-DEIMc stabilize this bad behaviour and both tend to a constant error of order $10^{-4}$ for $r \geq 10$.
In $\mathcal{I}_2$ (see Figure \ref{fig_fhn5}(b)) POD is not unstable with respect to $r$, while POD-DEIM is slightly oscillating for $r \geq 6$.  POD and PODc have almost the same trend and both achieve an error of order $10^{-8}$ for $r_{max}=12=R_2\approx\rho_{sol}$. While POD-DEIMc has an error of order $10^{-6}$ and does not improve its performance obtained on $[0,T]$. Hence, in terms of accuracy, the adaptive strategy improves only the accuracy of PODc that moves from an order of $10^{-7}$ to $10^{-8}$.
%

\begin{figure}[tb]
\centering
\begin{subfigure}{0.49 \textwidth}
\captionsetup{justification=centering}
\centering
\includegraphics[scale=0.48]{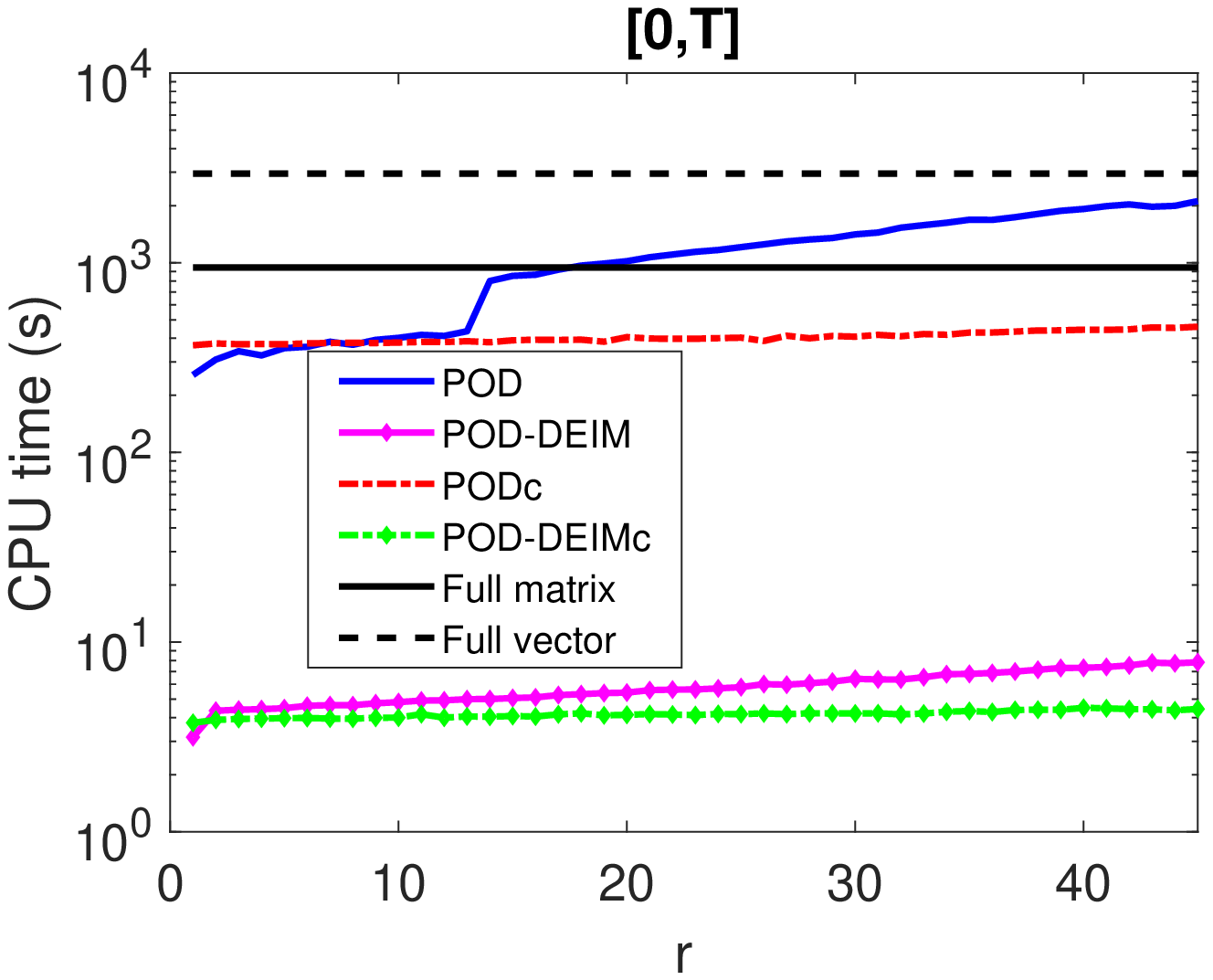}
\caption{}
\end{subfigure}
\begin{subfigure}{0.49 \textwidth}
\captionsetup{justification=centering}
\centering
\includegraphics[scale=0.48]{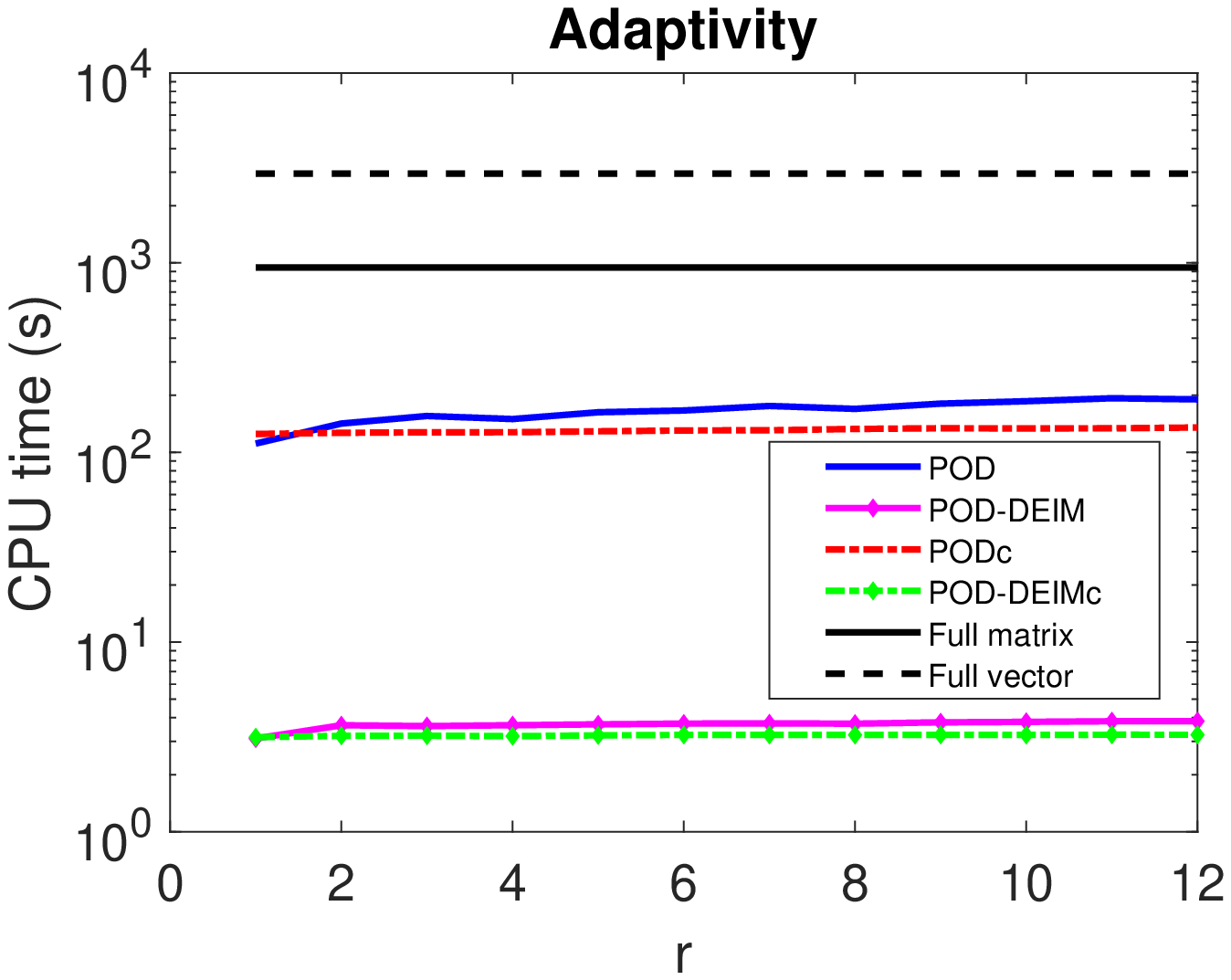}
\caption{}
\end{subfigure}
\captionsetup{justification=justified}
\caption{Test 1: FHN model. Computational costs of all the techniques proposed. (a) CPU time (s) for solving each reduced model in the whole time interval $[0,T]$; (b) total cost of the adaptive algorithms solving the reduced systems in zone $\mathcal{I}_1$ and $\mathcal{I}_2$ with $R_1 = 45$ and $R_2 = 12$. The computation of $\{\tilde{\u}_R, \tilde{\v}_R\}$ for the corrected systems is not included because is considered offline. It is clear that POD benefits from the adaptivity and results to be very competitive with respect to the full model, for any choices of $r$. In this case, its corrected version is slightly cheaper.}
\label{fig_fhn7}
\end{figure}

{\bf Computational cost: online stage.} To conclude we compare the results in terms of computational cost in the online stage. For increasing values of $r$, in Figure \ref{fig_fhn7}(a) we show the CPU time (in seconds) needed to solve the reduced models in the time interval $[0,T]$, while in Figure \ref{fig_fhn7}(b) the cost of the adaptive algorithms including the cost for solving the subsystems in $\mathcal{I}_i$, $i = 1,2$. As a reference, the black continuous line represents the cost of the full model solution when solved by the IMEX-Euler in matrix form, whereas the cost of the vector form is drawn by the black dashed line. It is easy to see that each reduced and corrected technique, not only stabilizes the classical POD and POD-DEIM, but turns out to be faster than them and than the full model also for larger values of $r$. Moreover, in the adaptive case the classical POD and PODc have a almost the same cost (see Figure \ref{fig_fhn7} (b)). We remind that in Figure \ref{fig_fhn7} the computation of $\{\tilde{\u}_R,\tilde{\v}_R\}$ is not included because it is considered offline. 
Table \ref{tab_adapt_fhn} reports the choices of $R$ in each time interval and the computational cost of $\{\tilde{\u}_R, \tilde{\v}_R\}$. The adaptive algorithm results to be competitive also in the offline stage, with a speed-up factor of 1.48.

\begin{table}[bp]
\centering
\begin{tabular}{ c c c|c c | c c | c c | c c }
 tol & $10^{-3}$ & $r_0$ & $10^{-4}$ & $r_0$&  $10^{-5}$& $r_0$ &  $10^{-6}$ & $r_0$ & $10^{-7}$ & $r_0$\\
\hline
PODc & $375.30$s & $\phantom{x}2$ & $371.68$s & $\phantom{x}5$ &  $378.02$s & $\phantom{x}8$ & $385.89$s & $13$ & - & - \\
PODc adaptive & $126.89$s & $\phantom{x}2$ & $127.68$s & $\phantom{x}3$ & $127.77$s & $\phantom{x}4$ & $129.19$s & $\phantom{x}5$ & $130.56$s & $\phantom{x}6$ \\
POD-DEIMc & $\phantom{xx}3.89$s & $\phantom{x}2$ & $\phantom{xx}3.97$s & $\phantom{x}5$ & $\phantom{xx}3.95$s & $\phantom{x}8$ & - &- & -&-  \\
POD-DEIMc adaptive & $\phantom{xx}3.21$s  & $\phantom{x}2$  & $\phantom{xx}3.21$s & $\phantom{x}3$ & $\phantom{xx}3.19$s  & $\phantom{x}4$  & - &- & -&- \\
\end{tabular}
\captionsetup{justification=justified}
\caption{Test 1: FHN model. CPU time needed by the PODc and POD-DEIMc to achieve a desidered tolerance for the relative error of the unknown $u$. The best gain is obtained by the PODc in adaptive version for $r = 6$, while the POD-DEIMc in adaptive and global versions is cheaper than PODc.}
\label{tab_fhn}
\end{table}

In Table \ref{tab_fhn} we report the CPU time needed in the online stage to achieve a desired accuracy $tol$ by the corrected algorithms. We also report for each $tol$ the values $r_0$ such that $\mathcal{E}(\u,r) \leq tol$ for $r \geq r_0$. We can deduce that PODc clearly benefits from the adaptive approach. In fact, the CPU time has a speed-up factor of about $3\times$ for any desired $tol$ and PODc achieves higher accuracy (see last column). 
Concerning the POD-DEIMc, for all tolerances we have almost the same cost and there is no difference if we apply or not the adaptivity. On the other hand, we have already discussed that the main difference in the computational cost lies in the offline stage, because of the different choices of $R$ in each time interval (see Table \ref{tab_adapt_fhn}). Not surprisingly, the POD-DEIMc in both form is faster than PODc. Finally, it is worth recalling (see Fig. 7) that the CPU time to solve the full model is $943$ seconds in matrix form and $2948$ seconds in vector form. This further stresses the efficiency of the surrogate models and of the stabilizing algorithms here proposed. \\
To conclude, by considering the cost in both offline and online stage, in this test the adaptive POD-DEIMc results to be the best method presenting an error $\mathcal{E}(\u,r) \leq 10^{-5}$ for $r \geq 4$. However, if a higher accuracy is desired, we can use both POD and PODc with an adaptive strategy and $r \geq 5$, with an increased computational cost in the online stage, but in any case competitive with respect to the full model. 

\subsection{Test 2: Schnakenberg model }\label{run2}
\begin{figure}[tbp]
\centering
\begin{subfigure}{0.34 \textwidth}
\captionsetup{justification=centering}
\centering
\hspace{-1cm}
\includegraphics[scale=0.38]{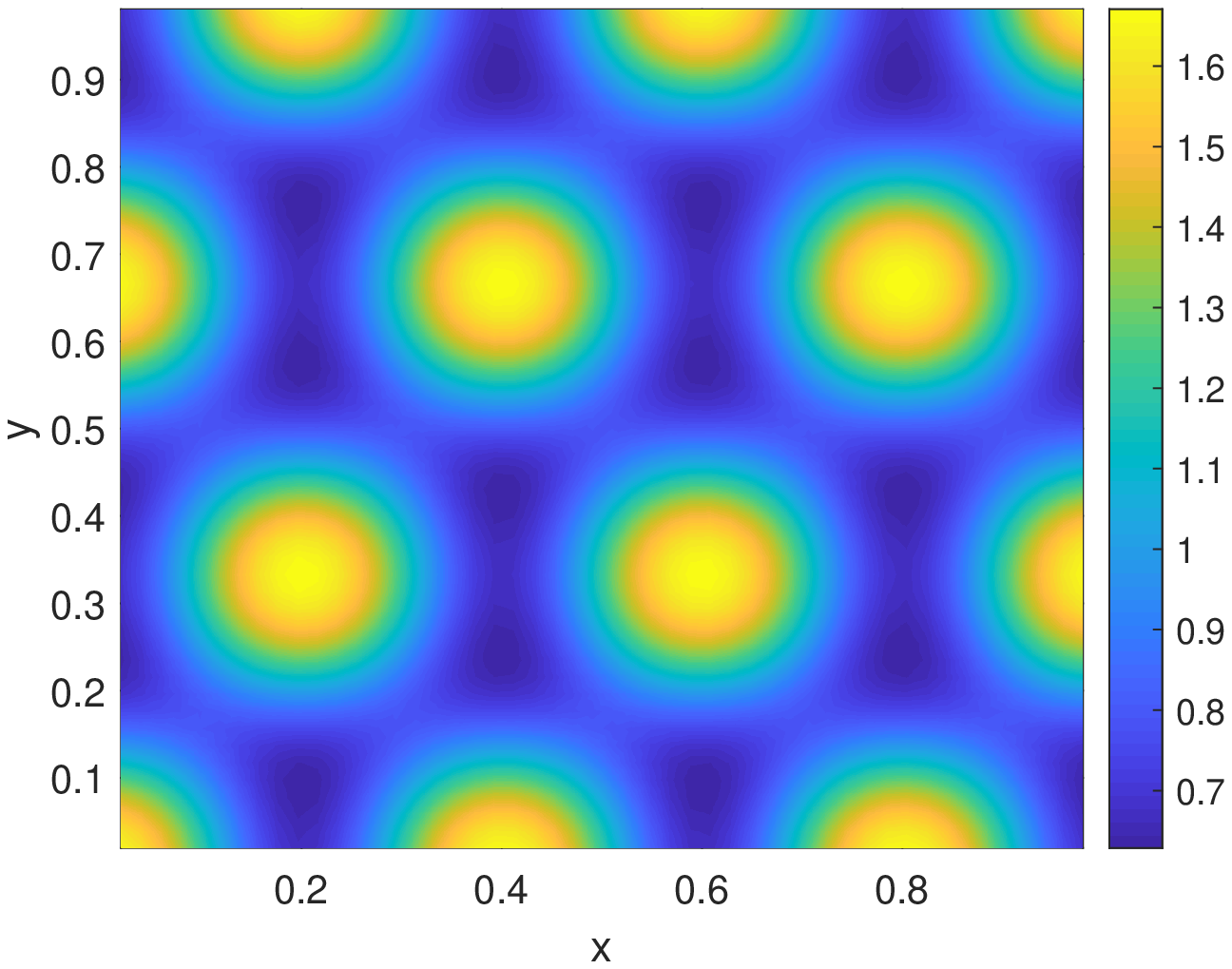}
\caption{}
\end{subfigure}
\begin{subfigure}{0.32 \textwidth}
\captionsetup{justification=centering}
\centering
\hspace{-0.5cm}
\includegraphics[scale=0.4]{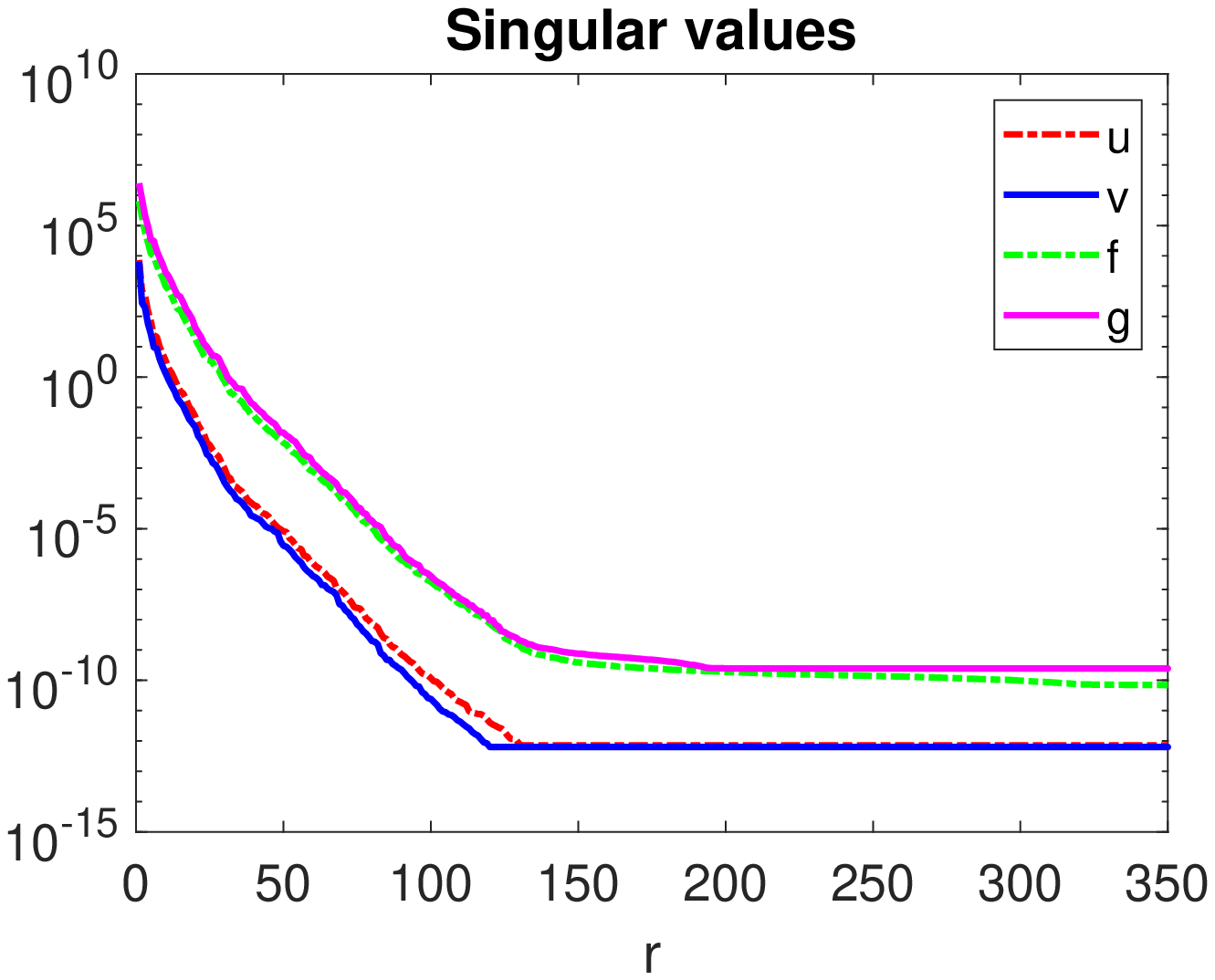}
\caption{}
\end{subfigure}
\begin{subfigure}{0.32 \textwidth}
\captionsetup{justification=centering}
\centering
\includegraphics[scale=0.4]{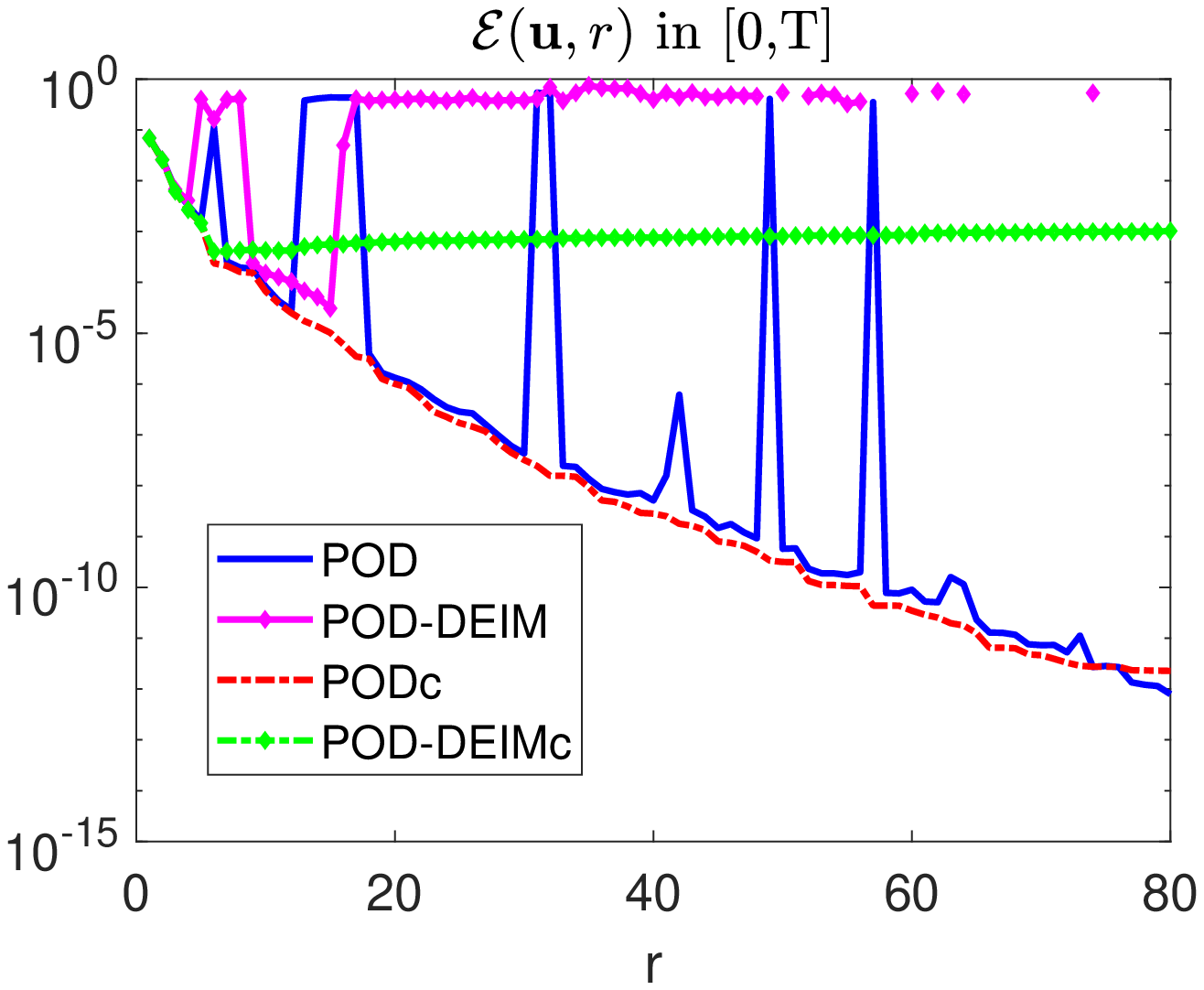}
\caption{}
\end{subfigure}
\captionsetup{justification=justified}
\caption{Test 2: Schnakenberg model. (a) Full model solution $u$ at the final time $T = 2$. (b) Singular values decay of the snapshot matrices. (c) Relative error $\mathcal{E}(\u,r)$ at the final time $T = 2$ for all MOR techniques. The correction terms are computed by choosing $R = \rho_{sol}=80$. DEIM is always applied with $\ell =  \rho_{kin}=90$. The POD error seems to decay, but exhibits big jumps for several values of $r$, whereas POD-DEIM is completely unstable. PODc and POD-DEIMc are able to stabilize these bad behaviours. As for the FHN model, POD-DEIMc shows for $r \geq 10$ a stagnation of the error around $10^{-3}$.}
\label{fig_schnak1}
\end{figure}

We consider the Schnakenberg RD-PDE model with kinetics in \eqref{Schnak_kin}, which unique homogeneous equilibrium $u_e = a+b$ and $v_e = \frac{b}{(a+b)^2}$ undergoes Turing instability.  
In \eqref{RDPDE}-\eqref{Schnak_kin} we choose the parameter values (see \cite{Madzva03,BMS21}):
$$d_u = 1, \ d_v = 10, \ a = 0.1,\ b = 0.9,\ \gamma = 1000$$ 
and as initial conditions a small random perturbations of $(u_e,v_e)$, 
$u_0(x,y) = u_e + 10^{-5} {\tt{rand}}(x,y), \ v_0(x,y) = v_e + 10^{-5} {\tt{rand}}(x,y).$
The spatial domain $\Omega=[0,1] \times [0,1]$ is discretized with $n_x = n_y = 50$ interior points ($n = n_x n_y = 2500$).
The full model is integrated in time by applying the IMEX-Euler scheme in matrix-oriented form (see Section \ref{sec:22}) on the time interval $[0,T]=[0,2]$ with time step $h_t = 10^{-4}$ such that $n_t = 20000$ solution maps are calculated for $u$ and $v$. For memory reasons, we construct the snapshot matrices by saving the solutions $\u_k$ and $\v_k$ only every four time steps, such that $ S_\u$, $S_\v \in \RR^{2500\times 5001}$. The obtained Turing pattern for $u$ at the final time $T=2$ is shown in Figure \ref{fig_schnak1}(a) and the singular values decay of the snapshot matrices \eqref{snap:mat} and \eqref{non_snap} in Figure \ref{fig_schnak1}(b).
The maximum rank of the snapshot matrices $S_\u$, $S_\v$  is $\rho_{sol}=80$, whereas for the kinetics $S_\f$, $S_\g$ is $\rho_{kin}=90$, in fact the singular values decay is slower than in the FitzHugh-Nagumo model.


{\bf Stabilization.} We first integrate the reduced models \eqref{pod_galerkin} and \eqref{deim_2sottospazi} on the full time interval $[0,T]$ with time step $h_t = 10^{-4}$ and for the projection dimension $1 \leq r \leq 80 = \rho_{sol}$. We apply DEIM with the parameter $\ell = 90=\rho_{kin}$ (see Section \ref{sec:3}). To build the corrected form PODc \eqref{podc_simpl} and POD-DEIMc \eqref{podeim_correction} we choose $R = 80 \approx \rho_{sol}$. In Figure \ref{fig_schnak1}(c) we show the relative error $\mathcal{E}(\u,r)$ for the unknown $u$ in \eqref{err:frob} for increasing values of the dimension of the reduced space $r$. Similar behaviour has been obtained for the variable $v$ (not reported).


It is clear that POD exhibits big jumps of high order of magnitudes for several values of $r$ and that POD-DEIM is completely unstable. These drawbacks are overcome if we consider the corresponding corrected approaches (drawn in dash-dot lines) that imply for both PODc and POD-DEIMc a decreasing error decay. Neverthless the error of the POD-DEIMc tends to be constant around $10^{-3}$ for $r \geq 10$. 

\begin{figure}[tbp]
\centering
\begin{subfigure}{0.49 \textwidth}
\captionsetup{justification=centering}
\centering
\includegraphics[scale=0.45]{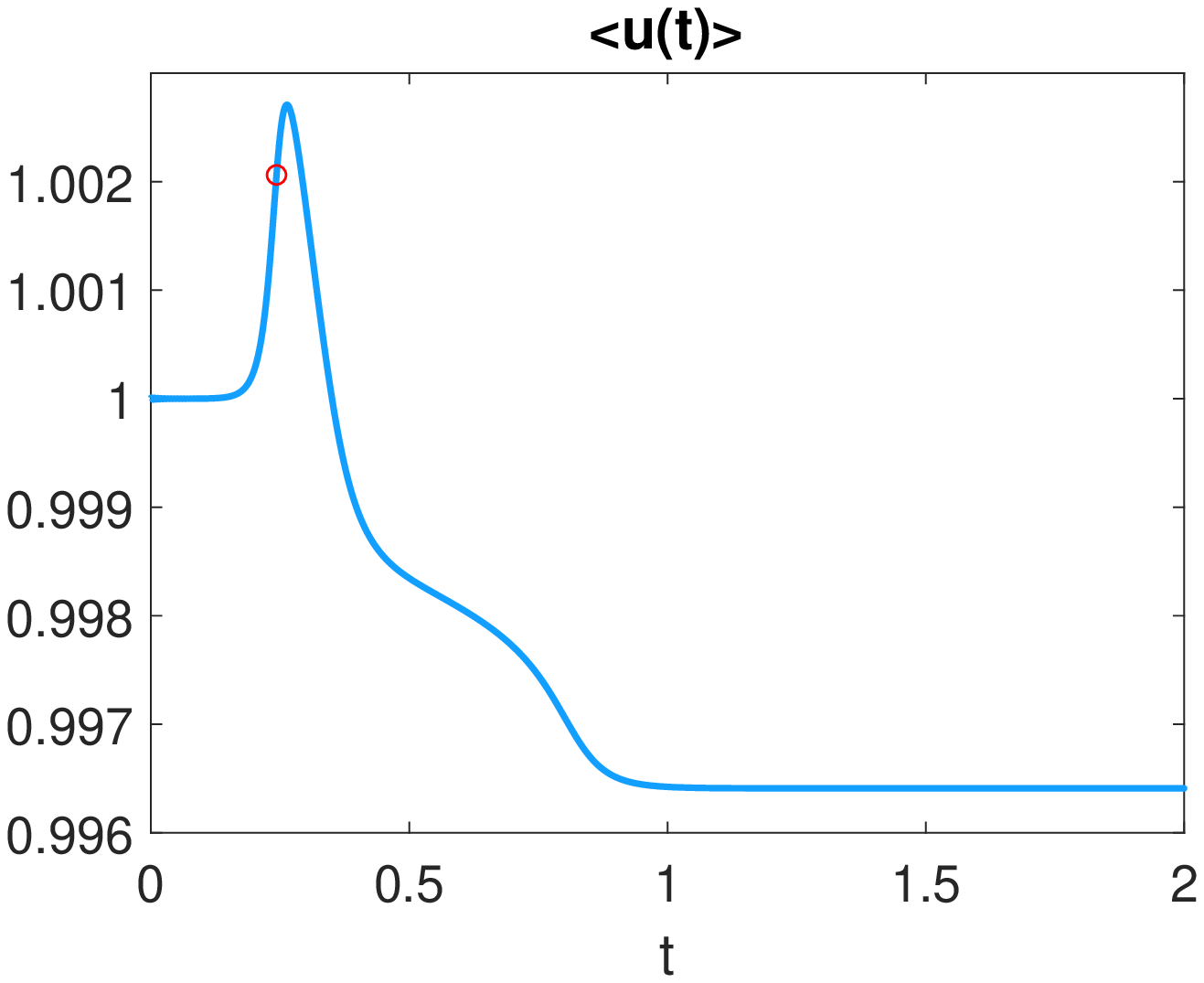}
\caption{}
\end{subfigure}
\begin{subfigure}{0.49 \textwidth}
\captionsetup{justification=centering}
\centering
\includegraphics[scale=0.45]{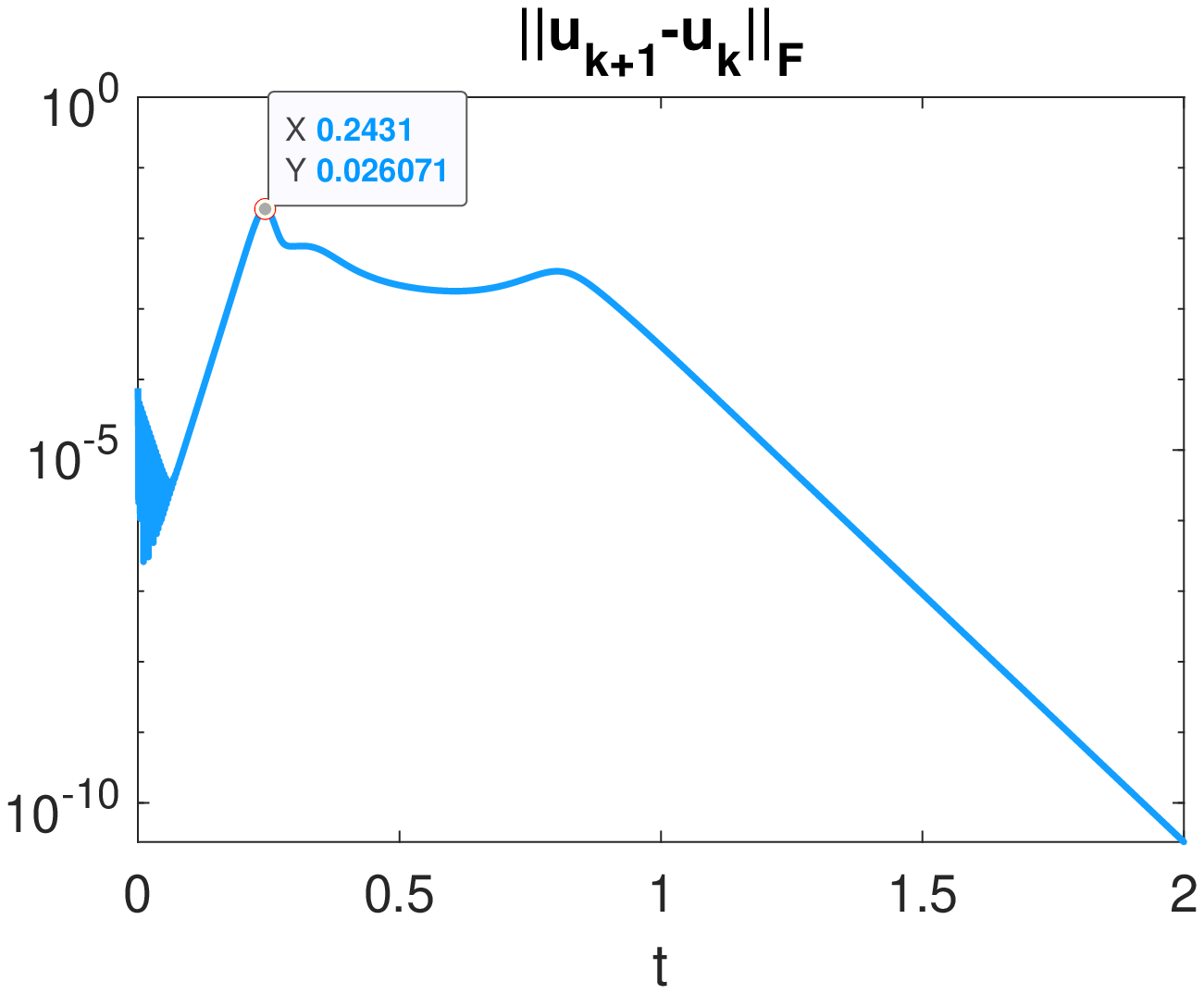}
\caption{}
\end{subfigure}
\captionsetup{justification=justified}
\caption{Test 2: Schnakenberg model. Time dynamics of the spatial mean $\langle u(t) \rangle$ (a) and of the increment $\delta_k = \|u_{k+1}-u_k \|_F$ (b) for the full model solution $u$. The time value $\tau$ to split the integration interval in $\mathcal{I}_1=[0, \tau]$ and $\mathcal{I}_2=[\tau,T]$ is indicated with a red 'o' symbol and also with a data tip in the right plot.}
\label{fig_schnak3}
\end{figure}

{\bf Adaptivity.} In the above corrected algorithms the cost to obtain $\{\tilde{ \u}_R,\tilde{ \v}_R \}$ for $R=80$ in the offline stage is $33.12$ seconds. Also for the Schnackenberg model our purpose is to reduce this cost by applying an adaptive strategy.
Hence, to identify the two subintervals $\mathcal{I}_1 = [0, \tau]$, $\mathcal{I}_2 = [\tau, T]$, from the snapshot matrix $S_\u$ we extract the spatial mean $\langle u(t) \rangle$ and the increment $\delta_k = \|\u_{k+1}-\u_k \|_F, \ k=0, \dots, n_t-1$, that are reported in Figure \ref{fig_schnak3}(a)-(b) respectively.


The maximum of the increment allows to identify $\tau = 0.2431$ to split the main interval (see red 'o' symbol). The initial conditions in $\mathcal{I}_2$ are those defined in \eqref{ic:i2} where $\u_{r_1}^{(1)}(\tau)$ and $\v_{r_1}^{(1)}(\tau)$ are the solutions of the PODc system in $\mathcal{I}_1$ with  $R_1 = 35 < \rho_{sol}$ and $r_1 = 10<R_1$.

We compare the adaptive and non adaptive correction techniques in terms of computational cost in the offline stage and in terms of the generated relative errors for increasing value $r$ of the reduced space.
In Table \ref{tab_adapt_schnak} we report the values $\ell$ used for DEIM, $R_i, i=1,2$ for the correction sizes and the CPU time in the offline stage to approximate $\{\tilde{ \u}_{R_i},\tilde{ \v}_{R_i} \}$. Also in this case, the computational load in the offline stage due to the correction is improved by the adaptive approach (the speed up factor is about 8).

\begin{table}[htb!]
\centering
\begin{tabular}{ c | c | c | c }
& $\ell$ (DEIM) & $R$ (correction) & CPU time for $\{\tilde{ \u}_R,\tilde{ \v}_R \}$\\ 
\hline
$[0,T]$ & $90$ & $80$ & $33.12$s\\
$\mathcal{I}_1$ & $50$ & $35$ & $\phantom{x}0.88$s \\
$\mathcal{I}_2$ & $14$ & $15$ & $\phantom{x}2.98$s \\
\end{tabular}
\captionsetup{justification=justified}
\caption{Test 2: Schnakenberg model. Values used to solve the reduced systems \eqref{podc_simpl} and \eqref{podeim_correction} in each time interval. The computational load in the offline stage for the correction is improved by the adaptive approach (the speed up factor is about 8).}
\label{tab_adapt_schnak}
\end{table}


\begin{figure}[t!]
\centering
\begin{subfigure}{0.49 \textwidth}
\captionsetup{justification=centering}
\centering
\includegraphics[scale=0.45]{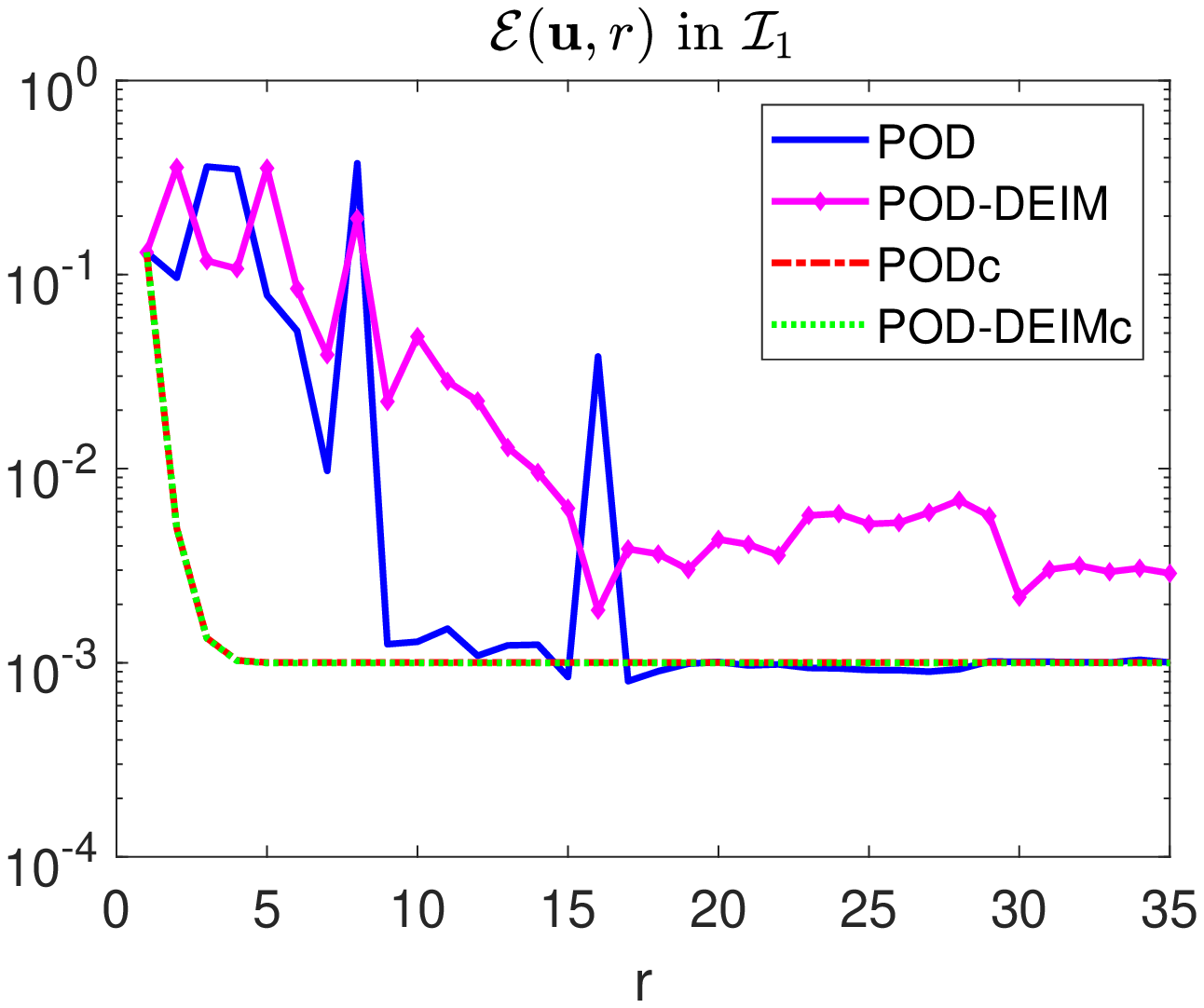}
\caption{}
\end{subfigure}
\begin{subfigure}{0.49 \textwidth}
\captionsetup{justification=centering}
\centering
\includegraphics[scale=0.45]{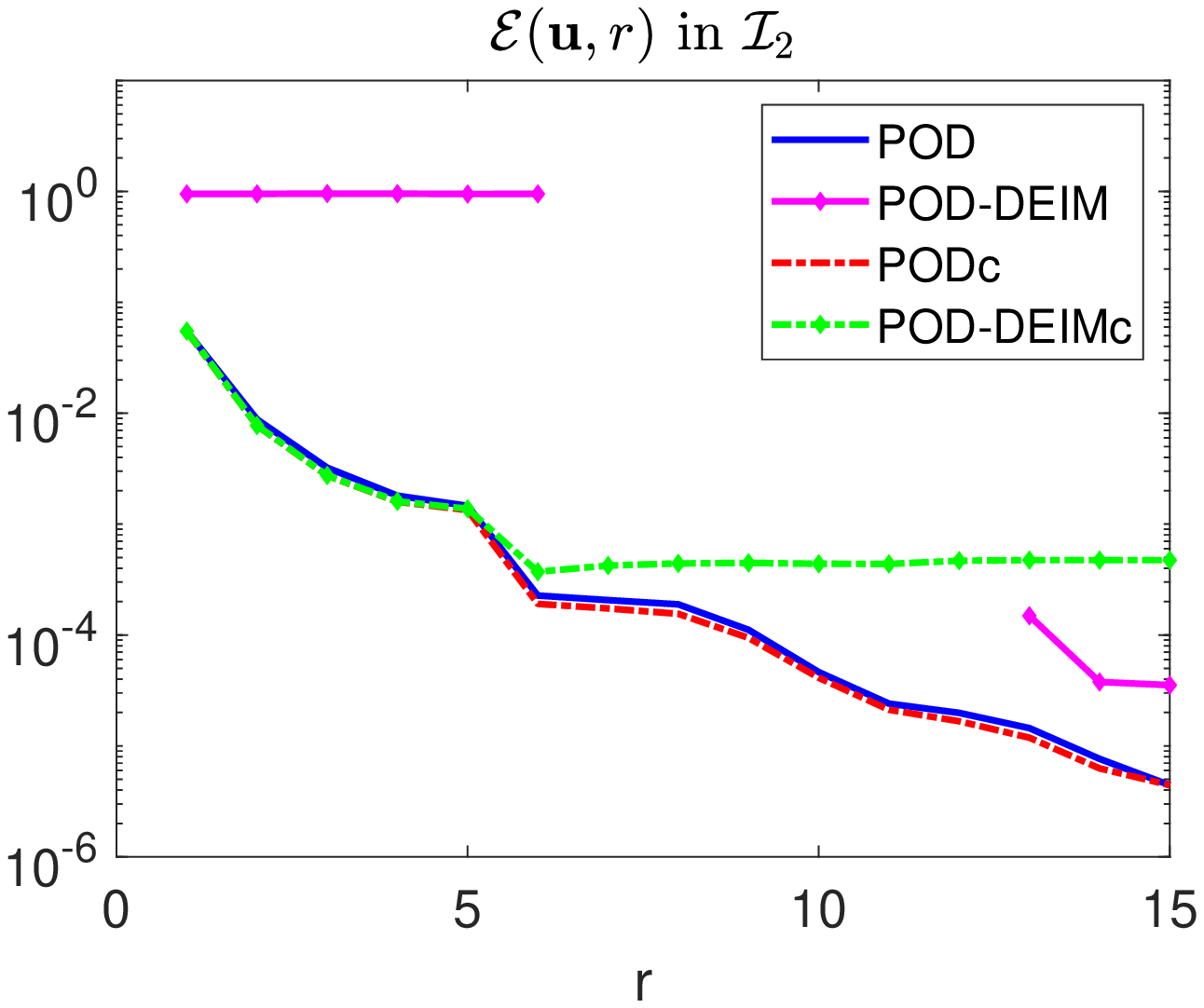}
\caption{}
\end{subfigure}
\captionsetup{justification=justified}
\caption{Test 2: Schnakenberg model adaptivity. (a) Zone $\mathcal{I}_1$: relative error $\mathcal{E}(\u,r)$ at $\tau$. (b) Zone $\mathcal{I}_2$: relative error $\mathcal{E}(\u,r)$ at the final time 
$T = 2$. In Table \ref{tab_adapt_schnak} are reported the values of $\ell, R_i, i=1,2$ used for the computation of DEIM and correction terms. In $\mathcal{I}_1$ PODc and POD-DEIMc  show the same decreasing error decay, they stabilizes the uncorrected algorithms and tends to $10^{-3}$ for $r \geq 5$. In $\mathcal{I}_2$: POD and its corrected counterpart have very similar decreasing behaviour; POD-DEIM is totally unstable, but POD-DEIMc stabilizes it showing a monotone decay that tends to $10^{-4}$ for $r \geq 6$.}
\label{fig_schnak4}
\end{figure}

The errors ${\cal{E}}(\u,r)$ on $\mathcal{I}_1$ and $\mathcal{I}_2$ are shown in Figure \ref{fig_schnak4}, (a) and (b) respectively. In $\mathcal{I}_1$: POD is unstable for $r \leq 20$;  POD-DEIM exhibits an oscillating behaviour;  their corrected counterparts have the same decreasing monotone decay and for $r \geq 5$  tend to be constant with an order of $ 10^{-3}$ (see Figure \ref{fig_schnak4}(a)). In $\mathcal{I}_2$: both POD and PODc have a similar monotone decreasing error that attains an order of $10^{-6}$ for $r = 15$; POD-DEIM is completely unstable, while the POD-DEIMc error decreases and for $r \geq 6$ tends to stagnate around $10^{-4}$ (see Figure \ref{fig_schnak4}(b)).

%

\begin{figure}[tbp!]
\centering
\begin{subfigure}{0.49 \textwidth}
\captionsetup{justification=centering}
\centering
\includegraphics[scale=0.45]{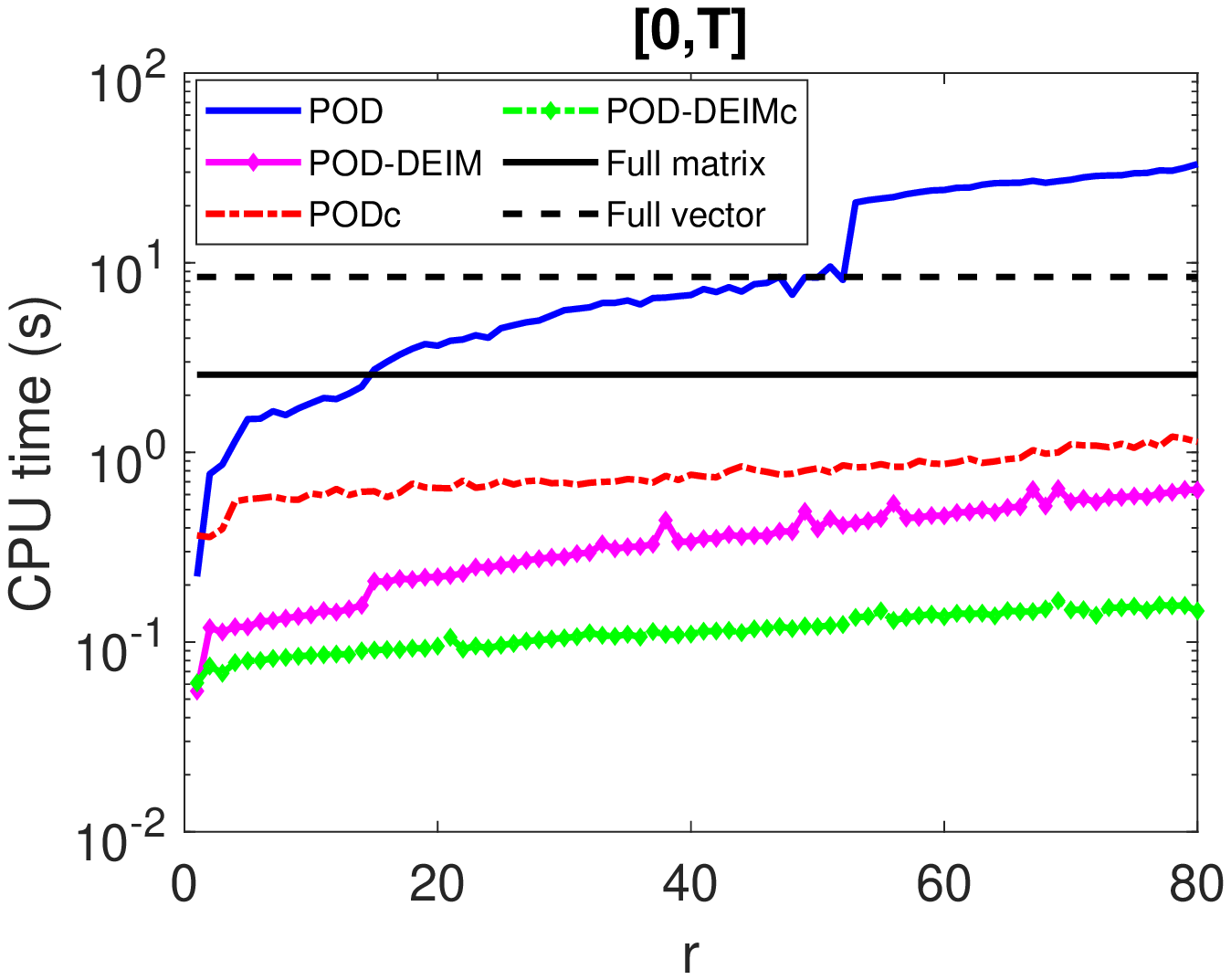}
\caption{}
\end{subfigure}
\begin{subfigure}{0.49 \textwidth}
\captionsetup{justification=centering}
\centering
\includegraphics[scale=0.45]{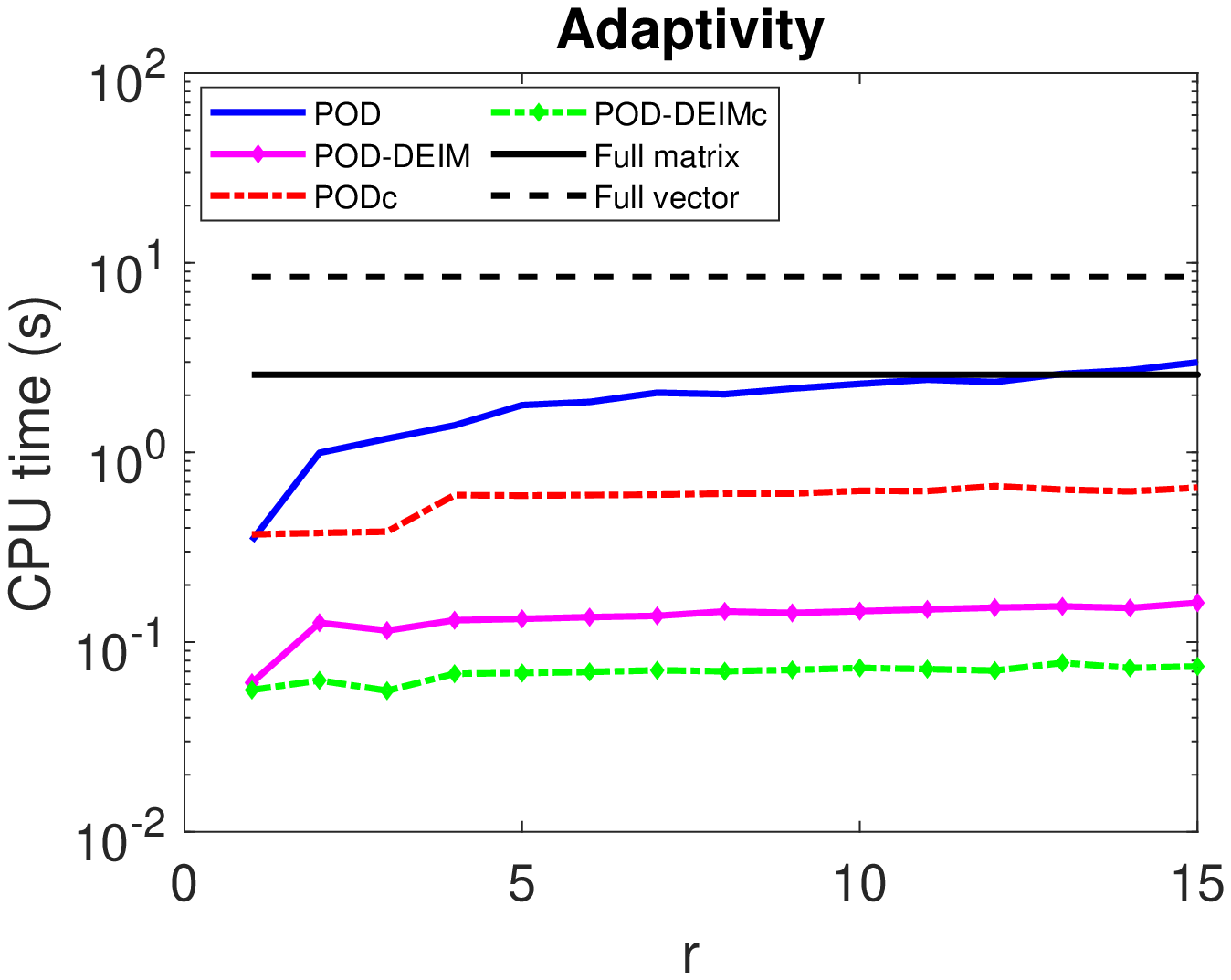}
\caption{}
\end{subfigure}
\captionsetup{justification=justified}
\caption{Test 2: Schnakenberg model. Comparison in terms of computational cost in the online stage. (a) CPU time (s) for solving the reduced models on $[0,T]$ and (b) with adaptivity on $\mathcal{I}_1$ and $\mathcal{I}_2$. (The computation of $\{\tilde{\u}_R, \tilde{\v}_R\}$ for the corrected systems is considered offline.)
Even if POD and PODc are stable with the same error behaviour, also in adaptive way (see Fig. 10) POD is too expensive. The best performance in terms of cost and errors is obtained by PODc: it is stable in both adaptive and non adaptive implementation, with sligthly lower cost in the adaptive way. Its total cost (considering also the offline contribution) improves if used in adaptive way.
We also note that POD-DEIM, although very fast, is completely unstable.}
\label{fig_schnak6}
\end{figure}

{\bf Computational cost: online stage.} As for the FHN model, we compare all the techniques proposed in terms of computational cost in the online stage with respect to the full model approximation by IMEX-Euler in matrix (black continuous line, $2.6$ seconds) and vector form (black dashed line, $8.4$ seconds). In Figure \ref{fig_schnak6}(a) we report the CPU time (seconds) for solving the reduced models in $[0,T]$ for $r \leq 80 = R$, while in the right panel the results in the case of adaptivity for $r \leq 15 = R_2 \approx \rho_{sol}$, where the CPU times to solve the two subsystems on $\mathcal{I}_1$ and $\mathcal{I}_2$ have been additioned. (Note that POD-DEIM, although very fast, is completely inaccurate, see Figures \ref{fig_schnak1}(c) and \ref{fig_schnak4}(b)). 
We find that, in all cases, the corrected methods are faster than the full model.

%
As final study, in Table \ref{tab_schnak} we report the CPU time needed to achieve a desired accuracy $tol$ by the corrected techniques. We also report for each $tol$ the values $r_0$ such that $\mathcal{E}(\u,r) \leq tol$ for $r \geq r_0$. The computational costs are very low and there is a further speed up by applying the POD-DEIMc. Nevertheless it results to be less accurate.

\begin{table}[tbp]
\centering
\begin{tabular}{ c c c | c c | c c | c c | c c  }
tol & $10^{-2}$ & $r_0$ & $10^{-3}$ & $r_0$ &  $10^{-4}$ & $r_0$ &  $10^{-5}$ & $r_0$ & $10^{-6}$ & $r_0$ \\
\hline
PODc & $0.40$s & $\phantom{x}3$	& $0.57$s & $\phantom{x}6$ & $0.61$s & $10$ & $0.58$s & $16$ & $0.65$s & $21$ \\
PODc adaptive  & $0.38$s & $\phantom{x}2$ & $0.60$s	& $\phantom{x}6$	& $0.61$s	& $\phantom{x}9$ & $0.62$s & $14$ & -  & - \\
POD-DEIMc   & $0.07$s & $\phantom{x}3$ & $0.08$s & $\phantom{x}6$ & - & - & - & - & - & - \\
POD-DEIMc adaptive   & $0.06$s & $\phantom{x}2$ & $0.07$s & $ \phantom{x}6$ & - & - & - & - & - & - \\
\end{tabular}
\captionsetup{justification=justified}
\caption{Test 2: Schnakenberg model. CPU time needed by the corrected techniques to achieve a desired accuracy $tol$ for the relative error of the unknown $u$ such that $\mathcal{E}(\u,r) \leq tol$ for $r \geq r_0$. The costs are very similar by applying or not an adaptive approach.}
\label{tab_schnak}
\end{table}

In conclusion, we can say that for the Schnackebenrg model the best performance in terms of cost and errors is obtained by PODc. It is stable in both adaptive and non adaptive implementation, can attain high accuracy for moderate values of $r$ and its total cost (considering also the offline contribution for the correction) improves if used in adaptive way.

\subsection{Test 3: DIB model}\label{run3}
As a final applicative example (see \cite{BMS21}), we consider again the DIB morpho-chemical model for electrodepostion with kinetics \eqref{DIB_kin} with parameter choice given in Sections \ref{sec:22} and \ref{sec:31}.

\begin{figure}[bt!]
\centering
\begin{subfigure}{0.49 \textwidth}
\centering
\includegraphics[scale=0.45]{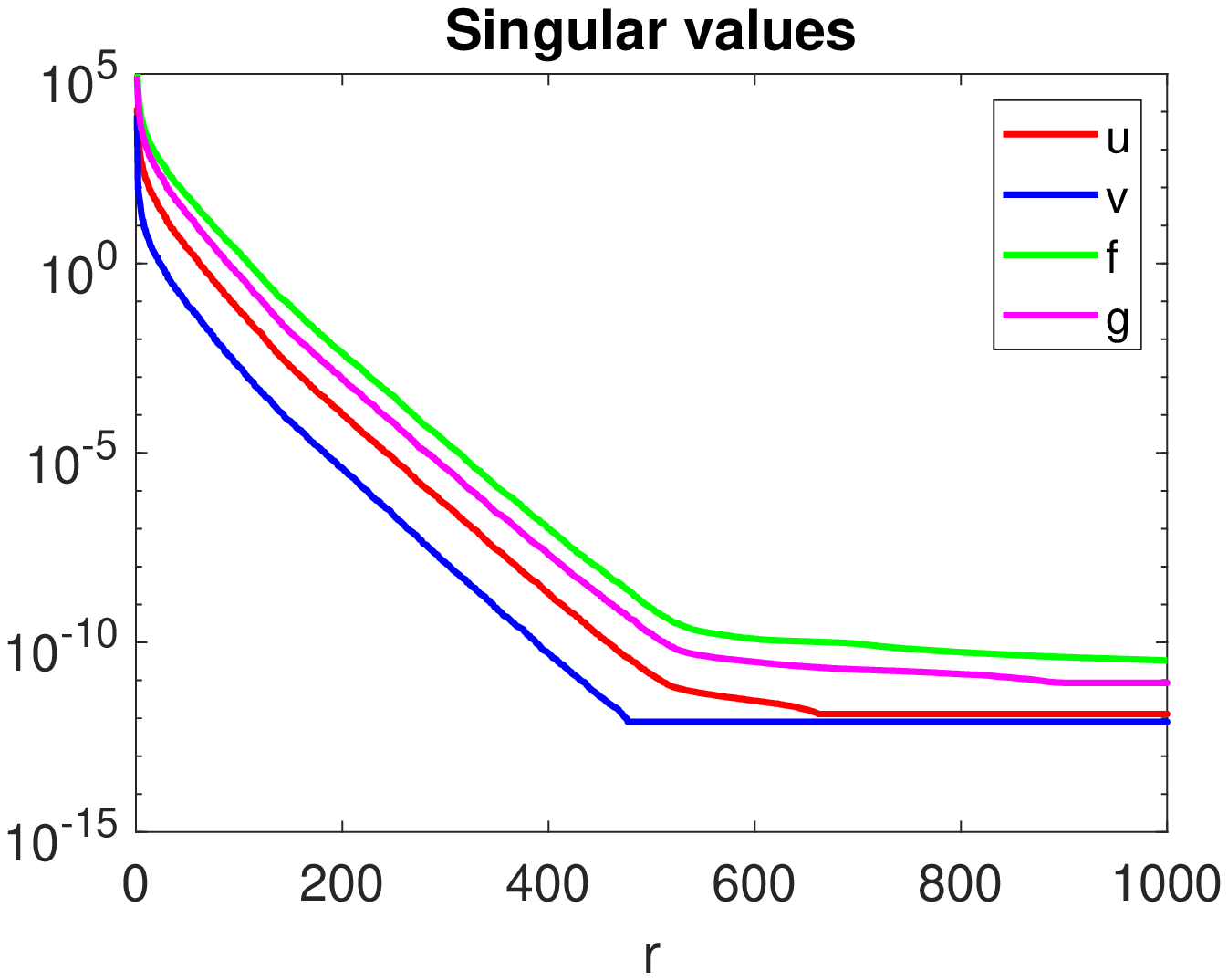}
\captionsetup{justification=centering}
\caption{}
\end{subfigure}
\begin{subfigure}{0.49 \textwidth}
\includegraphics[scale=0.45]{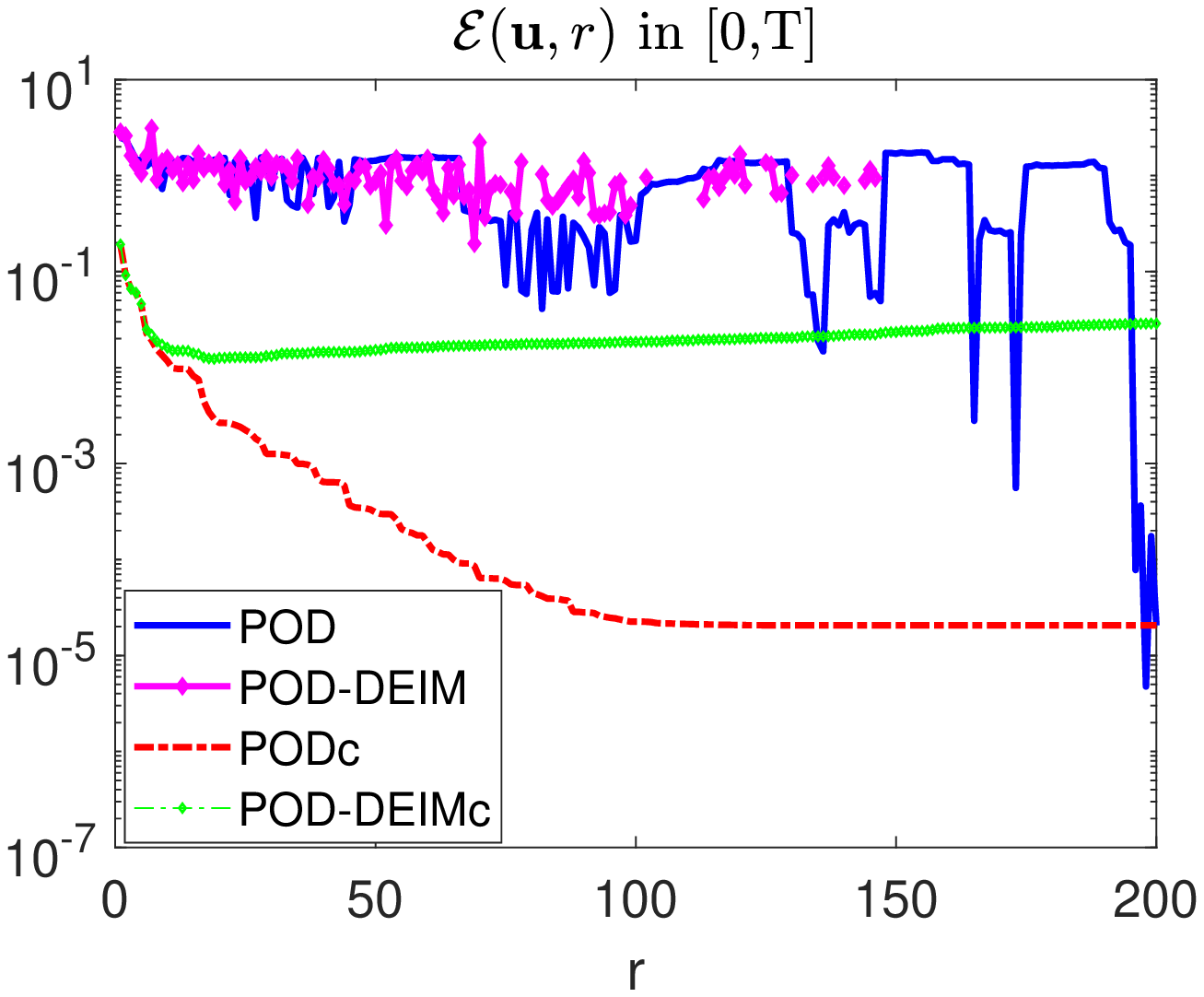}
\captionsetup{justification=centering}
\caption{}
\end{subfigure}
\captionsetup{justification=justified}
\caption{Test 3: DIB model. (a) Singular values decay of the snapshot matrices defined in \eqref{snap:mat} and \eqref{non_snap}. (b) Relative error $\mathcal{E}(\u,r)$ at the final time $T = 100$. The correction algorithms are applied by choosing $R = 200$, the DEIM techniques are applied with $\ell=363$. The standard POD technique shows a very irregular error behaviour; POD-DEIM is completely unstable. The instability is overcome by using the corrected models. As for the other models, POD-DEIMc exhibits for $r \geq 10$ a stagnation of the error around $10^{-2}$.}
\label{fig_dib1}
\end{figure}

The snapshot matrices \eqref{snap:mat} are $S_\u$, $S_\v \in \RR^{10000\times 25001}$. In Figure \ref{fig_dib1}(a), we report the singular values decay for both snapshot matrices \eqref{snap:mat} and \eqref{non_snap} for $r=1, \dots, 1000$. The decay is very slow with respect to the RD-PDE models studied in the previous cases and the ranks of the snapshot matrices are $\rho_{sol} = 342$ and $\rho_{kin} = 363$. We argue that these high values can depend from the more rich labyrinth structure of the Turing pattern expected for this parameter choice (see Figure \ref{fig_dib}(a)) .

{\bf Stabilization.} First of all we begin by solving each reduced model \eqref{pod_galerkin} and \eqref{deim_2sottospazi} in the time interval $[0,T]$ with time step $h_t = 10^{-3}$. We remind that DEIM is applied with $\ell = 363 = \rho_{kin}$ (see Section \ref{sec:31}). The corrected systems \eqref{podc_simpl} and \eqref{podeim_correction} are integrated by choosing $R = 200 < \rho_{sol}$. In Figure \ref{fig_dib1}(b) the relative errors for increasing values of $r \leq r_{max}=R = 200$ for POD and POD-DEIM are those reported in Figure \ref{fig:motex2} for the unknown $u$.
As already discussed, POD exhibits a highly erratic and somehow oscillating error behaviour for several values of $r$, whereas POD-DEIM results to be completely unstable (no line is shown).

From Figure \ref{fig_dib1}, it is evident that a stabilization occurs by using PODc and POD-DEIMc that avoid the previous drawbacks and present monotone error decays. PODc (for $r \geq 100$) tend to stagnate around $10^{-5}$, similarly POD-DEIMc around $10^{-2}$ but earlier for $r \geq 10$. We will show that the bad trend of POD-DEIMc will be improved by using the adaptive strategy. 


{\bf Adaptivity.} We show that also for the DIB model the adaptivity reduces the computational cost needed by correction in the offline stage. This cost on the whole time interval is about $2240$ seconds.
The goal is again to adapt the value of the correction and DEIM parameters $R$ and $\ell$ to the time regimes of the Turing dynamics.
The spatial mean $\langle u(t) \rangle$ and the increment of the numerical solution $u$ are reported in Figure \ref{fig_dib}(b)-(c), where the maximum of the increment at $\tau = 2.123$ allows to identify the two subintervals $\mathcal{I}_1 = [0, \tau]$ and $\mathcal{I}_2 = [\tau,T]$ where the adaptive MOR strategies can be applied. The snapshot matrices in $\mathcal{I}_1$ are constructed by taking the full model solutions $\u_k$, $\v_k$ every two time steps. This is motivated by the very small size of the reactivity zone and by the need to give enough information in this time interval. The initial conditions in $\mathcal{I}_2$ are those defined in \eqref{ic:i2} where $\u_{r_1}^{(1)}(\tau)$ and $\v_{r_1}^{(1)}(\tau)$ are the solutions of the PODc system in $\mathcal{I}_1$ with $R_1 = 41 \approx \rho_{sol}$ and $r_1 = 10 \ll R_1$.

Also for the DIB model, we compare the correction techniques in adaptive and non adaptive implementation in terms of computational cost in the offline stage and in terms of the relative errors $\mathcal{E}(\u,r)$.
In Table \ref{tab_adapt_dib} we list the choices of $\ell$ for DEIM, of $R_i, i=1,2$ for the correction spaces and the CPU time in the offline stage to approximate $\{\tilde{ \u}_{R_i},\tilde{ \v}_{R_i} \}$. 
We can deduce that also for the DIB model the cost in the offline stage due to the correction is reduced by the adaptive approach, in fact there is a speed up factor of 1.5.
\begin{table}[tbp]
\centering
\begin{tabular}{ c | c | c | c }
& $\ell$ (DEIM) & $R$ (correction) & CPU time for $\{\tilde{ \u}_R,\tilde{ \v}_R \}$\\
\hline
$[0,T]$ & $363$ & $200$ & $2239.9$s\\
$\mathcal{I}_1$ & $\phantom{x}60$ & $\phantom{x}41$ & $\phantom{xxx}5.8$s\\
$\mathcal{I}_2$ & $324$ & $150$ & $1523.4$s\\
\end{tabular}
\captionsetup{justification=justified}
\caption{Test 3: DIB model. MOR parameters used to solve the reduced systems \eqref{podc_simpl} and \eqref{podeim_correction} in each time interval. The advantage of the adaptive strategy becomes evident by comparing the computational times in the last column.}
\label{tab_adapt_dib}
\end{table}

\begin{figure}[t!]
\centering
\begin{subfigure}{0.49 \textwidth}
\centering
\includegraphics[scale=0.45]{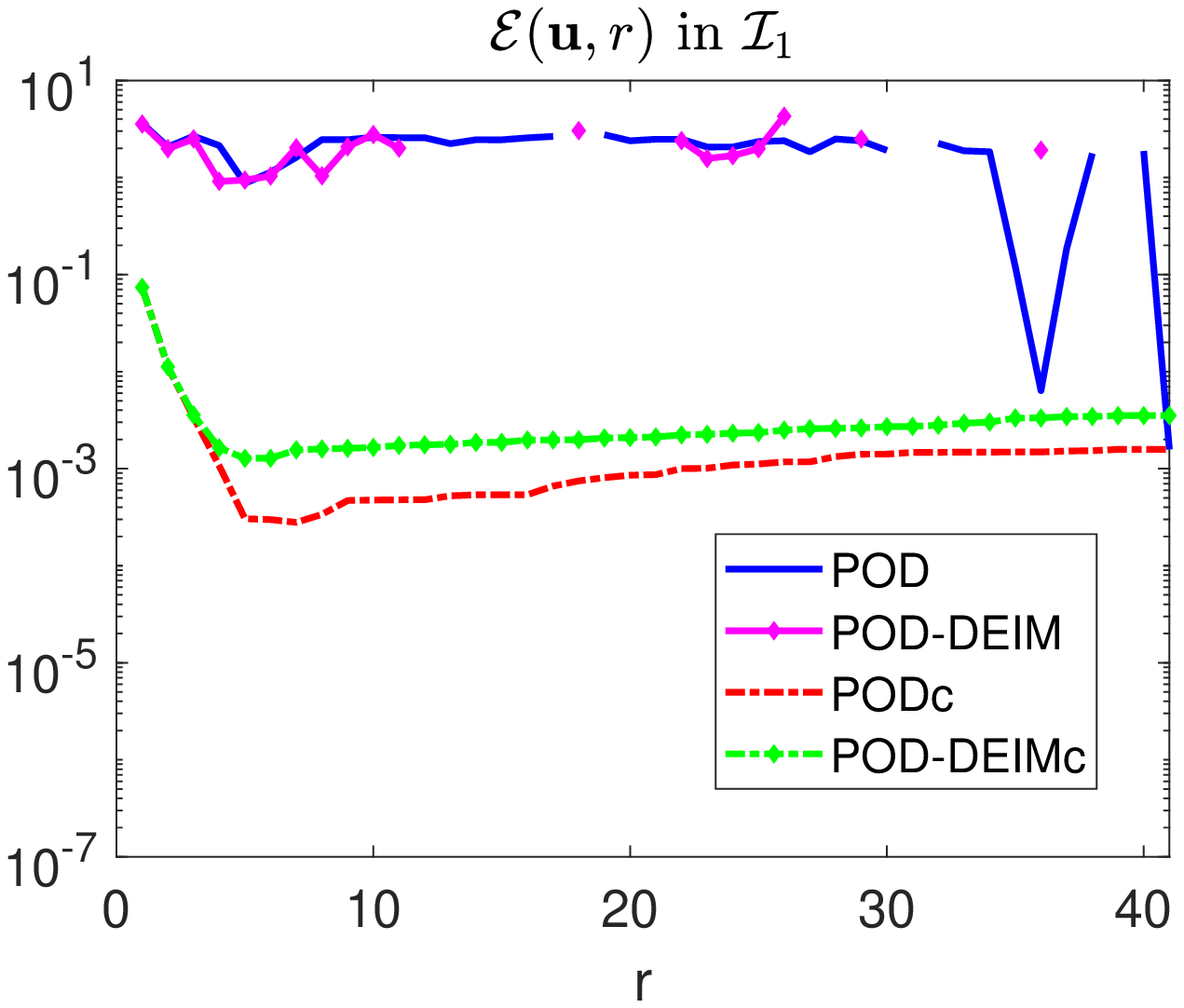}
\captionsetup{justification=centering}
\caption{}
\end{subfigure}
\begin{subfigure}{0.49 \textwidth}
\centering
\includegraphics[scale=0.45]{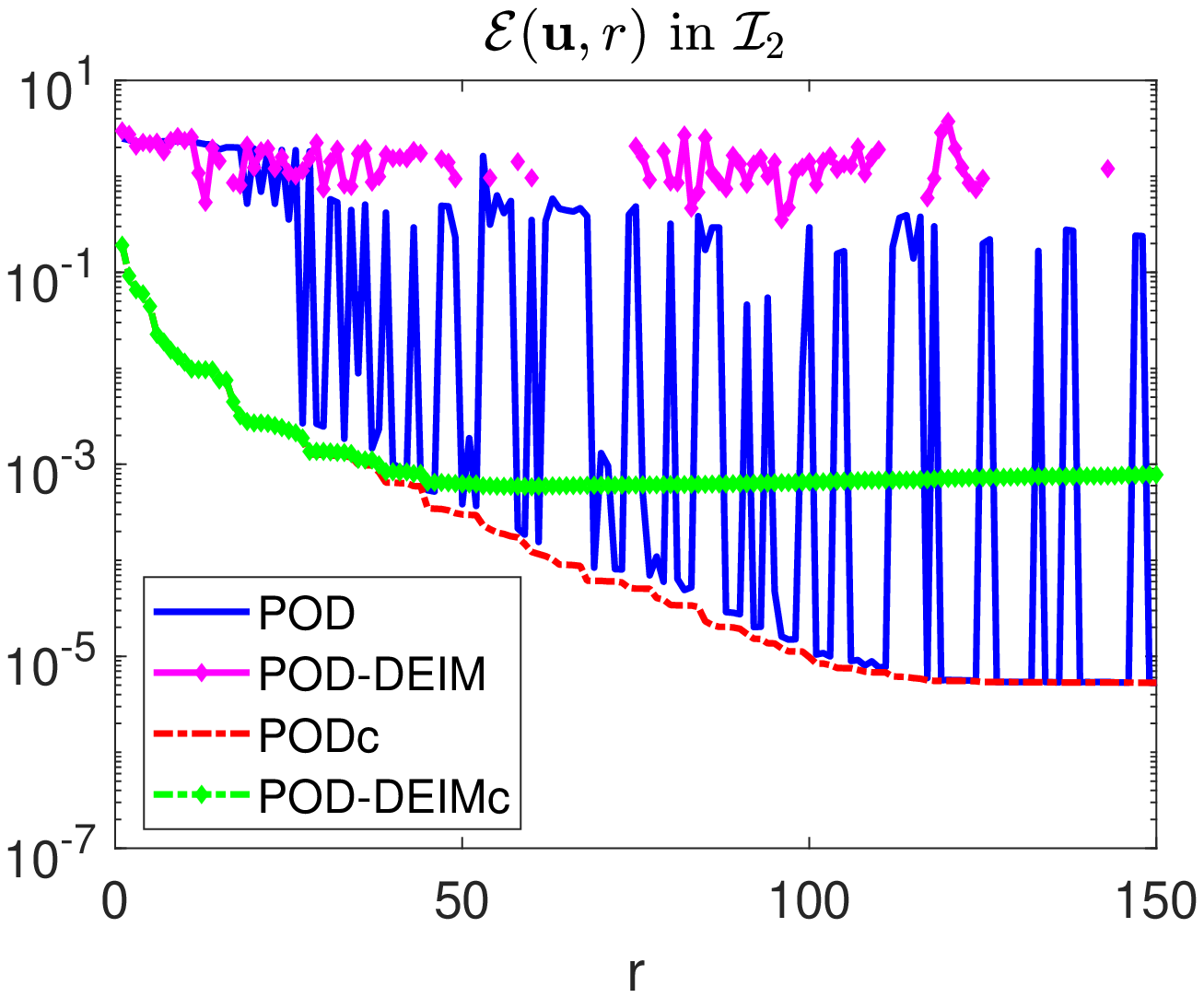}
\captionsetup{justification=centering}
\caption{}
\end{subfigure}
\captionsetup{justification=justified}
\caption{Test 3: DIB model adaptive algorithms. (a) Zone $\mathcal{I}_1$: relative error $\mathcal{E}(\u,r)$ at the final time $\tau$ for $u$. (b) Zone $\mathcal{I}_2$: relative error at the final time $T = 100$ for the unknown $u$. MOR parameter values for the application of POD-DEIM and correction are reported in Table \ref{tab_adapt_dib}. The classical POD and POD-DEIM are very unstable in both zones. PODc and POD-DEIMc stabilize this trend in $\mathcal{I}_1$ both with low accuracy, then in $\mathcal{I}_2$ they exhibit a decreasing error decay, even if POD-DEIMc for $r \geq 40$ tends to a constant error of order $10^{-3}.$}
\label{fig_dib4}
\end{figure}

The errors $\mathcal{E}(\u,r)$ by all algorithms are shown in Figure \ref{fig_dib4}(a) for the simulations on $\mathcal{I}_1$ and in Figure \ref{fig_dib4}(b) for $\mathcal{I}_2$.
The classical POD and POD-DEIM are unstable or completely erratic (POD) in both zones. The application of PODc and POD-DEIMc stabilizes these bad behaviours in both subdomains. In particular, in $\mathcal{I}_1$ the errors of PODc and POD-DEIMc are slightly increasing for $r \geq 8$ but bounded around $10^{-3}$. In $\mathcal{I}_2$ the POD-DEIMc has a similar trend with a constant error of order $10^{-3}$ for $r \geq 40$. Instead, PODc exhibits a monotone decay and reaches an error of order $10^{-6}$ for $r\geq 100$.
It is worth noting that in the first region $\mathcal{I}_1$ we do not pretend a very low error, because there the main goal is to obtain $\u_{r_1}^{(1)}(\tau)$ and $\v_{r_1}^{(1)}(\tau)$ as ``acceptable'' initial conditions for the subsystem in $\mathcal{I}_2$.

%

\begin{figure}[tbp]
\centering
\begin{subfigure}{0.49 \textwidth}
\centering
\includegraphics[scale=0.48]{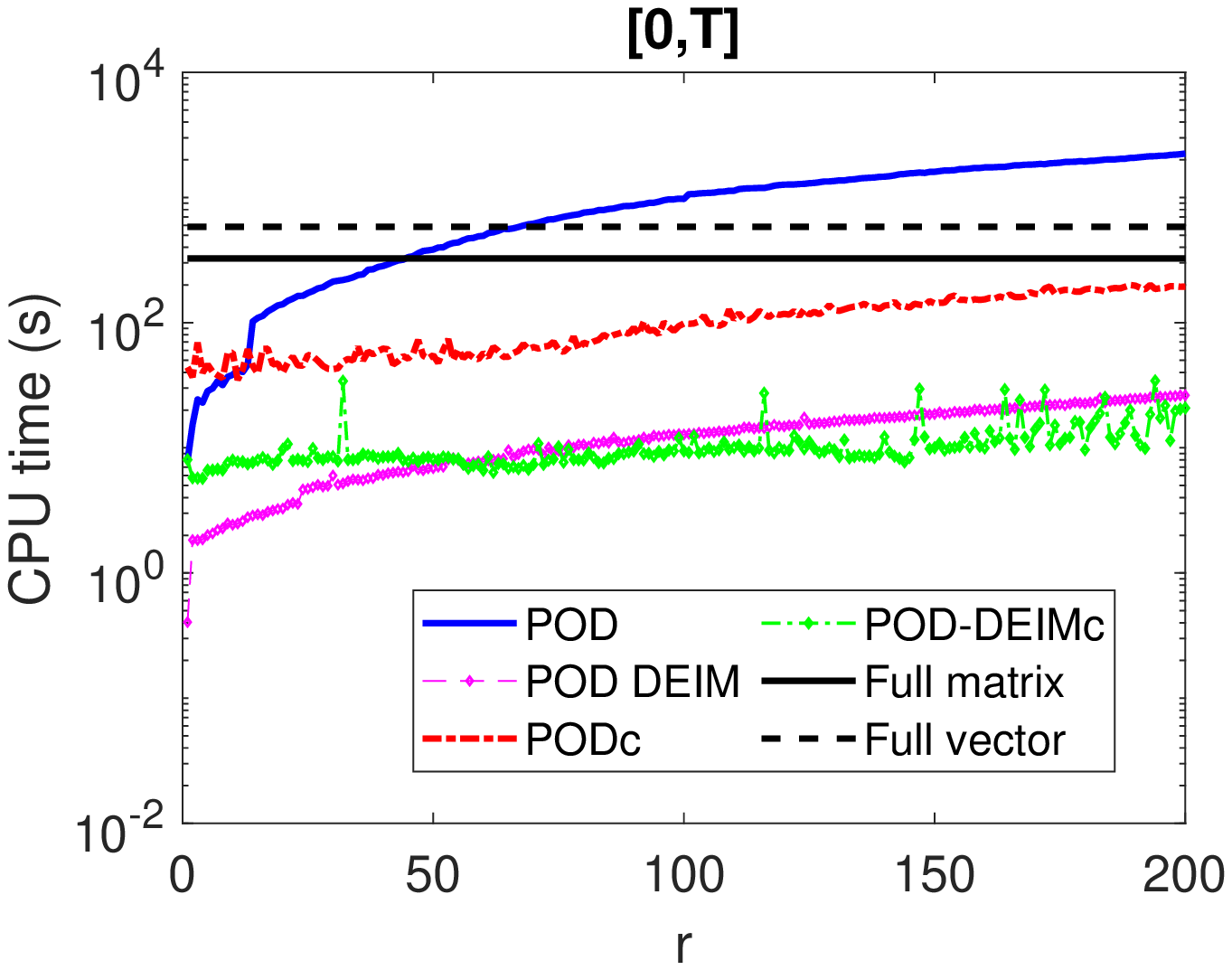}
\captionsetup{justification=centering}
\caption{}
\end{subfigure}
\begin{subfigure}{0.49 \textwidth}
\centering
\includegraphics[scale=0.48]{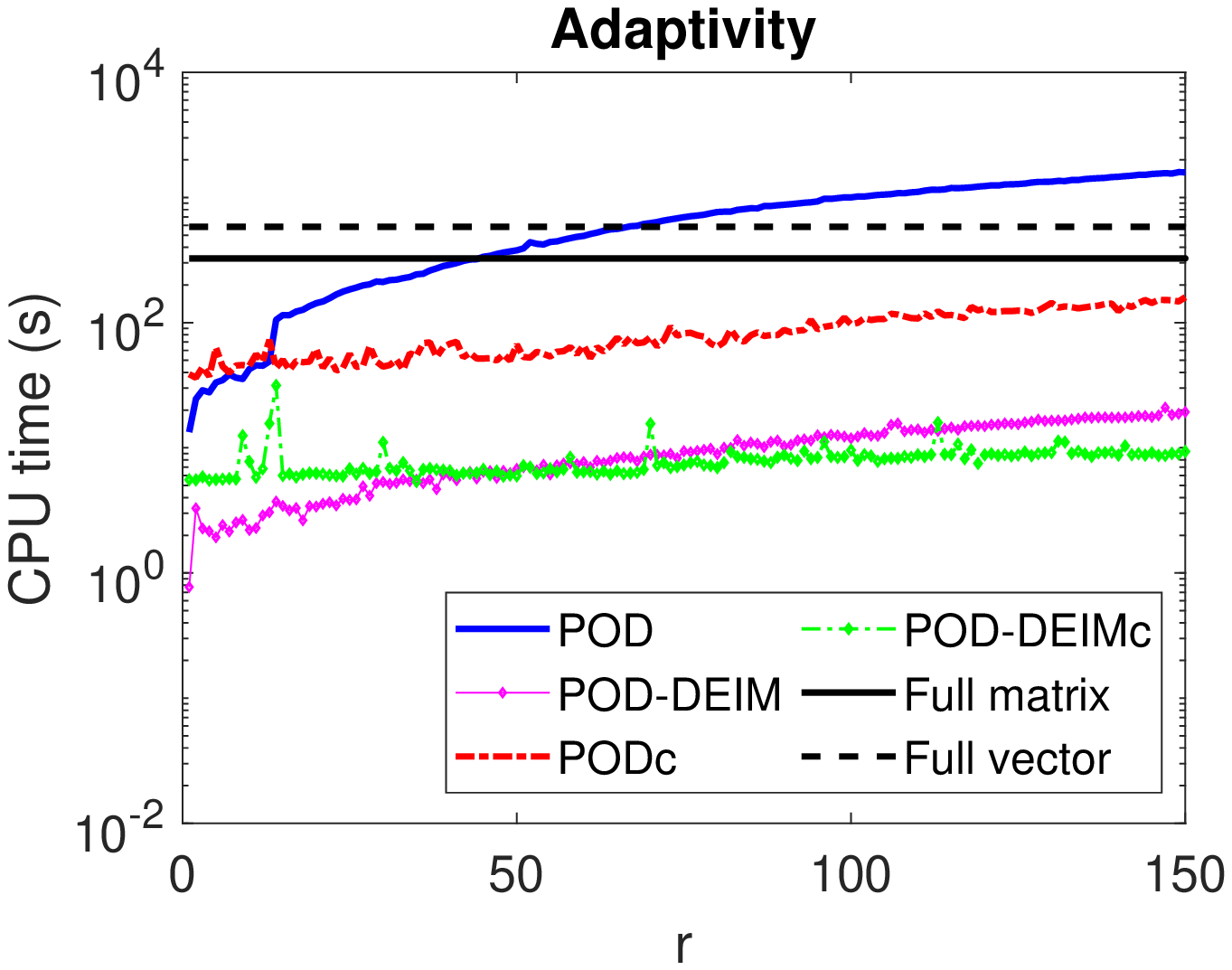}
\captionsetup{justification=centering}
\caption{}
\end{subfigure}
\captionsetup{justification=justified}
\caption{Test 3: DIB model. Online computational costs in CPU time (s). (a) Solving the reduced models in the time interval $[0,T]$; (b) same algorithms with adaptivity. The computation of $\{\tilde{\u}_R, \tilde{\v}_R\}$ for the corrected systems is considered offline and is reported in Table \ref{tab_adapt_dib}.}
\label{fig_dib6}
\end{figure}

{\bf Computational cost: online stage.} To decide which is the best MOR approach among those we proposed, we compare all techniques in terms of computational costs. It is worth noting that DIB is more demanding than the other RD-PDE models considered due to the well structured labyrinth pattern solution and the nonlinearity of the kinetics. In fact, to solve the full model the IMEX-Euler scheme in matrix form employed $326.3$ seconds and $582.7$ seconds in vector form (see Section \ref{sec:2}) (see the black continuous and dashed lines reported in Figure \ref{fig_dib6}). The computational costs of the corrected and classical POD and POD-DEIM algorithms are also reported in Figure \ref{fig_dib6}. The left panel concerns the simulations on the entire interval $[0,T]$ for $r \leq r_{max}=R = 200$, while the right one those for the adaptive strategy on both $\mathcal{I}_1$ and $\mathcal{I}_2$ until $r \leq r_{max}= R_2 = 150$.
It is easy to see that the PODc and POD-DEIMc result to be faster than the full model for all choices of $r$. We also note that although POD-DEIM is faster than the full model, it is very inaccurate and unstable (see Figures \ref{fig_dib1}(b) and \ref{fig_dib4}). 

In Table \ref{tab_dib} we compare the performances of the corrected methods with and without adaptivity in terms of the CPU time needed to reach a desired accuracy $tol$. We also report for each $tol$ the value $r_0$ such that $\mathcal{E}(\u,r) \leq tol$ for $r \geq r_0$. By using the adaptive approach we can achieve a better accuracy up to two orders of magnitude (for the POD-DEIMc), even though we have almost the same cost by applying or not an adaptive approach. As reported in Table \ref{tab_adapt_dib}, the adaptivity shows its advantages in the offline stage.
\begin{table}[btp]
\centering
\begin{tabular}{c c c | c c | c c | c c | c c }
 tol & $10^{-1}$ & $r_0$ & $10^{-2}$ & $r_0$ & $10^{-3}$ & $r_0$ & $10^{-4}$& $r_0$ & $10^{-5}$ & $r_0$ \\
\hline
PODc   & $36.74$s & $\phantom{x}2$ & $34.87$s & $11$ & $64.10$s & $35$ & $53.48$s & $65$ & - & - \\
PODc adaptive  & $36.32$s & $\phantom{x}2$ & $55.78$s & $11$ & $54.34$s & $35$ & $65.30$s & $64$ & $97.06$s & $100$ \\
POD-DEIMc  & $\phantom{x}5.76$s & $\phantom{x}2$ & - & - & - & - & - & - & - & -\\
POD-DEIMc adaptive  & $\phantom{x}5.54$s & $\phantom{x}2$ & $\phantom{x}6.80$s & $12$ & $\phantom{x}6.69$s & $38$ & - & - & - & - \\
\end{tabular}
\captionsetup{justification=justified}
\caption{Test 3: DIB model. CPU time needed to achieve a desired accuracy $tol$ for the corrected systems. The best gain is obtained by PODc adaptive for $r = 100$.}
\label{tab_dib}
\end{table}

In terms of speed up factor we want to emphasize that PODc, for $tol \leq 10^{-4}$ is about 8 times faster than the full problem in the matrix form. Even though POD-DEIMc is more economic, it is less accurate achieving an error of order $10^{-3}$ with the adaptive algorithm. To conclude, the best performance in terms of accuracy and efficiency is obtained by the adaptive PODc.


\section{Conclusions and future works}

In this paper, we have presented a new algorithm that stabilizes the well-known POD-DEIM algorithm and we have applied the new approach to coupled PDE systems with Turing type solutions. The idea is to add a correction term based on high-ranked POD solution and to introduce an adaptive version based on the time dynamics of the RD-PDEs that further improves the computational efficiency. We have found that both PODc and adaptive PODc improve the accuracy of the surrogate model with respect to the classical POD. Furthermore, the adaptive POD-DEIMc allows to obtain faster computations but less accurate approximations.

The present research has been motivated by the initial results proposed in \cite{BMS21} for Turing pattern approximation. In particular, for the DIB electrochemical model and for a limited range of parameters, POD has been applied to recover patterns of different morphology by using the same POD bases and a fixed reasonable $r$. The stabilization and computational efficiency shown in Section \ref{sec:5}, encourage the application of our corrected and adaptive MOR techniques in this direction. For example, some applications in this sense aims to: (i) obtain a many-query scenario of the possible patterns in the Turing region as in \cite{SLB19}; (ii) devise smart parameter identification techniques to compare numerical solutions of the DIB model with experimental data of battery life decay; (iii) construct in economic way model-based training sets in machine learning techniques for energetic applications.

\medskip
{\bf Acknowledgemens.} AA, AM and IS are members of the InDAM-GNCS activity group. IS acknowledges
the PRIN 2017 research Project (No. 2017KL4EF3) “Mathematics of active materials: from mechanobiology to smart devices.” 



\bibliographystyle{plain}
\bibliography{refs_ALLA.bib}
%
%
%
%
%
%
%
%
%
%
%
%
%
%
%
\end{document}